\documentclass[12pt,oneside,a4paper]{report}
\usepackage{latexsym} 
\usepackage{amssymb, amscd} 
\usepackage{makeidx}

\usepackage{amsmath}
\usepackage{amstext}

\newcounter{satz}

\newenvironment{thm}[2][\itshape]{\medskip\refstepcounter{satz}{\bf(\arabic{chapter}.\arabic{section}.\arabic{satz})~#2:\ }#1}{\vspace{0cm}}
\newcounter{appsatz}

\newenvironment{proof}{{\bf Proof:\ }}{\qed}
\newcommand{\qed}{{\ifhmode\unskip\nobreak\hfil\penalty50
\hskip1em \else\nobreak\fi \nobreak\mbox{}\hfil\mbox{$\square$}
\parfillskip=0pt \finalhyphendemerits=0 \par\vskip5pt}}

\DeclareMathOperator{\depth}{depth}
 \DeclareMathOperator{\pd}{pd}
\DeclareMathOperator{\coker}{coker}
\newcommand{\nzlink}{E(3p_1\dots,3p_\alpha,3q_1+2,\dots,3q_\beta+2)}
\newcommand{\linkad}{{\rm{Link}}_{\AD}}
\newcommand{\Link}{{\rm{Link}}}
\renewcommand{\char}{\rm{char }}
\newcommand{\Tor}{\rm{Tor}}
\newcommand{\Ass}{\rm{Ass}}
\newcommand{\dm}{\rm{dim}}
\newcommand{\Zdv}{\rm{Zdv}}
\newcommand{\e}{\varepsilon}
\newcommand{\lc}{ \left\{ }
\newcommand{\rc}{ \right\} }
\renewcommand{\b}{\beta}
\renewcommand{\a}{\alpha}
\newcommand{\D}{\Delta}
\newcommand{\G}{\Gamma}
\newcommand{\AD}{\Delta^*}
\newcommand{\Hred}{\widetilde{H}}

\newcommand{\Def}{Definition}
\newcommand{\Lem}{Lemma}
\newcommand{\Prop}{Proposition}
\newcommand{\Cor}{Corollary}

\newcommand{\ie}{\mbox{i.e.,~}}

\newfont{\Fettcal}{cmbsy10 scaled 1200}
\newfont{\Mzeichen}{bbm12}
\newfont{\MLzeichen}{bbm17}                
\newfont{\MLBzeichen}{bbmbx10 scaled 1700} 
\newfont{\MSzeichen}{bbm8}                 
\newfont{\MFzeichen}{bbm12 scaled 850}                 
\newfont{\MIzeichen}{bbm7 scaled 1150}                 
\newfont{\MTzeichen}{bbm5}                 
\newfont{\MathTzeichen}{cmsy5}             
\newfont{\MathSSzeichen}{cmsy6}             

\newfont{\fraktur}{eufm10 scaled 1200}     

\newcommand{\FI}{{\mbox {\fraktur F}}}
\newcommand{\GI}{{\mbox {\fraktur G}}}

\newfont{\frakturI}{eufm10 scaled 750}     

\usepackage{graphicx}

\usepackage {texdraw}
\usepackage{epsfig}
\input txdtools
\let\et=\etexdraw
\def\etexdraw{\drawbb\et}

\begin{document}
\pagenumbering{roman}
 \setcounter{page}{2}

\begin{center}
{\huge{\textbf{Betti Numbers of Graph}}}

{\huge{\textbf{Ideals}}}

{\large{\textbf{Sean Jacques}}}

{\large {Thesis submitted to the University of Sheffield for the
degree of Doctor of Philosophy

Department of Pure Mathematics

June 2004}}
\end{center}
\newpage

I would like to thank my supervisor Moty Katzman for his support
and guidance during the preparation of this thesis. I would also
like to thank the Engineering and Physical Sciences Research
Council for supporting this work.
 \newpage
{\huge{\textbf{Abstract}}}

In this thesis we investigate certain types of monomial ideals of
polynomial rings over fields. We are interested in minimal free
resolutions of the ideals (or equivalently the quotients of the
polynomial ring by the ideals) considered as modules over the
polynomial ring. There is no simple method of finding such
resolutions but there are formulae for finding the Betti numbers
(part of the information which comprises a minimal free
resolution) of certain types of monomial ideals.

Even with these formulae it is not in general possible to find
especially explicit or useful descriptions of the Betti numbers.
However we restrict our attention to those ideals which are
generated by square free monomials of degree 2. The purpose of
this is to associate these ideals with graphs. This provides a
link between algebraic objects, the monomial modules, and
combinatorial objects, the graphs. From this correspondence we
define as new numerical invariants of graphs Betti numbers and
projective dimension. This allows us to find results about the
Betti numbers and projective dimensions of monomial modules by
using combinatorial properties of graphs.

We find explicit descriptions of the Betti numbers of particular
families of monomial modules by using this link to such families
of graphs as complete, complete bipartite graphs and cycles and
the combinatorics of these graphs. We also draw conclusions about
the projective dimensions of certain graph ideals (such as lower
bounds) using this correspondence. We also look briefly at an
alternative approach to finding free resolutions using cellular
resolutions.

In the case of tree ideals we find a method of describing the
Betti numbers in terms of the Betti numbers of ideals associated
to subtrees. This also leads to a description of the projective
dimension of a tree ideal in terms of the projective dimensions of
these subtree ideals.

  \tableofcontents
\newpage

\newpage

\pagenumbering{arabic}
\setcounter{page}{1}

\chapter{Introduction}
\setcounter{satz}{0}
\section{Introduction}
In this thesis we will study graphs in association with algebraic
objects. For any (finite, simple) graph we will define a module
over a polynomial ring. The intention is to find information about
combinatorial objects, \ie graphs, by studying the corresponding
algebraic objects and vice versa.

We do this by considering minimal free resolutions of modules. In
general finding such resolutions may not be practical. However it
will be useful to examine segments of the information which
constitutes minimal free resolutions, namely Betti numbers and
projective dimensions.

In the first chapter we will describe precisely the objects we
will be using and some of the machinery we will be requiring. In
the second chapter we investigate the projective dimension of a
graph. This is a numerical invariant of a graph defined via the
associated module. We then go on to show that the Betti numbers of
certain families of graphs, complete, complete bipartite and
cycles can be computed explicitly. Finally we examine the Betti
numbers of trees and find a recursive formula for their
calculation.

\section{Graphs and Graph Ideals}

\begin{thm}[\rm]{\Def}
A graph $G$ is a collection of vertices, which we will often write
as $x_1,\dots,x_n$ (or sometimes for simplicity as just
$1,\dots,n$) and say $G$ has vertex set $V(G)= \{x_1,\dots,x_n\}$
(or $\{1,\dots,n\}$) along with a set of edges $E(G) \subseteq
V(G) \times V(G)$. If $\{x_i,x_j\} \in E(G)$ we will say $x_i$ and
$x_j$ are joined by an edge. Graphs will be represented
pictorially in the obvious way. For example, let $G$ be the graph
with vertices $x_1,x_2,x_3,x_4$ and edges
$\{x_1,x_2\},\{x_1,x_4\},\{x_2,x_3\},\{x_3,x_4\},\{x_2,x_4\}$.

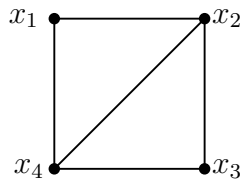
\begin{figure}[h]\label{fig1}
\begin{texdraw}
  \drawdim cm
  \textref h:R v:C
  \move(0 0)
  \move(6 0) \fcir f:0 r:0.075  \htext{$x_4 \ $}
  \lvec(8 2) \lvec(6 2) \lvec(6 0) \lvec (8 0)\lvec(8 2)
  \move(8 0) \fcir f:0 r:0.075  \move(8.5 0) \htext{$x_3$}
  \move(8 2) \fcir f:0 r:0.075  \move(8.5 2) \htext{$x_2$}
  \move(6 2) \fcir f:0 r:0.075  \move(5.8 2) \htext{$x_1$}
\linewd 0.01 \lpatt(0.2 0.2)
  \lpatt(0.01 100)
\end{texdraw}
\end{figure}

\end{thm}

\begin{thm}[\rm]{Definitions}\label{definitions}
\begin{enumerate}
\item We will say that two vertices in a graph are {\it
neighbours} if and only if they are joined by an edge. \item A
graph in which no vertex is joined to itself by an edge is said to
be {\it simple}. \item A {\it finite graph} is one with finitely
many vertices and finitely many edges. \item The {\it degree of a
vertex} is the number of edges to which it is joined. \item The
{\it degree of a graph} is the maximum of the degrees of its
vertices. \item A {\it terminal vertex} is a vertex which is
connected to at most one other vertex. \item A {\it subgraph} $G'$
of the graph $G$ is a graph such that $V(G') \subseteq V(G)$ and
$E(G') \subseteq E(G)$. \item Let $G$ be a graph with the vertex
set $V= \lc x_1,\dots,x_n \rc$. For $W \subseteq V$ we define the
{\it induced subgraph} of $G$ with vertex set $W$ to be the graph
with vertex set $W$ which has an edge between any two vertices if
and only if there is an edge between them in $G$.

For example, if $G$ is as in definition (\ref{fig1}) and
$W=\{x_1,x_2,x_4\}$ then the induced subgraph of $G$ with vertex
set $W$ is

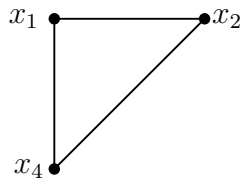
\begin{figure}[h]\label{fig2}
\begin{texdraw}
  \drawdim cm
  \textref h:R v:C
  \move(0 0)
  \move(6 0) \fcir f:0 r:0.075  \htext{$x_4 \ $}
  \lvec(8 2) \lvec(6 2) \lvec(6 0)
  \move(8.5 2) \htext{$x_2$} \move(8 2) \fcir f:0 r:0.075
  \move(6 2) \fcir f:0 r:0.075  \move(5.8 2) \htext{$x_1$}
\linewd 0.01 \lpatt(0.2 0.2)
  \lpatt(0.01 100)
\end{texdraw}
\end{figure}

\item A {\it path} from a vertex $x_1$ to a vertex $x_n$ is a
sequence of edges,
$$\{x_1,x_2\},\{x_2,x_3\},\dots\{x_{n-1},x_n\}.$$ Such a path is a
$cycle$ if $x_1,x_2,\dots,x_{n-1}$ are distinct and $x_1=x_n.$
 \item A graph will be said to be {\it connected} if any two of its vertices are
joined by a path of edges. A {\it connected component} is a
maximal connected subgraph. \item If $a\in V(G)$ then $G \setminus
\{a\}$ will denote the subgraph of $G$ which has vertex set $V(G)
\setminus \{a\}$ and all the edges of $G$ which do not feature
$a$. \item If $e=\{x,y\} \in E(G)$ then $G \setminus e$ will
denote the graph with the same vertex set as $G$ and edges $E(G)
\setminus e$.

\item A {\it{tree}} is a connected graph with no cycles. \item
 A {\it forest} is a graph whose connected components are all trees.

\end{enumerate}
\end{thm}

 Henceforth we will
simply say `graph' to mean a simple, finite graph with at most one
edge between any two edges, as all graphs discussed in this thesis
will be of this type.

 Given any graph we will
associate with it a polynomial ring and a monomial ideal (that is
an ideal generated by monomials) of that ring as follows:

\begin{thm}[\rm]{\Def}
 For a graph $G$ on vertices $x_1, \dots, x_n$ and a field $k$ we define $R_k(G)$ to be the polynomial
 ring over $k$ in the in $n$ indeterminants which we will also call $x_1,\dots,x_n$, \ie $R_k(G)=k[x_1, \dots,
 x_n]$, where $k$ is any field, and we define $I(G)$ to be the monomial ideal of $R_k(G)$
 generated by $\left\{x_ix_j \mid \{x_i,x_j\} \mathrm{\ is \ an \ edge \ of \ G}\right
 \}$. We call $I(G)$ the {\it graph ideal} of $G.$
\end{thm}

\begin{thm}[\rm]{Remark}
Every monomial ideal of $R_k=k[x_1,\dots,x_n]$ which is generated
by square free monomials of degree 2 may be interpreted as a graph
ideal. If $I=\langle x_{i_1}x_{j_1},\dots, x_{i_s}x_{j_s} \rangle
\subseteq R_k$ then $I=I(G)$ where $G$ is the graph with vertices
$x_1,\dots,x_n$ and edges $\{x_{i_1},x_{j_1}\},\dots,
\{x_{i_s},x_{j_s}\}$. Note that not all of the vertices
necessarily feature in any edge, \ie $G$ may have some isolated
vertices.
\end{thm}

\begin{thm}[\rm]{Definition}
Let $V=\{v_1,\dots,v_n\}$ be a finite set. A {\it simplicial
complex} $\G$ on $V$ is a set of subsets of $V$ such that
$\{v_1\},\dots,\{v_n\} \in \G$ and if $F \in \G$ and $G \subseteq
F$ then $G \in \G$. The elements of $\G$ are called $\it faces$
and a face $F$ has {\it dimension} $\dim F= |F|-1$ and the {\it
dimension} of $\G$, $\dim \G$, is the maximum of the dimensions of
its faces.
\end{thm}

\begin{thm}[\rm]{Definition}
Let $\G$ be a simplicial complex with vertex set $\{x_1,\dots,x_n
\}$. Let $R(\G)=k[x_1,\dots,x_n]$, for a field, $k$. The {\it
Stanley-Reisner ideal} of $\G$, denoted by $I(\G)$, is the ideal
of $R(\G)$ generated by all squarefree monomials $x_{i_1}\dots
x_{i_j}$ such that there is no face of $\G$ with vertices
$x_{i_1},\dots,x_{i_j}$. The {\it Stanley-Reisner Ring} (or {\it
face ring}) of $\G$ is defined to be the quotient ring
$k[\G]=R(\G)/I(\G)$.
\end{thm}

\begin{thm}[\rm]{Remark}\label{SRC}
Note that $I(G)=I(\Delta(G))$, the Stanley-Reisner ideal of the
simplicial complex $\Delta(G)$ which has faces
$$\left\{\left\{x_{i_1},\dots , x_{i_l} \right\}| \mathrm{ \ no \ }
\{x_{i_j},x_{i_k}\} \mathrm{\ is \ an \ edge \ of \ } G
\right\}.$$

Consequently graph ideals are a special case of Stanley-Reisner
ideals and we will henceforth write $   k[\D(G)]$ for
$R_k(G)/I(G)$. We attempt to use the combinatorial properties of
graphs to understand some of the algebraic properties of the
associated ideals and vice versa.
\end{thm}

\section{Minimal Free Resolutions}

We now define resolutions of modules and in particular minimal
free resolutions. For a graph $G$ we will be interested in
resolutions of $k[\D(G)]$ considered as an $R_k(G)$-module. These
resolutions and parts of the information they contain, such as
Betti numbers and projective dimensions, will be of significant
importance in what follows.

\begin{thm}[\rm]{\Def} A {\it free resolution} of an $R$- module $M$ is
a complex of free modules
$$\FI :  \dots \xrightarrow{\varphi_{n-1}} F_n \xrightarrow{\varphi_n}
\dots \xrightarrow{\varphi_2} F_1 \xrightarrow{\varphi_1} F_0 $$
which is exact and is such that $\coker \varphi_1=M$ . It is also
a {\it graded} free resolution if $R$ is a graded ring, the $F_i$
are graded free modules and the maps are homogeneous of degree 0.
A {\it finite free resolution of length $n$} is one in which
$F_{i}=0$ for all $i \geqslant n+1$  but $F_0,\dots,F_n$ are all non zero.
\end{thm}
If
$$\FI :  \dots \xrightarrow{\varphi_{n-1}} F_n \xrightarrow{\varphi_n}
\dots \xrightarrow{\varphi_2} F_1 \xrightarrow{\varphi_1} F_0 $$
and
$$\GI :  \dots \xrightarrow{\psi _{n-1}} G_n \xrightarrow{\psi_n}
\dots \xrightarrow{\psi_2} G_1 \xrightarrow{\psi_1} G_0 $$ are
free resolutions of $M$ such that for all $n$ there exists an
isomorphism of modules, $\theta_n:F_n \rightarrow G_n$, and for
all $n$ we have $\theta_{n-1}\varphi_n=\psi_n\theta_n$ then $\FI$
and $\GI$ will be said to be isomorphic resolutions.

\begin{thm}[\rm]{Remark} We will sometimes describe a free
resolution of $M$ slightly differently from the above and say that
an exact sequence with free modules $F_0,\dots,F_n,\dots$,

$$\FI :  \dots \xrightarrow{\varphi_{n-1}} F_n \xrightarrow{\varphi_n}
\dots \xrightarrow{\varphi_2} F_1 \xrightarrow{\varphi_1} F_0
\xrightarrow{\varphi_0} M \to 0$$
 is a free resolution of $M$.
\end{thm}

We shall focus our attention on minimal free resolutions.

\begin{thm}[\rm]{\Def}
A {\it minimal finite free resolution} of an $R$-module $M$ is one
with the smallest possible length and smallest possible rank for
each of the free modules. Minimal resolutions are unique up to
isomorphism. The rank of the $i$th free module in a minimal
resolution is called the $i$th {\it Betti number} of $M$.
\end{thm}

The following theorem guarantees the existence of graded finite
resolutions of the types of modules in which we are interested.

\begin{thm}{Theorem}(Hilbert Syzygy Theorem)\label{HST}
Let $R=k[x_1,\dots,x_n] $. \\Every finitely generated graded
$R$-module has a graded finite free resolution, of length no more
than $n$, by finitely generated free modules.
\end{thm}

We can use this idea of minimal free resolutions to define new
numerical invariants of graphs, using Betti numbers and also
projective dimension.

\begin{thm}[\rm]{\Def} \label{pddef}
The {\it projective dimension} of a module $M$, $\pd(M)$, is the
length of its minimal resolution (or $\infty $ if it has no finite
resolution). We define the projective dimension of a graph $G$ to
be the projective dimension of the $R_k(G)$-module $k[\D(G)]$, and
we will write $\pd^k(G)=\pd(k[\D(G)])$.
\end{thm}

The projective dimension of a graph will in general be dependent
on the characteristic of our choice of the field $k$. However in
many cases we will see that it is independent. If this is the case
we will simply write the projective dimension of $G$ as $\pd(G)$.
Results about the projective dimension of a graph are examined in
chapter 2.

\begin{thm}[\rm]{Remark}When we are
investigating minimal free resolutions of graph ideals it does not
matter which of the $R_k(G)$-modules $I(G)$ and $k[\D(G)]$ we
consider. This is because there is a bijective correspondence
between minimal free resolutions of $I(G)$ and minimal free
resolutions of $k[\D(G)]$. That is to say
$$ \FI \dots \to F_n \to
\dots \to F_1$$ is a minimal free resolution of $I(G)$ if and only
if
$$ \FI' \dots \to F_n \to
\dots \to F_1 \to R$$ is a minimal free resolution of $k[\D(G)]$.
Note that this implies $\pd(k[\D(G)])=\pd(I(G)) +1$.

\end{thm}

\section{Graded Betti Numbers}
Now we examine more closely Betti numbers by considering
$\mathbb{N}$-graded and $\mathbb{N}^n$-graded Betti numbers.
 We may also describe Betti numbers in terms of Tor. For an
$\mathbb{N}^n$-graded module $M$ over $R=k[x_1,\dots,x_n]$ ,
tensoring a minimal free resolution of $M$,
$$ 0\to R^{\beta_h} \to
\dots \to R^{\beta_1}, $$ with $k$ (considered as an $R$-module
via $k\cong R/\langle x_1,\dots,x_n \rangle$) produces a complex
$$ 0\to k^{\beta_h} \to
\dots \to k^{\beta_1}$$ where all the maps are zero. Taking $i$th
homology of this complex we obtain $\b_i(M)=\b_i^k(M)=\dim_k
\mathrm {Tor}_i^R(k,M)$.

For a (simple,  finite) graph $G$ we may write the minimal free
resolution of $k[\D(G)]$ as

$$ 0\to R^{\beta_h} \to
\dots \to R^{\beta_1} \to R^{\b_0}$$ where $R=R_k(G)$,
$\beta_i(k[\D(G)])=\dim_k(\mathrm{Tor}_i^R(k[\D(G)],k))$ and $\b_h
\neq 0$. Note that $\b_0(k[\D(G)])=1$. Observe that here we have
$\pd^k(G)=\pd(k[\D(G)])=h$ (and so $\pd^k(I(G))=h-1$). Note also
that $\mathrm{Tor}_i^R(k[\D(G)],k)$ inherits the $\mathbb{N}^n$-
grading. For $\mathbf{a} \in \mathbb{N}^n$ let
$(\mathrm{Tor}_i^R(k[\D(G)],k))_\mathbf{a}$ denote the $\mathbf{a}
$ graded component. This allows us to define
${\mathbb{N}}^n$-graded Betti numbers and $\mathbb{N}$-graded
Betti numbers which will in some circumstances be easier to handle
than the total Betti numbers.

\begin{thm}[\rm]{\Def} \label{bettor}
For $\mathbf{a} \in {\mathbb{N}}^n$ the $i$th Betti number of
degree $\mathbf{a}$ of $R/I$ is $$\b_{i,\mathbf{a}}^k= {\dim}_k(
\mathrm{Tor}_i^R(R/I,k))_{\mathbf{a}}.$$ For $d \in \mathbb{N}$
the $i$th Betti number of degree $d$ of $R/I$ is
$$\b_{i,d}^k =
\sum_{|\mathbf{a}|=d}{\dim}_k(\mathrm{Tor}_i^R(R/I,k))_{\mathbf{a}}
=\sum_{|\mathbf{a}|=d} \b^k_{i,\mathbf{a}}$$
            where, for $\mathbf{a}=(a_1,\dots,a_n)\in \mathbb{N}^n$,
$|\mathbf{a}|=a_1+\dots+a_n$.
\end{thm}

It will often be convenient to consider these graded Betti numbers
and recover the (total) Betti numbers by noting that
$$\b_i^k=\sum_{d \in \mathbb{N}}\b_{i,d}^k=\sum_{\mathbf{a} \in
\mathbb{N}^n}\b_{i,\mathbf{a}}^k.$$

\begin{thm}[\rm]{\Def}
We will write $\b_i^k(G)=\b_i^k(k[\D(G)])$ and say this is the
$i$th Betti number of the graph $G$. We will use similar
terminology for $\mathbb{N}$- and $\mathbb{N}^n$- graded numbers
of $k[\D(G)]$.
\end{thm}

\begin{thm}[\rm]{Remark}
The Betti numbers are in general going to be dependent on the
characteristic of the field $k$. However we shall see later that
for many of the families of graphs we deal with (cycles, complete
graphs and trees, for example) the $R_k(G)$-module $k[\D(G)]$ has
Betti numbers which are independent of the choice of $k$. When
this is this case we will simply write $\b_i(G)$, $\ \b_{i,d}(G)$
and $\b_{i,\mathbf{a}}(G)$ for the $i$th Betti numbers of $G$
(where $d \in \mathbb{N}$ and $\mathbf{a} \in \mathbb{N}^n$).
\end{thm}

\section{Hochster's Formula}
We now present a result of Hochster's which relates the Betti
numbers of the Stanley-Reisner ring $k[\D]$ (considered as a
module over $k[x_1,\dots,x_n]$) to the simplicial homology of
$\D$. This will be of much use when calculating the Betti numbers
of graph ideals. More extensive details may be found in \cite{BCP}
or \cite{MS}. We first introduce the following simplicial complex.

\begin{thm}[\rm]{Definition}
Let $I$ be a monomial ideal of $R=k[x_1,\dots,x_n].$ For
$\mathbf{b} \in \mathbb{N}^n$, define
$$K_{\mathbf{b}}(I)=\{ F \subseteq \{1,\dots,n\} \mid
\mathbf{x}^{\mathbf{b} -F} \in I\}.$$ Each face $F \in
K_{\mathbf{b}}(I)$ is identified with its characteristic vector in
$\{0,1\}^n$, \ie $F$ is identified with the element of $\{0,1\}^n$
which has $j$th component $1$ if and only if $j\in F$ for
$j=1,\dots,n$. Therefore if $\mathbf{b}=(b_1,\dots,b_n)$ then the
$j$th component of $\mathbf{b} -F$ is $b_j -1$ if $j\in F$ and
$b_j$ if $j \notin F$.
\end{thm}

The Betti numbers of $I$ can be found from $K_{\mathbf{b}}(I)$ as
follows.

\begin{thm}{Theorem}\cite{BCP}\label{bettK}
The ${\mathbb{N}^n}$-graded Betti numbers of $I$ are
$$\b_{i,\mathbf{b}}(I)= \dim_k \Hred_{i-1}(K_{\mathbf{b}}(I);k).$$
\end{thm}

\begin{proof}
We use the Koszul complex $K_{\bullet}$ which is a minimal free
resolution of the $R$-module $k\cong R/\langle x_1,\dots,x_n
\rangle$.
$$K_{\bullet}: \ \ 0 \to \wedge^nV \to \cdots \to \wedge^1V\to\wedge^0V \to
0,$$ where $V=\{ \lambda_1x_1+\dots+\lambda_nx_n \mid
\lambda_1,\dots,\lambda_n \in k\}$. The Betti number
$\b_{i,\mathbf{b}}= \dim_k ({\Tor}_i^R(I,k))_{\mathbf{b}}$ may be
found by taking the $i$th homology of the degree $\mathbf{b}$ part
of the complex
$$I \otimes K_{\bullet}: \ \ 0 \to I \otimes\wedge^nV \to \cdots \to I\otimes\wedge^1V\to I\otimes\wedge^0V \to
0.$$ There is a basis of $(I \otimes \wedge ^iV)_{\mathbf{b}}$
which consists of all the elements $$ \mathbf{x^b}/x_{j_1}\cdots
x_{j_i} \otimes x_{j_1}\wedge \cdots \wedge x_{j_i}$$ such that
$\mathbf{x^b}/x_{j_1}\cdots x_{j_i} \in I$. There is a one to one
correspondence between these expressions and the $(i-1)$-faces
$\{j_1,\dots,j_i\}$ of $K_{\mathbf{b}}(I)$. We see that $(I
\otimes K_{\bullet})_{\mathbf{b}}$ is the augmented chain complex
of $K_{\mathbf{b}}(I)$ and hence that $\b_{i,\mathbf{b}}(I)=
\dim_k \Hred_{i-1}(K_{\mathbf{b}}(I);k).$
\end{proof}

We now introduce the Alexander dual of a simplicial complex.

\begin{thm}[\rm]{\Def}\label{addef}
Let $\Delta$ be a simplicial complex with vertex set
$[n]=\{1,\dots,n\}$. The Alexander Dual, $\AD$, is defined to be
the simplicial complex with faces $$\left\{ F\subseteq
\left\{1,\dots , n \right\} |\left\{1, \dots, n \right\}
\backslash F \mbox{ is not a face of } \Delta \right\}.$$
\end{thm}

\begin{thm}[\rm]{\Def}
Let $\D$ be a simplicial complex with vertex set $\{1,\dots,n \}$
and let $V \subseteq \{1,\dots,n \}$. The restriction of $\D$ to
$V$ is the subcomplex of $\D$ $$\D_V=\{ F \in \D \mid F \subseteq
V \}.$$ If $\mathbf{b}\in\{0,1\}^n$ then the restriction of $\D$
to $\mathbf{b}$ is $\D_{\mathbf{b}}=\D_V$ where $V$ is the subset
of $\{1,\dots,n \}$ identified with $\mathbf{b}$, \ie if
$\mathbf{b}=(b_1,\dots,b_n)$ then $j \in V$ if and only if
$b_j=1$.
\end{thm}

When $I$ is a squarefree monomial ideal, \ie when $R/I=k[\D]$ for
some simplicial complex $\D$  theorem (\ref{bettK}) may be
reformulated as follows.

\begin{thm}{Theorem}\cite{H}(Hochster's Formula)\label{HochstersFormula}
Let $k[\D]=R(\D)/I(\D)$ be the Stanley-Reisner ring of the
simplicial complex $\D$. The non-zero Betti numbers of $k[\D]$ are
only in squarefree degrees $\mathbf{b}$ and may be expressed as
$$\b_{i,\mathbf{b}}(k[\D])=\dim_k \Hred_{|\mathbf{b}|-i-1}(\D_{\mathbf{b}};k).$$
Hence the i{\rm th} total Betti number may be expressed as
$$\b_i(k[\D])= \sum_{V \subseteq \{1,\dots,n\}} \dim_k
\Hred_{|V|-i-1}(\D_V;k).$$
\end{thm}

\begin{proof}
Let $\mathbf{b}=(b_1,\dots,b_n)$. If $b_j>1$ for some $j$ then
$K_{\mathbf{b}}(I)$ is a cone over $j$ and therefore
$\b_{i,\mathbf{b}}(k[\D])=\b_{i-1,\mathbf{b}}(I)=\dim_k
\Hred_{i-2}(K_{\mathbf{b}}(I);k)=0$ by Theorem (\ref{bettK}). If
$\mathbf{b}\in \{0,1\}^n$ and $V$ is the subset of $\{1,\dots,n\}$
identified with $\mathbf{b}$ (\ie if $\mathbf{b}=(b_1,\dots,b_n)$
then $j \in V$ if and only if $b_j=1$) then $\D_{\mathbf{b}}$ is
the Alexander dual of $K_{\mathbf{b}}(I)$. \ie $F\in
K_{\mathbf{b}}(I) \iff \mathbf{x}^{\mathbf{b}-F} \in I \iff F
\subseteq V {\rm and}\ V \setminus F \notin \D \iff F \in
(\D_{\mathbf{b}})^*.$ We now use the fact that, for any simplicial
complex $X$ with $m$ vertices $\Hred_i(X;k) \cong
\Hred_{m-i-3}(X^*;k)$ along with Theorem (\ref{bettK}) to obtain
$\b_{i,\mathbf{b}}(k[\D])=\dim_k
\Hred_{i-2}(K_{\mathbf{b}}(I);k)=\dim_k
\Hred_{|\mathbf{b}|-i-1}(\D_{\mathbf{b}};k).$ The total Betti
number is found by summing over all squarefree vectors, noting
there is a one to one correspondence between these vectors and
subsets of $\{1,\dots,n\}.$
\end{proof}





Hochster's formula can be used to find explicit descriptions of
the Betti numbers of $k[\D(G)]$ when the simplicial complex
$\D(G)$ has reduced homology which is easily understood. In
chapter 5 will we see this is the case for complete and complete
bipartite graphs. For other types of graphs it will be necessary
to use Theorem (\ref{ER}). This is a reformulation of Theorem
(\ref{HochstersFormula}) (see \cite{ER} for further details) in
terms of reduced homology of new simplicial complexes arising from
$\Delta$.

\begin{thm}[\rm]{\Def}\label{link}
The link of a face $F$ of $\Delta$, ${\Link}_{\Delta} F$, is the
simplicial complex consisting of all faces, $G$, which satisfy
\begin{enumerate}
\item $G \in \Delta$; \item $G \cup F \in \Delta$; \item $G \cap F
= \emptyset $.
\end{enumerate}
\end{thm}

\begin{thm}{Theorem}\cite{ER}\label{ER}
The $\mathbb{N}$-graded Betti numbers of $k[\D]$ may be expressed
by

$$\b_{i,d}(k[\D])=\sum_{F \in \AD, |F|=n-d} \dim_k \Hred_{i-2}({\Link_{\Delta^*}} F;k)
,$$ and the total Betti numbers by
$$\b_i(k[\D])=\sum_{F \in \AD} \dim_k \Hred_{i-2}({\Link_{\Delta^*}} F;k) .$$
\end{thm}

\begin{proof}
We use Hochster's formula (\ref{HochstersFormula}),
$$\b_{i,d}(k[\D])=\sum_{|\mathbf{b}|=d} \dim_k
\Hred_{|\mathbf{b}|-i-1}(\D_{\mathbf{b}};k).$$

For $\mathbf{b}$ appearing in this sum, with $\mathbf{b}$
corresponding to $\tau \subseteq \{1,\dots,n\}$, let
$F=\{1,\dots,n\} \setminus \tau$ (so $|F|=n-d$). If $\tau$ is a
face of $\D$ then $\D_{\tau}$ is a simplex and will have no
reduced homology. Hence it may be assumed that $\tau$ is not a
face of $\D$. From the definition (\ref{addef}) we see that $F$ is
a face of $\AD$. It is therefore enough to show that
$\Hred_{i-2}(\linkad F;k) \cong \Hred_{|\tau|-i-1}(\D_{\tau};k)$.
This follows as the complexes $\linkad F$ and $(\D_{\tau})^*$ are
isomorphic and $\Hred_{i-2}((\D_{\tau})^*;k) \cong
\Hred_{|\tau|-i-1}(\D_{\tau};k)$. The total Betti number is found
by summing all graded Betti numbers.
\end{proof}

The above theorem finds Betti numbers from the reduced homology
modules of links of faces of the Alexander Dual. For graph ideals
we will see later (in chapter 6) how we can find a correspondence
between the subgraphs and the links of faces of the Alexander
dual.

\chapter{Projective Dimension}
\setcounter{satz}{0}
\section{Projective Dimensions of Subgraphs}

In this chapter we look at various properties of the projective
dimension of graphs, as defined in (\ref{pddef}).

\begin{thm}[\rm]{Remark}
It is worth noting that there really are graphs whose Betti
numbers (and so projective dimension) depend on the characteristic
of the choice of field.
\end{thm}

\begin{thm}[\rm]{Definition}
Let $\G$ be a simplicial complex. The {\it barycentric
subdivision} of $\G$ is a simplicial complex which has as vertices
the non-empty faces of $\G$ and faces $\{ \{f_1,\dots,f_r\} |
f_1,\dots,f_r \in \G \ \mathrm { and } \ f_1 \subset f_2 \subset
\dots \subset f_r \}$.
\end{thm}

\begin{thm}[\rm]{Example}
We show here how to construct a graph with Betti numbers dependent
on $\mathrm{char\ } k$. Let $\G$ be the simplicial complex shown
in figure 2.1, see \cite{BH} p.236. This is a triangulation of the
real projective plane $\mathbb{P}^2$. The 1st reduced homology
module of this is
$$
\Hred_1(\mathbb{P}^2;k)= \left\{
\begin{array}{lll}
k &\mathrm{if}& \mathrm{char} \ k =2 \\
0 &\mathrm{if}& \mathrm{char} \ k \neq 2.
\end{array}
\right.
$$
\begin{figure}[h]
\begin{texdraw}
  \drawdim cm
  \textref h:R v:C
  \move(0 0)
  \move(6 -1) \fcir f:0 r:0.075
 \lvec(9 -1) \fcir f:0 r:0.075
 \lvec(11 2) \fcir f:0 r:0.075
 \lvec(9 5) \fcir f:0 r:0.075
 \lvec(6 5) \fcir f:0 r:0.075
 \lvec(4 2) \fcir f:0 r:0.075
\lvec(6 -1) \fcir f:0 r:0.075 \lfill f:0.93 \lvec(6.75 1.25) \fcir
f:0 r:0.075  \lvec(8.25 1.25) \fcir f:0 r:0.075 \lvec(7.5 2.75)
\fcir f:0 r:0.075  \lvec(6.75 1.25) \fcir f:0 r:0.075 \lvec(4 2)
\fcir f:0 r:0.075 \move (7.5 2.75) \lvec(6 5) \fcir f:0 r:0.075
\move (7.5 2.75) \lvec(9 5) \fcir f:0 r:0.075 \move (7.5 2.75)
\lvec(11 2) \fcir f:0 r:0.075 \move (6.75 1.25) \lvec(6 5) \fcir
f:0 r:0.075
 \move (8.25 1.25) \lvec(6 -1) \fcir f:0 r:0.075
  \move (8.25 1.25) \lvec(9 -1) \fcir f:0 r:0.075
   \move (8.25 1.25) \lvec(11 2) \fcir f:0 r:0.075
\move (6 -1.3) \htext{ $x_1$} \move (9.2 -1.3) \htext{ $x_2$}
\move (3.8 2) \htext{ $x_3$}\move (11.6 2) \htext{ $x_3$} \move (6
5.3) \htext{ $x_2$}\move (9 5.3) \htext{ $x_1$} \move (6.6 1.1)
\htext{ $x_4$}\move (8.8 1.1) \htext{ $x_5$}\move (7.73 3.2)
\htext{ $x_6$}
 \linewd 0.01 \lpatt(0.2 0.2)
  \lpatt(0.01 100)
\end{texdraw}
\caption {$\G$  }
\end{figure}
We now form a new simplicial complex $\D$ by taking the
barycentric subdivision of $\G$. We may write $\D=\D(G)$ where $G$
is the graph with the same vertices as $\D$, that is the faces of
$\G$, and edges $\{f,g\}$ where $f,g \in \G$ with $f \nsubseteq g$
and $g\nsubseteq f$. This Graph has $31$ vertices. We now consider
the $29$th Betti number of degree $31$ using the formula of
(\ref{HochstersFormula}). Let $[31]=\{1,\dots,31\}$ and note that
$\D_{[31]}=\D$.
\begin{eqnarray*}
\b^k_{29,31}(G)&=&\sum_{V \subseteq [31], \ |V|=31} \dim_k
\Hred_{|V|-29-1}(\D_V;k)\\
&=&\dim_k\Hred_{31-29-1}(\D_{[31]};k)\\
&=&\dim_k \Hred_1(\D;k)\\
&=&\dim_k \Hred_1(\mathbb{P}^2;k)\\
&=& \left\{\begin{array}{lll}
           1 &\mathrm{if}& \  \char \  k=2\\
           0 &\mathrm{if}& \  \char \ k \neq2.
           \end{array} \right.
\end{eqnarray*}

\end{thm}

We now investigate how the projective dimension of a graph is
affected by some simple transformations of the graph, such as
deleting an edge. In what follows $G$ is a graph on $n$ vertices
and $R_k(G)=k[x_1,\dots,x_n]$.

\begin{thm}{\Prop}
Let $G$ be a graph with an edge $e=\{a,b\}$ such that $b$ is a
terminal vertex. Then $\pd^k(G \setminus \{b\}) \leqslant
\pd^k(G).$
\end{thm}

\begin{proof}
First note that $V(G) \supseteq V(G \setminus \{b\})=V(G)
\setminus \{b\}$. Let $p=\pd^k (G \setminus \{b\})$. We now use
Theorem (\ref{HochstersFormula}) to see that
\begin{eqnarray*}
\b_p(G)&=&\sum_{W \subseteq V(G)} \dim_k \Hred_{|W|-p-1}(\D_W;k)\\
       &\geqslant&\sum_{W \subseteq V(G) \setminus \{b\}} \dim_k \Hred_{|W|-p-1}(\D_W;k)\\
       &=&\b_p(G \setminus \{b\}) \neq 0.
\end{eqnarray*}
Therefore we may conclude that $\pd^k(G) \geqslant \pd^k(G
\setminus\{b\})$.
\end{proof}

Here we supply an example to show that with the situation as above
it is possible either to have $\pd(G \setminus \{b\}) < \pd(G)$ or
to have $\pd(G \setminus \{b\}) = \pd(G)$.

\begin{thm}[\rm]{Example}
Let $G_1,G_2$ and $G_3$ be the following graphs

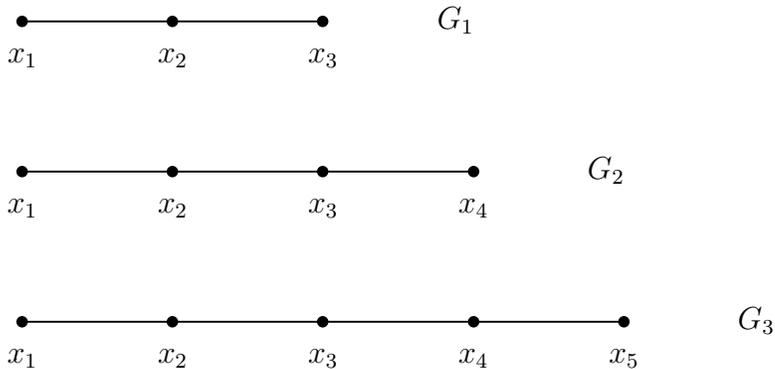
\begin{figure}[h]\label{fig3}
\begin{texdraw}
  \drawdim cm
  \textref h:R v:C
  \move(0 0) \fcir f:0 r:0.075 \lvec(2 0)\move(0.2 -.5)\htext{$x_1$}
  \move(2 0) \fcir f:0 r:0.075 \lvec(4 0)\move(2.2 -.5)\htext{$x_2$}
  \move(4 0) \fcir f:0 r:0.075 \move(4.2 -.5)\htext{$x_3$} \move(6
  0) \htext{$G_1$}
  \move(0 -2) \fcir f:0 r:0.075 \lvec(2 -2)\move(0.2 -2.5)\htext{$x_1$}
  \move(2 -2) \fcir f:0 r:0.075 \lvec(4 -2)\move(2.2 -2.5)\htext{$x_2$}
  \move(4 -2) \fcir f:0 r:0.075 \lvec(6 -2)\move(4.2 -2.5)\htext{$x_3$}
  \move(6 -2) \fcir f:0 r:0.075 \move(6.2 -2.5)\htext{$x_4$}
  \move(8 -2) \htext{$G_2$}
  \move(0 -4) \fcir f:0 r:0.075 \lvec(2 -4)\move(0.2 -4.5)\htext{$x_1$}
  \move(2 -4) \fcir f:0 r:0.075 \lvec(4 -4)\move(2.2 -4.5)\htext{$x_2$}
  \move(4 -4) \fcir f:0 r:0.075 \lvec(6 -4)\move(4.2 -4.5)\htext{$x_3$}
  \move(6 -4) \fcir f:0 r:0.075 \lvec(8 -4)\move(6.2 -4.5)\htext{$x_4$}
  \move(8 -4) \fcir f:0 r:0.075 \move(8.2 -4.5)\htext{$x_5$}
  \move(10 -4) \htext{$G_3$}

\linewd 0.01 \lpatt(0.2 0.2)
  \lpatt(0.01 100)
\end{texdraw}
\end{figure}
These graphs are all trees and we shall see later that they
therefore have projective dimension independent of the
characteristic of our choice of field. The projective dimensions,
which can be found from Corollary (\ref{linedim}), are $\pd
(G_1)=2$, $\pd (G_2)=2$ and $\pd(G_3)=3$. So we have $\pd(G_3
\setminus \{x_5\}) <\pd (G_3)$ and $\pd(G_2 \setminus
\{x_4\})=\pd(G_2)$.
\end{thm}

\begin{thm}[\rm]{Remark}The above example raises further questions
about the effect that adding or deleting edges has on the
projective dimension of a graph. For instance, what happens if an
edge is subdivided? By this we mean placing a new vertex in the
middle of an edge thereby replacing a single edge with two edges.
Whenever a new graph is produced in this way is its projective
dimension necessarily at least as big as that of the initial
graph? This has not be proved nor has a counter example been
produced.

\end{thm}

We now provide an example to show that this operation of
subdividing an edge can cause an arbitrarily large increase in
projective dimension.

    \begin{thm}[\rm]{Example}
For $r,s \geqslant1$ let $G$ be the graph shown in figure 2.2. We
will use some facts which will be proved in chapters 5 and 9 to
find the projective dimension of $G$. First we note that $G$ is a
tree and so by Theorem (\ref{BettiTrees}) its projective dimension
is independent of the field $k$. We now use Theorem
(\ref{pdtrees}) to obtain
\begin{eqnarray*}
\pd G &=& \max \{ \pd(G \setminus \{x_r\}),\pd (G \setminus
\{x_1,\dots,x_r,a,b\}) +(r+1) \}\\
      &=& \max \{\pd(G \setminus \{x_r\}),r+1\}\\
      & &\mbox{(as $G \setminus
\{x_1,\dots,x_r,a,b\}$ is a graph with no edges) }\\
      &=& \max \{ \pd (G \setminus \{x_r,x_{r-1}\}), r, r+1 \}\\
      &\dots& \\
      &=& \max \{\pd (G \setminus \{x_1,\dots, x_r\}),r+1\}\\
      &=& \max \{s+1,r+1\}\\
      & & \mbox{(by Proposition (\ref{Star}) as $G \setminus \{x_1,\dots, x_r\}$ is the star
      graph $S_{s+1}$).}
\end{eqnarray*}
Now we calculate the projective dimension of the graph formed by
subdivision of the edge $\{a,b\}$. Let $H$ be the graph pictured
in figure 2.3.
\begin{figure}
\begin{texdraw}
  \drawdim cm
  \textref h:R v:C
  \move(0 0)
  \move(4 0) \fcir f:0 r:0.075 \lvec (6 0 ) \fcir f:0 r:0.075
  \move(4 0 ) \lvec (2.5 1.5) \fcir f:0 r:0.075
  \move(4 0 ) \lvec (2.5 -1.5) \fcir f:0 r:0.075
  \move(6 0 ) \lvec (7.5 1.5) \fcir f:0 r:0.075
  \move(6 0 ) \lvec (7.5 -1.5) \fcir f:0 r:0.075
  \move (4.2 -0.2) \htext{$a$}
  \move (5.8 -0.2) \htext{$b$}
  \move (2.2 1.5)  \htext {$x_1$}
  \move (2.2 -1.5)  \htext {$x_r$}
  \move (8 1.5)  \htext {$y_1$}
  \move (8 -1.5)  \htext {$y_s$}
  \move (2.5 .75) \fcir f:0 r:0.02
  \move (2.5 0) \fcir f:0 r:0.02
  \move (2.5 -.75) \fcir f:0 r:0.02
  \move (7.5 .75) \fcir f:0 r:0.02
  \move (7.5 0) \fcir f:0 r:0.02
  \move (7.5 -.75) \fcir f:0 r:0.02

  \linewd 0.01 \lpatt(0.2 0.2)
  \lpatt(0.01 100)
\end{texdraw}
\caption{$G$}
\end{figure}

\begin{figure}
\begin{texdraw}
  \drawdim cm
  \textref h:R v:C
  \move(0 0)
  \move(4 0) \fcir f:0 r:0.075 \lvec (5 0.7) \fcir f:0 r:0.075
  \lvec (6 0) \fcir f:0 r:0.075
  \move(4 0 ) \lvec (2.5 1.5) \fcir f:0 r:0.075
  \move(4 0 ) \lvec (2.5 -1.5) \fcir f:0 r:0.075
  \move(6 0 ) \lvec (7.5 1.5) \fcir f:0 r:0.075
  \move(6 0 ) \lvec (7.5 -1.5) \fcir f:0 r:0.075
  \move (4.2 -0.2) \htext{$a$}
  \move (5.8 -0.2) \htext{$b$}
  \move (2.2 1.5)  \htext {$x_1$}
  \move (2.2 -1.5)  \htext {$x_r$}
  \move (8 1.5)  \htext {$y_1$}
  \move (8 -1.5)  \htext {$y_s$}
  \move (2.5 .75) \fcir f:0 r:0.02
  \move (2.5 0) \fcir f:0 r:0.02
  \move (2.5 -.75) \fcir f:0 r:0.02
  \move (7.5 .75) \fcir f:0 r:0.02
  \move (7.5 0) \fcir f:0 r:0.02
  \move (7.5 -.75) \fcir f:0 r:0.02
\move(5 0.4) \htext{$c$}
  \linewd 0.01 \lpatt(0.2 0.2)
  \lpatt(0.01 100)
\end{texdraw}
\caption{$H$}
\end{figure}
Again we note this is a tree and so has projective dimension
independent of our choice of base field for the polynomial ring
$R(H)$. By Theorem (\ref{pdtrees})
\begin{eqnarray*}
\pd H &=& \max \{ \pd (H \setminus \{x_r\}), \pd ( H \setminus
\{x_1,\dots,x_r,a,c\}) +(r+1) \}\\
      &=& \max \{ \pd (H \setminus \{x_r\}),  s +r+1 \}\\
      & &\mbox{(using Proposition (\ref{Star}) as $H \setminus
\{x_1,\dots,x_r,a,c\}$ is the star graph $S_{s}$)}\\
      &=& \max \{ \pd(H \setminus \{x_r,x_{r-1}\}),s+r,s+r+1\}\\
      &\dots& \\
      &=& \max \{ \pd(H \setminus \{x_1,\dots,x_r\}),s+r+1\}.\\
\end{eqnarray*}
From the first part of this example we see that $\pd(H \setminus
\{x_1,\dots,x_r\})=\max \{s+1,1+1\}=s+1$. Hence $\pd H = \max
\{s+1,s+r+1\}=s+r+1$. This demonstrates that subdividing
 an edge
can lead to a large growth in projective dimension. In this
example if we have $r=s$ then the subdivision of $\{a,b\}$
(almost) doubles the projective dimension.

\end{thm}

\section{Projective Dimension of Unions of Graphs}

We will now see that the projective dimension of a graph is the
sum of the projective dimensions of its connected components.

\begin{thm}{\Prop}\label{pdsum}
If $G$ is the disjoint union of the two graphs $G_1$ and $G_2$
then $\pd G =\pd G_1 + \pd G_2$.
\end{thm}

\begin{proof}
Suppose that $G_1$ has vertices $x_1,\dots,x_a$ and $G_2$ has
vertices $y_1,\dots,y_b$. Let
$$\begin{array}{l}
R_1=R(G_1)=k[x_1,\dots,x_a] \ \mathrm{and}\   I_1=I(G_1) \subseteq R_1, \\
R_2=R(G_2)=k[y_1,\dots,y_b] \ \mathrm{and}\   I_2=I(G_2)\subseteq R_2,  \\
\ R=R_k(G)=k[x_1,\dots,x_a,y_1,\dots,y_b].
\end{array}$$
We now select a maximal homogeneous regular sequence,
$f_1,\dots,f_r$, on $R_1/I_1$ and a maximal homogeneous regular
sequence, $g_1,\dots,g_s$, on $R_2/I_2$. We now show that the
images of $f_1,\dots,f_r,g_1,\dots,g_s$ in $R$ form a maximal
regular sequence on the $R$-module $R_1/I_1 \otimes_k R_2/I_2.$ We
first note that there is an injection
$$\frac{R_1}{I(G_1)+\langle f_1,\dots
f_{i-1}\rangle} \xrightarrow{f_i} \frac{R_1}{I(G_1)+ \langle
f_1,\dots,f_{i-1} \rangle}$$ which is multiplication by $f_i$ for
$i=1,\dots,r$. Tensoring with the $k$-vector space $R_2/I_2$ will
preserve this injection as $R_2/I_2$ is a free and therefore a
flat $k$-module. Hence the images of $f_1,\dots,f_r$ in $R$ (in
the obvious order) are a regular sequence on $R_1/I_1 \otimes_k
R_2/I_2$. Similarly we have an injection
$$\frac{R_2}{I(G_2)+\langle g_1,\dots g_{i-1}\rangle}
\xrightarrow{g_i} \frac{R_2}{I(G_2)+ \langle g_1,\dots,g_{i-1}
\rangle}.$$ This injection is preserved when we tensor over $k$
with $R_1/(I_1 +\langle f_1,\dots,f_r\rangle).$ Therefore the
images in $R$ of $f_1,\dots,f_r,g_1,\dots g_s$ form a regular
sequence on $R_1/I_1 \otimes_k R_2/I_2.$

We claim that this is a maximal regular sequence. Let
$$T_1=R_1/(I_1+\langle f_1,\dots,f_r \rangle),\
T_2=R_2/(I_2+\langle g_1,\dots,g_s \rangle)$$ and $T=T_1 \otimes_k
T_2$. As $f_1,\dots,f_r$ is a maximal regular sequence on
$R_1/I_1$ we have $\depth_{m_1} T_1=0$ where $m_1=\langle
x_1,\dots,x_a \rangle \subseteq R_1.$ Hence $m_1 \subseteq
{\Zdv}_{R_1} (T_1)=\bigcup_{P \in {\Ass}(T_1)}P.$ By prime
avoidance and maximality of $m_1$ we conclude that $m_1$ is an
associated prime of $T_1$. So there must be some non-zero $\a \in
T_1$ such that ${\rm Ann}_{R_1} \a =m_1$. Similarly there is a
non-zero $\b \in T_2$ such that ${\rm Ann}_{R_2} \b=m_2$, where
$m_2=\langle y_1,\dots,y_b \rangle \subseteq R_2$. We now consider
$\a \otimes \b \in T$, which must be non-zero. To see this,
suppose that $\{v_i\}_{i \in I}$ is a basis for the $k$ vector
space $T_1$ and that $\{w_j\}_{j \in J}$ is a basis of $T_2$ so
that $\{v_i \otimes w_j \}_{(i,j)\in I \times J} $ is a basis of
$T$. We can write $\a \otimes \b =\sum_i a_iv_i \otimes \sum_j
b_jw_j=\sum_{(i,j)}(a_ib_j)(v_i \otimes w_j) \neq 0$. We note that
the ideal $m=\langle x_1,\dots,x_a,y_1,\dots,y_b \rangle$ of $R$
is contained in the annihilator of $\a \otimes \b$. Also ${\rm
Ann}_R (\a \otimes \b) \neq R$ because $\a \otimes \b \neq 0$. As
$m$ is maximal we deduce that $m={\rm Ann}_R (\a \otimes \b)$ and
is an associated prime of $T$. This is equivalent to the statement
$\depth_m T=0$ and so we conclude that
$f_1,\dots,f_r,g_1,\dots,g_s$ is a maximal regular sequence on
$T$.

Now we note that $R_1/I_1 \otimes_k R_2/I_2 \cong R/(I_1R +I_2R)
=k[\D(G)]$ and use the formula $\pd M + \depth M = \depth S$, (see
\cite{BH}), for a polynomial ring $S$ and a finitely generated
$S$-module $M$. We now deduce
\begin{eqnarray*}
\pd^k(G)&=&\pd(R/I)\\
      &=&\depth R - \depth(R/I)\\
      &=&\depth R- \left(\depth (R_1/I_1) + \depth(R_2/I_2)
      \right)\\
      &=&\depth R - \left(\depth R_1 - \pd(R_1/I_1) + \depth R_2 -
      \pd(R_2/I_2) \right)\\
      &=& r+s -(r- \pd(R_1/I_1) +s - \pd(R_2/I_2) )\\
      &=& \pd(R_1/I_1) + \pd(R_2/I_2)\\
      &=& \pd^k(G_1) +\pd^k(G_2).
\end{eqnarray*}

\end{proof}

\begin{thm}{\Cor}
Let $G$ be a graph whose connected components are $G_1,\dots,G_m$.
The projective dimension of $G$ is the sum of the projective
dimensions of $G_1,\dots,G_m$, \ie $\pd (G) =\sum_{j=1}^m
\pd(G_j).$
\end{thm}

\begin{proof}
This follows from Proposition (\ref{pdsum}). We can write $G$ as
the disjoint union of the two graphs $G_1$ and $(G_2 \cup \dots
\cup G_m)$ to see that $ \pd(G)=\pd(G_1) + \pd(G_2 \cup \dots \cup
G_m).$ Carry on in this way to obtain
$\pd(G)=\pd(G_1)+\dots+\pd(G_m)$.
\end{proof}

The Hilbert Syzygy Theorem (\ref{HST}) provides a useful bound for
the projective dimension of the modules in which we shall be
interested. The ideal $I(G)$ is homogeneous with respect to the
$\mathbb{N}^n$-grading on $R_k(G)$. (i.e. $x_i$ has degree $e_i$,
the $i$th standard basis vector of $\mathbb{Z}^n$). This grading
is inherited by the $R_k(G)$-module $k[\D(G)]$. The Hilbert Syzygy
Theorem guarantees the existence of finite free resolutions of
$k[\D(G)]$. In particular it guarantees the existence of finite
free resolutions of length no larger then the number of vertices
of $G$ and so the projective dimension can be no larger than the
number of vertices of $G$. \ie it ensures that $\pd^k(G) \leqslant
|V(G)|$.

\chapter{Cellular Resolutions}
\setcounter{satz}{0} In this chapter we consider cellular
resolutions of monomial ideals. Typically these are far from
minimal. We show that even for some apparently simple graph ideals
there may be no minimal cellular resolution.

We now briefly describe this construction. For a more detailed
account see \cite{BS} and \cite{BPS}.
\section{Definition of Cellular Resolutions}

Let $R=k[x_1,\dots,x_n]$ and let $M$ be a monomial ideal over $R$
with generating set $\{m_j =\textbf{x}^{\textbf{a}_j} \mid j \in J
\}$, for some index set $J$.

Let $X$ be a simplicial complex with vertex set $J$ and an
incidence function $\e(F,F')$ on pairs of faces. The function $\e$
takes values in $\{0,1,-1\}$, is such that $\e(F,F')=0$ except
when $F'$ is a facet of $F$, $\e(\{j\},\emptyset)=1$ for all $j\in
J$, and
$$\e(F,F_1)\e(F_1,F') + \e(F,F_2)\e(F_2,F')=0$$ for any codimension
2 face $F'$ of $F$ and $F_1$ and $F_2$ are the facets of $F$ which
contain $F'$. The cellular complex $\textbf{F}_X$ is the
$R$-module $\oplus_{F \in X,F \neq \emptyset}RF$ where $RF$ is the
free $R$-module with one generator $F$ in degree $\textbf{a}_F$
with differential $$\partial F =\sum_{F' \in X, F' \neq \emptyset}
\e(F,F')\frac{m_F}{m_{F'}} F'.$$

\begin{thm}[\rm]{Definition}
The degree of a face $F$ of the complex $X$ is the exponent vector
of the monomial $m_F:= {\rm{lcm}}\{m_j \mid j \in F \}$. Let
$\preceq$ be the partial order on $\mathbb{Z}^n$ where
$(\a_1,\dots,\a_n) \preceq (\b_1,\dots,\b_n)$ if and only if $\a_i
\leqslant \b_i$ for all $i=1,\dots,n$. We define $X_{\preceq \a}$
to be the subcomplex of $X$ on the vertices of degree $\preceq
\a$.
\end{thm}

\begin{thm}{\Prop}\label{cell}\cite{BS}
The complex $\textbf{F}_X$ is a free resolution of $M$ if and only
if $X_{\preceq \textbf{b}}$ is acyclic over $k$ in all degrees
$\textbf {b}$. If this is the case then $\textbf{F}_X$ is called a
cellular resolution of $M$.
\end{thm}

Typically we may find many possible cellular resolutions for
graphs ideals. Unfortunately there may be no cellular resolution
for a graph ideal  which is also a minimal resolution, as the
following example demonstrates.

\section{A Graph with no Minimal Cellular \\ Resolution}
\begin{thm}[\rm]{Example}
We consider the cycle on four vertices $C_4$.
\begin{figure}[h]\label{cycle}
\begin{texdraw}
  \drawdim cm
  \textref h:R v:C
  \move(-6 0)
  \move(0 0) \fcir f:0 r:0.075  \htext{$4 \ $} \lvec(2 0) \fcir f:0 r:0.075
  \lvec(2 2) \fcir f:0 r:0.075   \lvec(0 2) \fcir f:0 r:0.075
  \move(0 2) \fcir f:0 r:0.075  \htext{$1 \ $} \lvec(0 0)
  \move(2.5 2) \htext{$2$} \move(2.5 0) \htext{$3$}

  \linewd 0.01 \lpatt(0.2 0.2) \lpatt(0.01 100)

\end{texdraw}
\caption{$C_4.$}
\end{figure}

Let $R=R(C_4)$. The $R$-module $I(C_4)$ has minimal free
resolution $0 \to R \to R^4 \to R^4 $. We now show that there can
be no simplicial complex which provides us with a cellular
resolution the same as this minimal resolution.

Suppose that $X$ is a simplicial complex which does this. As we
are assuming this gives a cellular resolution which coincides with
the minimal resolution we can see how many faces of each dimension
$X$ must have. Hence $X$ has four vertices, four edges and one
2-dimensional face. Therefore $X$ must be as shown in figure
$3.2$. Each vertex of $X$ will have degree corresponding to one of
the minimal generators of $I(C_n)$.
\begin{figure}[]\label{fig3}
\begin{texdraw}
  \drawdim cm
  \textref h:R v:C
\move(-5 0)
 \move (0 0) \lvec (2 0) \lvec(0 2) \lvec(0 0) \lfill f:0.8
 \move (2 0) \lvec(4 0)
 \move(0 -0.25) \htext{$a$}
 \move(2 -0.25) \htext{$c$}
 \move(0 2) \htext{$b$}
 \move(4 -0.250) \htext{$d$}
\linewd 0.01 \lpatt(0.2 0.2)
  \lpatt(0.01 100)
\end{texdraw}
\caption {$X$.}
\end{figure}
For example suppose $a$ has degree $(0,1,1,0)$ (corresponding to
$x_2x_3$), $b$ has degree $(1,1,0,0)$ (corresponding to $x_1x_2$),
$c$ has degree $(0,0,1,1)$ (corresponding to $x_3x_4$) and $d$ has
degree $(1,0,0,1)$ (corresponding to $x_1x_4$). In this case
$X_{\preceq (1,1,0,1)}$ is not acyclic, it is the two disjoint
vertices $b$ and $d$. We will have a similar result whichever way
we assign degrees to the vertices of $X$. So by Theorem
(\ref{cell}) there is no simplicial complex which gives a cellular
resolution which is also a minimal resolution in this case.
\end{thm}

This example suggests that the problem of finding minimal free
resolutions of graph ideals is more difficult than the problem of
finding their Betti numbers. However we may use a cellular
resolution to make the following statement about Betti numbers of
a monomial ideal.

\section{Results about Betti Numbers Using \\ Cellular Resolutions}

\begin{thm}{Theorem}\label{greater2i}
Let $I$ be an ideal of $R=k[x_1,\dots,x_n]$ which is generated by
squarefree monomials all of which have degree not exceeding $d \in
\mathbb{N}$. Then $\b_{i,j}=0$ whenever $j>di$.
\end{thm}

\begin{proof}
We use Taylor's resolution (see \cite{BPS})). This is the cellular
resolution $\textbf{F}_X$ of $I$ where $X$ is the full simplex
with vertex set a minimal set of generators of $I$. By Proposition
(\ref{cell}) this is a free resolution of $I$ because $X_{\preceq
\mathbf{b}}$ is acyclic for all degrees $\mathbf{b}$ as $X$ is a
simplex.

We consider the ${\mathbb{N}}$-graded degree $j$ part of
$\mathbf{F}_X$ at the $i$th position. This is the direct sum of
the free $R$-modules generated by the $(i-1)$-faces of $X$ which
have degree $j$. The rank of this free module must be at least
$\b_{i,j}(I)$. Note that a face of $X$ has $\mathbb{N}$-graded
degree $j$ if and only if it has $\mathbb{N}^n$-graded degree
$\mathbf{a}=(a_1,\dots,a_n)$ such that $a_1+\dots+a_n=j$. Suppose
that $G \in X$ is such a face. The $\mathbb{N}^n$ degree of $G$ is
an element of $\{0,1\}^n$, the exponent vector of the lowest
common multiple of its vertices. At most $di$ of the coordinates
of this vector are $1$ (this is when the vertices of $G$ are all
of degree $d$ and are pairwise coprime). Hence the
$\mathbb{N}$-graded degree is at most $di$ and so $\b_{i,j}=0$ for
all $j>di$.
\end{proof}

We can express the $i$th Betti number of degree $2i$ in terms of
graph theoretical properties as follows.

\begin{thm}{Theorem}
For $i \geqslant 1$ and $G$ any graph $\b_{i,2i}$ is the number of
induced subgraphs of $G$ which are $i$ disjoint edges.
\end{thm}

\begin{proof}
As above we let ${\mathbf{F}}_{X}$ be Taylor's resolution of $I$
$${\mathbf{F}}_X : \cdots \to F_{i+1} \xrightarrow{{\partial}_{i+i}} F_i \xrightarrow{{\partial}_{i}}
F_{i-1} \to\cdots$$ We can define the differential $\partial_i$ by
$$
\begin{array}{lll}
\partial_i\{e_{j_1},\dots,e_{j_i}\}= \sum_{k=1}^i(-1)^k
\frac{{\mathrm{lcm}}\{e_{j_1},\dots ,
e_{j_i}\}}{{\mathrm{lcm}}\{e_{j_1},\dots
e_{j_{k-1}},e_{j_{k+1}},\dots ,e_{j_i}\}}\{e_{j_1},\dots
e_{j_{k-1}},e_{j_{k+1}},\dots ,e_{j_i}\}.
\end{array}
$$
 So the differential
$\overline{\partial_i}$  in $\mathbf{F}_X \otimes_R R/M$, where
$M=\langle x_1,\dots,x_n \rangle$, is defined by
$$\overline{\partial_i}\{e_{j_1},\dots,e_{j_i}\}=
\sum (-1)^k \{e_{j_1},\dots e_{j_{k-1}},e_{j_{k+1}},\dots
,e_{j_i}\},$$ where the sum is over all edges $e_{j_k}$ such that
the vertices of $e_{j_k}$ are in $e_{j_1},\dots
e_{j_{k-1}},e_{j_{k+1}},\dots ,e_{j_i}.$ The generators of $F_i$
in degree $2i$ are the sets of $i$ edges of $G$ which feature $2i$
vertices, \ie sets of $i$ disjoint edges of $G$. The image of such
a generator in $\mathbf{F}_X \otimes_R R/M$ will be in the kernel
of $\overline{\partial_i}$. If the disjoint edges
$\{e_{j_1},\dots,e_{j_i}\}$ do not form an induced subgraph of
$G$, \ie if there is another edge $e$ whose vertices are in
$\{e_{j_1},\dots,e_{j_i}\}$, then it is in the image of
$\overline{\partial_{i+1}}$ as
$$\overline{\partial_{i+1}}\{e_{j_1},\dots,e_{j_i},e\}=\{e_{j_1},\dots,e_{j_i}\}.$$
If the disjoint edges $\{e_{j_1},\dots,e_{j_i}\}$ form an induced
subgraph of $G$ then $\{e_{j_1},\dots,e_{j_i}\}$ cannot be in the
image of $\overline{\partial_{i+1}}$. If
$\{e_{j_1},\dots,e_{j_i}\}$ were in the image of
$\overline{\partial_{i+1}}$ then we would have
$\overline{\partial_{i+1}}\{e_{j_1},\dots,e_{j_i},e\}=\{e_{j_1},\dots,e_{j_i}\}$
for some $e$ which has at least one vertex not in
$\{e_{j_1},\dots,e_{j_i}\}$ since these edges form an induced
subgraph of $G$. However this implies that
$\frac{{\mathrm{lcm}}\{e_{j_1},\dots ,
e_{j_i},e\}}{{\mathrm{lcm}}\{e_{j_1},\dots e,\dots ,e_{j_i}\}}\in
M$ and so $\overline{\partial_{i+1}}\{e_{j_1},\dots,e_{j_i},e\}$
is zero. Now we see that $\b_{i,2i}= \dim_k H_i(\mathbf{F}_X
\otimes R/M)_{2i}$ is the number of induced subgraphs of $G$ which
are $i$ disjoint edges.

\end{proof}

\chapter{Hochster's Theorem}
\setcounter{satz}{0} In this chapter we look at some results about
Betti numbers and projective dimensions of graph ideals which are
consequences of Hochster's Theorem (\ref{HochstersFormula}).

\section{Betti Numbers of Induced Subgraphs}
We now see that the $i$th Betti number of an induced subgraph
cannot exceed the $i$th Betti number of the larger graph for all
$i$.

\begin{thm}{\Prop}\label{subbet}
If $H$ is an induced subgraph of $G$ on a subset of the vertices
of $G$ then
$$\b_{i,d}^k(H) \leqslant \b_{i,d}^k(G)$$ for all $i$.
\end{thm}

\begin{proof}
Suppose that $W$ is the vertex set of $G$. The vertex set of $H$
is $W \setminus S$ for some $S\subseteq W$. Every subset of $W
\setminus S$ is also a subset of $W$. Hence
\begin{eqnarray*}
\b_{i,d}^k(H) &=& \sum_{V \subseteq(W \ \backslash \ S),|V|=d}
\mathrm{\dim}_k \Hred_{|V|-i-1}(\Delta_V;k)\\
 &\ & \\
 &\leqslant& \ \ \  \sum_{V
\subseteq W, |V|=d} \mathrm{\dim}_k
\Hred_{|V|-i-1}(\Delta_V;k)\\
&=&\b_{i,d}^k(G).
\end{eqnarray*}
\end{proof}

\begin{thm}{\Cor}
If $H$ is an induced subgraph of $G$ then the total Betti numbers
of $H$ do not exceed those of $G$.
\end{thm}

\begin{proof}
This follows from Proposition (\ref{subbet}) as
$$\b^k_i(H)=\sum_{d\in \mathbb{N}} \b^k_{i,d}(H) \leqslant
\sum_{d\in \mathbb{N}} \b^k_{i,d}(G)=\b^k_i(G).$$
\end{proof}

\begin{thm}{\Cor}\label{pdsub}
If $H$ is an induced subgraph of $G$ the $\pd^k(H) \leqslant
\pd^k(G).$
\end{thm}

\begin{proof}
This follows from Proposition (\ref{subbet}) as $\b^k_p(G)
\geqslant \b^k_p(H) >0$ where $p=\pd^k(H)$.
\end{proof}

\section{A Lower Bound for Projective Dimension}

Here we find a lower bound for the projective dimensions of graphs
which depends on certain types of induced subgraph.

\begin{thm}[\rm]{\Def}
For a graph $G$ let $G^c$ denote the (simple, finite) graph on the
same vertices as $G$ with any two vertices joined by an edge if
and only if they are not joined by an edge in $G$.
\end{thm}

For example, if $G$ is the graph on vertices $\{x_1,x_2,x_3,x_4\}$
with edges $\{x_1,x_2\},\{x_1,x_3\}$ and $\{x_3,x_4\}$ then $G^c$
is the graph on the same vertices with edges
$\{x_1,x_4\},\{x_2,x_3\}$ and $\{x_2,x_4\}$.
\begin{figure}[h]
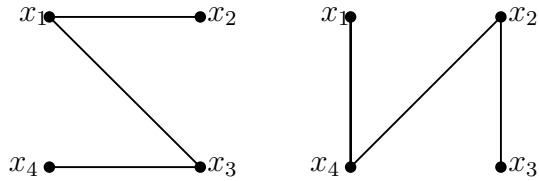
\label{fig3}
\begin{texdraw}
  \drawdim cm
  \textref h:R v:C
  \move(0 0)
  \move(4 0) \fcir f:0 r:0.075  \htext{$x_4 \ $}
  \lvec(6 0) \lvec(4 2) \lvec(6 2)
  \move(6 0) \fcir f:0 r:0.075  \move(6.5 0) \htext{$x_3$}
  \move(6 2) \fcir f:0 r:0.075  \move(6.5 2) \htext{$x_2$}
  \move(4 2) \fcir f:0 r:0.075  \move(4 2) \htext{$x_1$}

  \move(8 0) \fcir f:0 r:0.075  \htext{$x_4 \ $}
  \lvec(8 2) \lvec(8 0) \lvec(10 2) \lvec (10 0)
  \move(10 0) \fcir f:0 r:0.075  \move(10.5 0) \htext{$x_3$}
  \move(10 2) \fcir f:0 r:0.075  \move(10.5 2) \htext{$x_2$}
  \move(8 2) \fcir f:0 r:0.075  \move(8 2) \htext{$x_1$}
  \linewd 0.01 \lpatt(0.2 0.2)
  \lpatt(0.01 100)
\end{texdraw}
\caption{$G$ and $G^c$}
\end{figure}

\begin{thm}{Theorem}\label{pbound}
If $H$ is an induced subgraph of $G$ such that $H^c$ is
disconnected then $$\pd^k(G) \geqslant |V(H)| -1.$$
\end{thm}

\begin{proof}
As $H^c$ is disconnected we can choose a vertex, $x_1$ say, in one
connected component of $H^c$ and another vertex, $x_n$ say, in a
different connected component of $H^c$. Now suppose that
$\D_{V(H)}$ is connected. There must be a path of edges of
$\D_{V(H)}$ from $x_1$ to $x_n$, \ie there must be a sequence of
edges
$$\{x_1,x_2\},\{x_2,x_3\},\dots,\{x_{n-2},x_{n-1}\},\{x_{n-1},x_n\}$$ in $\D_{V(H)}$. This
implies that none of
$\{x_1,x_2\},\{x_2,x_3\},\dots,\{x_{n-1},x_n\}$ are edges of $H$
and so they must all be edges of $H^c$. However this implies that
$x_1$ and $x_n$ are in the same connected component of $H^c$, a
contradiction. Hence $\D_{V(H)}$ is not connected.

From Theorem (\ref{HochstersFormula})
\begin{eqnarray*}
\b_i(G) &=& \sum_{W \subseteq V(G)} {\rm {dim}}_k
\Hred_{|W|-i-1}(\Delta_W;k)\\
&\geqslant&{\rm {dim}}_k
\Hred_{|V(H)|-i-1}(\Delta_{V(H)};k)\\
\end{eqnarray*}
as $V(H) \subseteq V(G)$. When $i=|V(H)|-1$ this becomes
$$\b_{|V(H)|-1}(G) \geqslant {\rm{dim}}_k \Hred_0(\Delta_{V(H)};k) >
0$$ because $\D_{V(H)}$ is disconnected. As the $(|V(H)|-1)$th
Betti number is non-zero we must have
$$\pd(G) \geqslant |V(H)| -1.$$
\end{proof}

\begin{thm}{Theorem}\label{p=v-1}
If $G$ is a graph such that $G^c$ is disconnected then
$$\pd^k(G)=|V(G)|-1.$$
\end{thm}

\begin{proof}
Let $p=\pd(G)$ and let $n=|V(G)|$. By Theorem (\ref{pbound}) $p
\geqslant n -1$, as $G$ is certainly an induced subgraph of
itself. We may use the Hilbert syzygy Theorem (\ref{HST}) to see
that $p\leqslant n$, as $k[\D(G)]$ is a finitely generated graded
$R_k(G)$-module. Therefore the only possible values $p$ may take
are $n$ or $n-1$. Suppose it is the case that $p=n$ and hence that
$\b_p(G) \neq 0$ . Using (\ref{HochstersFormula}) we see that
 $$
 \b_p(G)=\sum_{V \subseteq [n]} {\rm{dim}}_k
\Hred_{|V|-p-1}(\Delta_V;k). $$

For any $V \subseteq [n]$ appearing in the above sum we must have
$|V|\leqslant n$. If $|V| < n$ then $|V|-p-1=|V|-n-1 <-1$. This
implies that $\Hred_{|V|-p-1}(\Delta_V;k)=0$. If, on the other
hand, $|V|=n$ (and hence $\Delta_V=\Delta(G)$) then
$|V|-p-1=|V|-n-1 =-1$. The reduced homology group
$\Hred_{-1}(\Delta;k) \neq 0 $ if and only if $\Delta(G) =
\emptyset$. However, the graph which is such that $\Delta(G) =
\emptyset$ is the graph with no vertices which we do not allow.
Because all the summands in the above formula are zero we must
have $\b_p=0$, a contradiction. Hence $p=n-1.$

\end{proof}

\begin{thm}[\rm]{\Def}\label{Kn}
\begin{enumerate}
\item The Complete Graph, $K_n$, (for $n \geqslant 2$) is the
graph on the $n$ vertices $x_1,\dots,x_n$ which has edges
$\{x_i,x_j\}$ for all $i$ and $j$ such that $1 \leqslant i < j
\leqslant n.$ \item The Complete Bipartite graph, $K_{n,m}$, ($n,m
\geqslant 1$)is the graph which consists of vertices
$x_1,\dots,x_n,y_1,\dots,y_m$ and the edges $\{x_i,y_j\}$ for all
$i$ and $j$ such that $1\leqslant i \leqslant n$ and $1\leqslant j
\leqslant m$.
\end{enumerate}
\end{thm}

\begin{thm}{\Cor}Let $G$ and $H$ be graphs. Let $K_n$ denote the
complete graph on $n$ vertices and let $K_{n,m}$ denote the
complete bipartite graph on $n+m$ vertices.

\begin{enumerate}
\item If $K_{n}$ is an induced subgraph of $G$ then $\pd(G)
\geqslant n-1$. \item If $K_{n,m}$ is an induced subgraph of $H$
then $\pd(H) \geqslant n+m-1.$
\end{enumerate}
\end{thm}

\begin{proof}
(i) The graph $K_n^c$ is disconnected (it is comprised of $n$
isolated vertices) so the result follows from (\ref{pbound}).

(ii) The graph $K_{n,m}$ is also disconnected (it is the graph
with two components, $K_n$ and $K_m$) so again the result follows
from (\ref{pbound})
\end{proof}

\begin{thm}{\Cor}\label{pdcomp}
\begin{enumerate} \item For the complete graph on $n$ vertices, $K_n$, we have
$\pd(K_n)=n-1$. \item For the complete bipartite graph on $n+m$
vertices, $K_{n+m}$, we have $ \pd(K_{n,m})=n+m-1 .$
\end{enumerate}
\end{thm}

\begin{proof} As we observed above $K_n^c$ and $K_{n,m}^c$ are
disconnected. The results follow from (\ref{p=v-1}).
\end{proof}

Notice that the projective dimensions of $K_n$ and $K_{n,m}$ are
independent of our choice of field. We shall see that the Betti
numbers of the complete and complete bipartite graphs are also
independent of the field in the next chapter.

\chapter{Complete and Bipartite Graphs}
\setcounter{satz}{0}

In this chapter we find explicit descriptions of the Betti numbers
of complete graphs, complete bipartite graphs and some related
families of graphs. In each of these cases we will see that the
Betti numbers (and therefore the projective dimensions) do not
depend on our choice of base field of the polynomial ring
$R_k(G)$.
\section{Betti Numbers of Complete Graphs}

\begin{figure}[h]
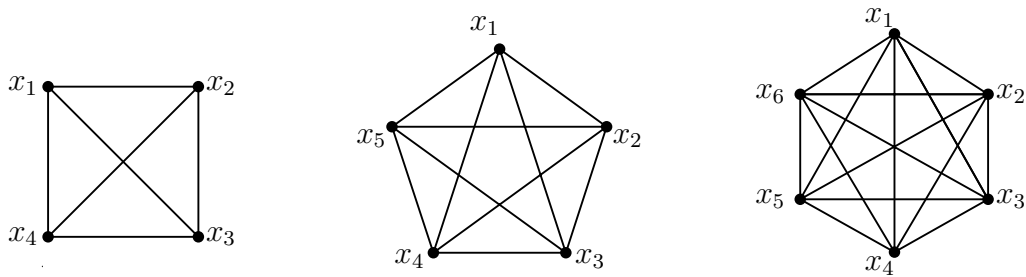
\label{picKn}
\begin{texdraw}
  \drawdim cm
  \textref h:R v:C
  \move(0 0) \fcir f:0 r:0.075  \htext{$x_4 \ $} \lvec(2 0)
  \lvec(2 2) \lvec(0 2) \lvec(0 0) \lvec(2 2)
  \move(0 2) \fcir f:0 r:0.075  \htext{$x_1 \ $} \lvec(2 0)
  \move(7.9 1.3) \htext{$x_2$} \move(7.4 -0.3) \htext{$x_3$} \move(5 -0.3)
  \htext{$x_4$} \move(4.5 1.3)\htext{$x_5$}
  \move(2 0) \fcir f:0 r:0.075  \move(2.5 0) \htext{$x_3$}
  \move(2 2) \fcir f:0 r:0.075  \move(2.5 2) \htext{$x_2$}
  \move(6 2.8) \htext{$x_1$}
  \move (6 2.5) \fcir f:0 r:0.075 \lvec(7.4265844774 1.463525492)
  \lvec(6.881677878 -0.213525491)\lvec(5.118322422
  -0.213525491)\lvec(4.573415226 1.463525492)\lvec(6 2.5)\lvec(6.881677878
  -0.213525491)\lvec(4.573415226 1.463525492)\lvec(7.4265844774
  1.463525492)\lvec(5.118322422 -0.213525491)
  \move (7.4265844774 1.463525492) \fcir f:0 r:0.075
  \move (6.881677878 -0.213525491) \fcir f:0 r:0.075
  \move (5.118322422 -0.213525491) \fcir f:0 r:0.075 \lvec(6 2.5)
  \move (4.573415226 1.463525492)  \fcir f:0 r:0.075
  \move (10 0.5) \fcir f:0 r:0.075 \lvec(10 1.9) \lvec (12.5 1.9)
  \lvec(11.25 -0.2) \lvec(10 1.9)
  \move (10 1.9) \fcir f:0 r:0.075 \lvec (11.25 2.7) \lvec(12.5
  0.5) \lvec (10 0.5) \lvec(11.25 2.7)
  \move (12.5 0.5) \fcir f:0 r:0.075 \lvec (11.25 -0.2)
  \move (12.5 1.9) \fcir f:0 r:0.075 \lvec (12.5 0.5)
  \move (11.25 2.7) \fcir f:0 r:0.075 \lvec (12.5 1.9)
  \move (11.25 -0.2) \fcir f:0 r:0.075 \lvec(10 0.5)\lvec(12.5
  1.9)\lvec(10 1.9) \lvec(12.5 0.5) \lvec(11.25 2.7) \lvec(11.25
  -0.2)
  \move (11.25 2.9) \htext{$x_1$}
  \move (13 1.9) \htext{$x_2$}
  \move (13 0.5) \htext{$x_3$}
  \move (11.25 -0.4) \htext{$x_4$}
  \move (9.8 0.5) \htext{$x_5$}
  \move (9.8 1.9) \htext{$x_6$}
\linewd 0.01 \lpatt(0.2 0.2)
  \lpatt(0.01 100)
\end{texdraw}
\caption{$K_4$, $K_5$ and $K_6$}
\end{figure}

The complete graph $K_n$ is defined in (\ref{Kn}). For a field $k$
let $R=R_k(K_n)=k[x_1,\dots,x_n]$ and
$$I=I(K_n)=\langle x_i x_j \mid 1\leqslant i<j \leqslant n\rangle .$$
By noting that $\Delta
(K_n)=\left\{\left\{1\right\},\dots,\left\{n\right\}\right\}$, the
simplicial complex which is just $n$ isolated  zero dimensional
faces, and applying Theorem (\ref{HochstersFormula}) we can
describe the Betti numbers of $K_n$.

\begin{thm}{Theorem} \label{complete}
The $\mathbb{N}$-graded Betti numbers of the complete graph with
$n$ vertices are independent of the characteristic of $k$ and may
be written
$$\beta_{i,d}(K_n)= \left\{ \begin{array}{lll}
                     i{n \choose {i+1}} &\mathrm{if}& d=i+1\\
                     0                  &\mathrm{if}& d \neq i+1.
                     \end{array}
                     \right.$$

\end{thm}

\begin{proof}From Theorem (\ref{HochstersFormula}) we have

$$\b_{i,d}(K_n) = \sum_{V \subseteq X, |V|=d} {\dm}_k \Hred_{|V|-i-1}(\Delta_V;k),$$
where $X=\lc x_1,\dots,x_n \rc$ and $\Delta=\Delta(K_n) =\lc \lc
x_1 \rc , \dots ,\lc x_n \rc \rc$. For any $V \subseteq X$, such
that $V \neq \emptyset$, $\Delta_V$ is a simplicial complex which
comprises of isolated zero dimensional faces. The only possibly
non-zero reduced homology groups for such simplicial complexes are
those in the 0th position, $\Hred_0(\Delta_V;k)$. Because,
regardless of our choice of the field $k$, ${\dm}_k
\Hred_0(\Delta_V;k)$ is one less than the number of connected
components of $\Delta$ we obtain ${\dm}_k
\Hred_0(\Delta_V;k)=|V|-1.$ So the only contribution to
$\b_{i,d}(K_n)$ is of $|V|-1=d-1$ whenever $d-i-1=0$, \ie when
$d=i+1$. There are ${n\choose d}={ n \choose {i+1}}$ choices for
subsets of $X$ of size $d$, each contributing $d-1=i$ to
$\b_{i,d}(K_n)$. Hence $\b_{i,d}(K_n)= i{n \choose {i+1}}$ when
$d=i+1$ and $\b_{i,d}(K_n)$ is zero otherwise.

\end{proof}

\begin{thm}{\Cor}
The total Betti numbers of the complete graph are
$$\b_i(K_n)=i{n \choose{i+1}}.$$
\end{thm}

\begin{proof}
This follows from (\ref{complete}) by summing all the
$\mathbb{N}$-graded $i$th Betti numbers.
\end{proof}

\begin{thm}[\rm]{Remark}
Note that (\ref{pdcomp}) part $(i)$ is corollary of this theorem
as $\b_{n-1}(K_n)=(n-1){n\choose n}=n-1 \neq 0$ but $\b_n(K_n)=n{
n\choose {n+1}}=0$.
\end{thm}

\section{Betti Numbers of Complete Bipartite Graphs}
We now consider the complete bipartite graphs, as defined in
(\ref{Kn}).

\begin{figure}[h]
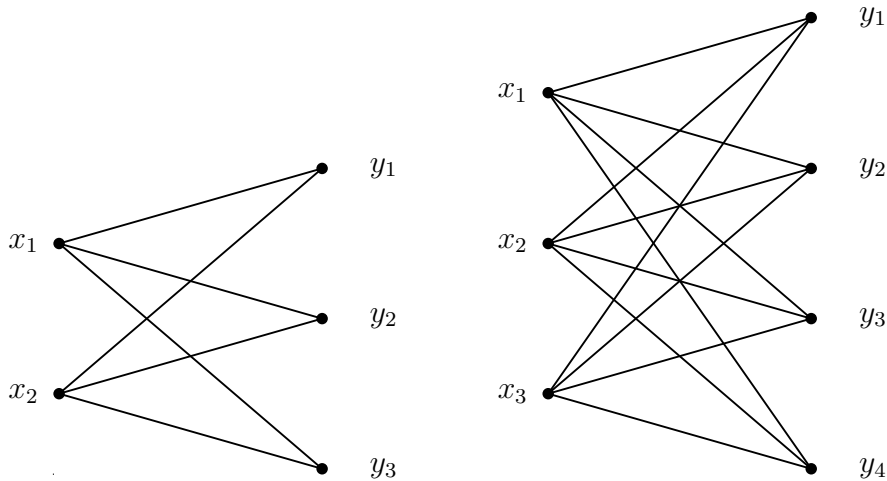
\label{CompBi}
\begin{texdraw}
  \drawdim cm
  \textref h:R v:C
  \move(0 1) \fcir f:0 r:0.075  \htext{$x_2\ \ $} \lvec(3.5 0)
    \move(0 1) \lvec(3.5 4)
  \move(0 3) \fcir f:0 r:0.075  \htext{$x_1\ \ $} \lvec(3.5 2)
  \move(3.5 0) \fcir f:0 r:0.075  \lvec(0 3) \move(4.5 0)\htext{$y_3$}
  \move(3.5 2) \fcir f:0 r:0.075   \lvec(0 1)\move(4.5 2) \htext{$y_2$}
  \move(3.5 4) \fcir f:0 r:0.075   \lvec(0 3)\move(4.5 4)\htext{$y_1$}

  \move(6.5 5) \fcir f:0 r:0.075 \htext{$x_1 \ \ $} \lvec(10 0)
  \move(6.5 5) \lvec(10 2)
  \move(6.5 5) \lvec(10 4)
  \move(6.5 5) \lvec(10 6)

  \move(6.5 1) \fcir f:0 r:0.075  \htext{$x_3\ \ $} \lvec(10 0)
  \move(6.5 1) \lvec(10 4)
  \move(6.5 3) \fcir f:0 r:0.075  \htext{$x_2\ \ $} \lvec(10 2)
  \move(10 0) \fcir f:0 r:0.075  \lvec(6.5 3) \move(11 0)\htext{$y_4$}
  \move(10 2) \fcir f:0 r:0.075   \lvec(6.5 1)\move(11 2) \htext{$y_3$}
  \move(10 4) \fcir f:0 r:0.075   \lvec(6.5 3)\move(11 4)\htext{$y_2$}
  \move(10 6) \fcir f:0 r:0.075 \lvec(6.5 1)
  \move(10 6) \lvec(6.5 3)
  \move(11 6) \htext{$y_1$}

  \linewd 0.01 \lpatt(0.2 0.2)
  \lpatt(0.01 100)
\end{texdraw}
\caption{$K_{2,3}$ and $K_{3,4}$}
\end{figure}

For a field $k$ let $R_k(K_{n,m})=k[x_1,\dots,x_n,y_1,\dots,y_m]$
and
$$I=I({K_{n,m}})=\langle x_iy_j \mid 1\leqslant i\leqslant n, 1
\leqslant j \leqslant m\rangle.$$ The simplicial complex
$\Delta(K_{n,m})$ is the disjoint union of two simplices, one of
dimension $n-1$, the other of dimension $m-1$. As in the case of
the complete graph we use this simple description of the
associated simplicial complex with (\ref{HochstersFormula}) to
find the Betti numbers.

\begin{thm}{Theorem} \label{bipartite}
The $\mathbb{N}$-graded Betti numbers of the complete bipartite
graph with $n+m$ vertices are independent of the characteristic of
$k$ and may be written
$$
\beta_{i,d}(K_{n,m})= \left\{ \begin{array}{lll}
          \sum_{j+l=i+1,\ j,l\geqslant1} {n\choose j} {m \choose
          l}&\mathrm{if}& d=i+1\\
          0 &\mathrm{if}& d \neq i+1.
          \end{array} \right.
$$
\end{thm}

\begin{proof}

Again we use Theorem (\ref{HochstersFormula}) to see that

$$\b_{i,d}(K_{n,m})=\sum_{V \subseteq X \cup Y,\ |V|=d}
{\dm}_k\Hred_{|V|-i-1}(\Delta_V;k)$$

where $X=\lc x_1  ,\dots,  x_n \rc $, $Y=\lc  y_1  ,\dots,  y_m
\rc  $ and $\Delta= \Delta(K_{n,m})$, the disjoint union of an
$(n-1)$-dimensional simplex and an $(m-1)$-dimensional simplex.
Suppose $V \neq \emptyset$. If $V \subseteq X$ or $V \subseteq Y$
then $\Delta_V$ is a simplex and hence has zero reduced homology
everywhere. Consequently we may suppose that $V \cap X \neq
\emptyset$ and $V \cap Y \neq \emptyset$ and therefore $\Delta_V$
is the disjoint union of two simplices. This implies that
$$\Hred_j(\Delta_V;k)=
\left\{
\begin{array}{ll}
k & \ \mathrm{if} \  j = 0\\
0 & \ \mathrm{if} \  j\neq 0.
\end{array}
\right.$$ So we will have a contribution (of 1) to
$\b_{i,d}(K_{n,m})$ if and only if $|V|-i-1=d-i-1$, \ie if
$d=i+1$. Therefore the only non-zero $\mathbb{N}$-graded Betti
numbers are of the form $\b_{i,i+1}$ and to find their value we
must count the subsets of $X \cup Y$, which are disjoint from
neither $X$ nor $Y$, and contain $i+1$ elements. We obtain a
suitable subset by choosing 1 element from $X$ and $i$ from $Y$.
There are ${n \choose 1}{ m \choose i}$ subsets of this type.
Similarly we may choose 2 from $X$ and $i-1$ from $Y$, providing a
further ${n \choose 2}{ m \choose {i-1}}$, and so on to obtain

$${n \choose 1}{m \choose i} + {n \choose 2}{m \choose {i-1}}+
 \dots +{n \choose i}{m \choose 1}= \sum_{j+l=i+1,\ j,l\geqslant1} {n\choose j} {m \choose l}.$$

\end{proof}

\begin{thm}{\Cor}
The total Betti numbers of the complete bipartite graph are
independent of choice of field and are
$$\b_i(K_{n,m})=\sum_{j+l=i+1,\ j,l\geqslant1} {n\choose j} {m \choose
l}.$$
\end{thm}

\begin{proof}
This follows from Theorem (\ref{bipartite}) by summing the graded
Betti numbers.
\end{proof}

\begin{thm}[\rm]{Remark}
Similar to the case for complete graphs, we see that
(\ref{pdcomp}) part $(ii)$ is a corollary of the above;
$\b_{n+m-1}(K_{n,m})={n \choose n}{m \choose m}=1$ and
$\b_{n+m}(K_{n,m})=0$.
\end{thm}

\section{Betti Numbers of Complete Multipartite Graphs}

We may generalize the concept of complete bipartite graphs to
complete multipartite graphs and use similar methods to determine
their Betti numbers.

\begin{thm}[\rm]{\Def}
We define the complete multipartite graph, $K_{n_1,\dots,n_t}$, as
follows. The vertices of  $K_{n_1,\dots,n_t}$ are $x^{(i)}_j$ for
$1 \leqslant i \leqslant t$ and $ 1 \leqslant j \leqslant n_i$.
The edges of $K_{n_1,\dots,n_t}$ are $\{x^{(i)}_j,x^{(k)}_l\}$ for
all $i,j,k,l$ such that $ i \neq k$.

\end{thm}

\begin{figure}[h]\label{Multi}
\begin{texdraw}
  \drawdim cm
  \textref h:R v:C

  \move(1 2) \fcir f:0 r:0.075  \htext{$z_2\ \ $} \lvec(3.5 4.5)
  \lvec(1 3) \lvec(3.5 0.5) \lvec(1 2) \lvec(7 5.5) \lvec (1 3)
  \lvec(7 2.5) \lvec(1 2) \lvec(7 -0.5) \lvec(1 3) \lvec(3.5 4.5)
  \lvec(7 5.5) \lvec(3.5 0.5) \lvec(7 2.5) \lvec(3.5 4.5) \lvec(7
  -0.5) \lvec(3.5 0.5)
  \move(1 3) \fcir f:0 r:0.075  \htext{$z_1\ \ $}
  \move(3.5 0.5) \fcir f:0 r:0.075  \htext{$x_2\ \ $}
  \move(3.5 4.5) \fcir f:0 r:0.075  \htext{$x_1\ \ $}
  \move(7 -0.5) \fcir f:0 r:0.075   \move(8 -0.5)\htext{$y_3$}
  \move(7 2.5) \fcir f:0 r:0.075  \move(8 2.5) \htext{$y_2$}
  \move(7 5.5) \fcir f:0 r:0.075  \move(8 5.5) \htext{$y_1$}
  \linewd 0.01 \lpatt(0.2 0.2)
  \lpatt(0.01 100)
\end{texdraw}
\caption{$K_{2,2,3}$}
\end{figure}

We have $R=R(K_{n_1,\dots,n_t})=k[x^{(i)}_j \mid i=1,\dots,t \ ,
j=1,\dots,n_i]$ and $$I=I(K_{n_1,\dots,n_t})=\langle
x^{(i)}_jx^{(k)}_l \mid i\neq k \rangle.$$ The simplicial complex
$\Delta (K_{n_1,\dots,n_t})$  is defined by its maximal faces
which are
$$\lc x^{(1)}_1,\dots,x^{(1)}_{n_1} \rc, \dots , \lc
x^{(t)}_1,\dots x^{(t)}_{n_t} \rc . $$ \ie
$\Delta(K_{n_1,\dots,n_t})$ is the disjoint union of $t$ simplices
of dimensions $$n_1-1,\dots,n_t-1.$$

\begin{thm}{Theorem}\label{multi}
The $\mathbb{N}$-graded Betti numbers of the complete multipartite
graph $K_{n_1,\dots,n_t}$ are independent of the characteristic of
$k$ and may be written
$$\b_{i,i+1}(K_{n_1,\dots,n_t}) = \sum_{l=2}^{t} (l-1)\sum_{\alpha_1+ \dots +\alpha_l=i+1, \ j_1< \dots <j_l,
  \a_1,\dots,\a_l \geqslant 1} {n_{j_1}
\choose {\alpha_1}} \dots {n_{j_l} \choose {\alpha_l}}.$$ and
$\b_{i,d}(K_{n_1,\dots,n_t})=0$ for $d \neq i+1$.
\end{thm}

\begin{proof}
Let $X^{(i)}= \lc x^{(i)}_j \mid j=1,\dots,n_i \rc$ for $i=1,\dots
t$ and let $X=\bigcup_{i=1}^t X^{(i)}.$ For $\emptyset \neq V
\subseteq X$ the simplicial complex $\Delta_V$ is the disjoint
union of at most $t$ simplices. So, as in the case for bipartite
graphs, the only homology module which may possibly be non-zero is
the 0th one.
$${\dim}_k \Hred_j(\Delta_V;k)=
\left\{
\begin{array}{ll}
(\rm{no.\  of \ connected \ components \ of} \ \Delta_V) -1  & \ if \  j = 0\\
0 & \ if \  j\neq 0.
\end{array}
\right.$$ Using Theorem (\ref{HochstersFormula}) we see that the
only possibly non-zero Betti numbers are those of the form
$\b_{i,i+1}(K_{n_1,\dots,n_t})$. We consider the contribution from
those choices of $V \subseteq X$ for which $\Delta_V$ has $l$
connected components. The $i+1$ elements of $V$ must be chosen
from precisely $l$ of the sets $X^{(1)},\dots X^{(t)}$. Counting
all the possibilities we obtain
$$
\sum_{\alpha_1+ \dots +\alpha_l=i+1, \ j_1< \dots <j_l} {n_{j_1}
\choose {\alpha_1}} \dots {n_{j_l} \choose {\alpha_l}}$$ such
subsets of $X$. Each of these contribute $(l-1)$ to
$\b_i(K_{n_1,\dots,n_t})$. Summing from $l=2$ to $t$ provides the
result.

\end{proof}

\begin{thm}{\Cor} The total Betti numbers of the complete
bipartite graph are independent of choice of field and are
$$\b_i(K_{n_1,\dots,n_t}) = \sum_{l=2}^{t} (l-1)\sum_{\alpha_1+ \dots +\alpha_l=i+1, \ j_1< \dots <j_l} {n_{j_1}
\choose {\alpha_1}} \dots {n_{j_l} \choose {\alpha_l}}.$$
\end{thm}

\begin{proof}
This follows from (\ref{multi}) by summing the graded Betti
numbers.
\end{proof}

\section{Betti Numbers of Star Graphs}
\begin{thm}[\rm]{\Def}
The Star Graph, $S_n$, is the graph with the $n+1$ vertices
$x,y_1,\dots,y_n$ and the $n$ edges $\{x,y_1\},\dots ,\{x,y_n\}$.
\end{thm}

\begin{figure}[h]
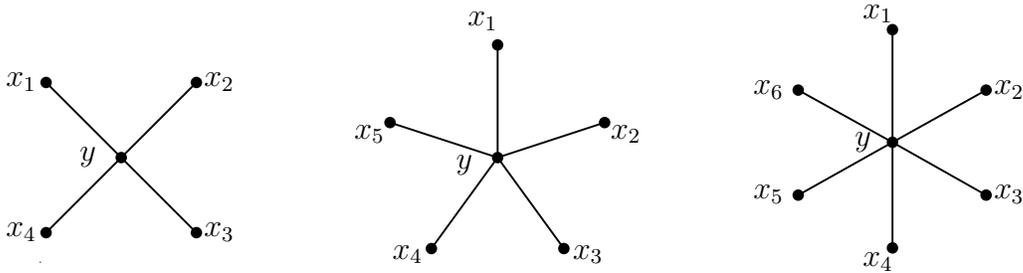
\label{Stars}
\begin{texdraw}
  \drawdim cm
  \textref h:R v:C
  \move(0 0) \fcir f:0 r:0.075  \htext{$x_4 \ $} \lvec(1 1)
  \move(0 2) \fcir f:0 r:0.075  \htext{$x_1 \ $}
  \lvec(1 1)
  \move(1 1) \fcir f:0 r:0.075 \move (.8 1)\htext{$y \ $}
  \move(7.9 1.3) \htext{$x_2$} \move(7.4 -0.3) \htext{$x_3$} \move(5 -0.3)

  \htext{$x_4$} \move(4.5 1.3)\htext{$x_5$}
  \move(2 0) \fcir f:0 r:0.075  \lvec(1 1)\move(2.5 0) \htext{$x_3$}
  \move(2 2) \fcir f:0 r:0.075  \lvec(1 1)\move(2.5 2) \htext{$x_2$}
  \move(6 2.8) \htext{$x_1$}
  \move (6 2.5) \fcir f:0 r:0.075 \lvec(6 1)
  \move (7.4265844774 1.463525492) \fcir f:0 r:0.075 \lvec(6 1)
  \move (6.881677878 -0.213525491) \fcir f:0 r:0.075 \lvec(6 1)
  \move (5.118322422 -0.213525491) \fcir f:0 r:0.075 \lvec(6 1)
  \move (4.573415226 1.463525492)  \fcir f:0 r:0.075 \lvec(6 1)
  \move (6 1) \fcir f:0 r:0.075 \move(5.8 .9)\htext{$y \ $}

  \move (10 0.5) \fcir f:0 r:0.075 \lvec(11.25 1.2)
  \move (10 1.9) \fcir f:0 r:0.075 \lvec(11.25 1.2)
  \move (12.5 0.5) \fcir f:0 r:0.075 \lvec(11.25 1.2)
  \move (12.5 1.9) \fcir f:0 r:0.075 \lvec(11.25 1.2)
  \move (11.25 2.7) \fcir f:0 r:0.075 \lvec(11.25 1.2)
  \move (11.25 -0.2) \fcir f:0 r:0.075 \lvec(11.25 1.2)
  \move (11.25 1.2) \fcir f:0 r:0.075 \move(11.1 1.2)\htext{$y \ $}
  \move (11.25 2.9) \htext{$x_1$}
  \move (13 1.9) \htext{$x_2$}
  \move (13 0.5) \htext{$x_3$}
  \move (11.25 -0.4) \htext{$x_4$}
  \move (9.8 0.5) \htext{$x_5$}
  \move (9.8 1.9) \htext{$x_6$}
  \linewd 0.01 \lpatt(0.2 0.2) \lpatt(0.01 100)

\end{texdraw}
\caption{$S_4$, $S_5$ and $S_6$}
\end{figure}

Let $R=k[x,y_1,\dots,y_n]$ and let $I=I(S_n)=\langle
xy_1,\dots,xy_n \rangle$.

\begin{thm}{Theorem}\label{Star}
The Betti numbers for the star graph with $n$ vertices, $S_n$, are
independent of k and may be expressed as
$$\beta_i(S_n)=
 {n \choose i}. $$It follows that $\pd(S_n)=n.$
\end{thm}

\begin{proof}
The star graph, $S_n$, is in fact just the complete bipartite
graph $K_{1,n}$. By Theorem \ref{bipartite}
$$\b_i(S_n)=\b_i(K_{1,n})=\sum_{j+l=i+1,\ j\geqslant 1, \ l \geqslant 1} {1 \choose j}{n \choose
l}={1 \choose 1}{n \choose i}={n \choose i}
$$
\end{proof}

The formula of Hochster, (\ref{HochstersFormula}), has given
explicit descriptions of Betti numbers of these particular
families of graphs. This is because the simplicial complexes
associated with these graphs have had reduced homology which is
easy to understand and calculate. Typically this is not the case
for a general graph.

\chapter{Links of Subgraphs}
\setcounter{satz}{0} In this chapter we will reinterpret Theorem
(\ref{ER}) for graph ideals which will make use of the structure
of graphs. We first need a useful way to describe the links of
faces of the Alexander Dual. This will be achieved by establishing
a correspondence between induced subgraphs of a graph and the
faces of the Alexander Dual.

\section{Links and Induced Subgraphs}
\begin{thm}[\rm]{\Def}
Let $a_1, \dots , a_s$ be subsets of the finite set $V$. Define
$\e (a_1,\dots,a_s;V)$ to be the simplicial complex which has
vertex set $\bigcup_{i=1}^s (V \setminus a_i)$ and  maximal faces
$V\setminus a_1,\dots, V\setminus a_s$.
\end{thm}

Notice that, with the above notation, the simplicial complex $\e(
\{1,2\};\{1,2\})$ has no vertices and is in fact the complex which
has only one face, namely the empty set (we write this simplicial
complex simply as $\emptyset$).

In what follows $I=I(G)$ will be a graph ideal (for some graph
$G$) of the polynomial ring $R_k(G)=k[x_1, \dots, x_n]$,
$\Delta=\Delta(I)$, the Stanley Reisner complex of $I$. The field
$k$ will remain unspecified throughout this chapter. Although as
noted before Betti numbers will in general be dependent on the
characteristic of $k$, this will not affect the results which
follow.

\begin{thm}[\rm]{Remark} We will sometimes write $\D(G)$ for
$\D(I(G))$ and $\AD(G)$ for $(\D(G))^*=(\D(I(G)))^*.$
\end{thm}

\begin{thm}[\rm]{\Def}
For a face, $F$, of the simplicial complex $\Gamma$ we write
$\Gamma \setminus F$ to denote the subcomplex of $\Gamma$ which
comprises all the faces of $\Gamma$ which do not intersect $F$.
Note that $V(\Gamma \setminus F)=V(\Gamma) \setminus F$ (where
$V(-)$ denotes the vertex set of a simplicial complex).

\end{thm}


\begin{thm}{\Lem}\label{LinkLem1}
Let $G$ be a graph with vertex set $V$ and let $\D=\D(G)$. For a
subset $F$ of $V$ we have $F\in \AD$ if and only if $V\setminus F$
contains an edge of $G$.
\end{thm}

\begin{proof}
This follows from the definition of $\D(G)$ (\ref{SRC}) and the
definition of the Alexander Dual $\AD(G)$ (\ref{addef}). We see
that $F$ is a face of $\AD$ if and only if $V \setminus F$ is not
a face of $\D$ which is in turn true if and only if $V \setminus
F$ contains an edge of $G$.
\end{proof}

\begin{thm}{\Prop}\label{linkepsilon}
Let $F\in \AD=\AD(G)$. Suppose that $e_1,\dots,e_r$ are all the
edges of $G$ which are disjoint from $F$. Then $\linkad
F=\e(e_1,\dots,e_r;W)$ where $W=V(G) \setminus F$.
\end{thm}

\begin{proof}
We first show that the set $(V(G) \setminus F) \setminus e_i$ is a
face of $\linkad F$ for $i=1,\dots,r$, using Lemma
(\ref{LinkLem1}).
\begin{enumerate}
\item $(V(G) \setminus F) \setminus e_i \in \AD$ because $V(G)
\setminus ((V(G) \setminus F) \setminus e_i )=F \cup e_i$ contains
the edge $e_i$ of $G$. \item $((V(G) \setminus F) \setminus e_i)
\cup F \in \AD$ because $V(G) \setminus (((V(G) \setminus F)
\setminus e_i) \cup F)=e_i$ contains an edge of $G$.\item $((V(G)
\setminus F) \setminus e_i) \cap F = \emptyset$
\end{enumerate}
By (\ref{link}) $(V(G) \setminus F) \setminus e_i$ is a face of
$\linkad F$. So all the maximal faces of $\e(e_1,\dots,e_r;V(G)
\setminus F)$ are faces of $\linkad F$. As $\linkad F$ is a
simplicial complex all the faces of $\e(e_1,\dots,e_r;V(G)
\setminus F)$ must be faces of $\linkad F$.
Therefore $\e(e_1,\dots,e_r;V(G) \setminus F) \subseteq  \linkad
F$.

Now suppose that $J \in \linkad F$. Therefore $V(G) \setminus (F
\cup J)$ contains an edge of $G$ (which is disjoint from $F$), \ie
it contains $e_i$ for some $i\in \{1,\dots,r\}$. Also $F \cap J =
\emptyset$, so we see that $J \subseteq (V(G) \setminus F)
\setminus e_i$ which is a face of $\e(e_1,\dots,e_r;V(G) \setminus
F)$. Hence $J \in \e(e_1,\dots,e_r; V(G) \setminus F)$ and
therefore $\linkad F \subseteq \e(e_1,\dots,e_r; V(G) \setminus
F)$ and we conclude $$\linkad F = \e(e_1,\dots,e_r; V(G) \setminus
F).$$
\end{proof}

We now show how we may associate faces of the Alexander Dual with
subgraphs of $G$.

\begin{thm}{\Prop}
There is a bijection between the faces of $\AD(G)$ and the set of
induced subgraphs of $G$ which have at least one edge.

$$F\in \AD(G) \xrightarrow{\varphi} The \ induced \ subgraph \ of\ G\ on\ vertices \ V(G) \setminus
F.
$$

\end{thm}

\begin{proof}
Suppose that $H$ is an induced subgraph of $G$ and that $e$ is an
edge of $H$ ( we are assuming that $H$ has at least one edge). Let
$F= V(G) \setminus V(H)$. Notice that $F\in \AD(G)$ because $V(G)
\setminus F =V(G) \setminus (V(G) \setminus V(H))=V(H)$ contains
the edge $e$ of $G$. Now we observe that $H= \varphi(F)$.
Therefore $\varphi$ is surjective.

Now suppose that $F \in \AD(G)$ and $F' \in \AD$ are such that
$\varphi(F)=\varphi(F')=H$. Thus $H$ must be the induced subgraph
of $G$ on vertices $V(G) \setminus F=V(G) \setminus F'$. Hence
$F=F'$ and $\varphi$ is injective.
\end{proof}

For example, let $G$ be the following graph,
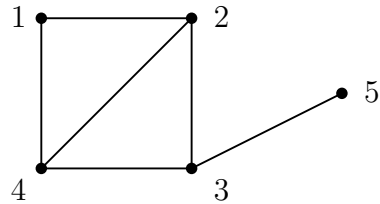
\begin{figure}[h]\label{fig4.1}
\begin{texdraw}
  \drawdim cm
  \textref h:R v:C
  \move(0 0)
  \move(6 0) \fcir f:0 r:0.075
  \lvec(8 2) \lvec(6 2) \lvec(6 0) \lvec (8 0)\lvec(8 2)
  \move(8 0) \fcir f:0 r:0.075  \move(8.5 -.3) \htext{$3$}
  \move(8 2) \fcir f:0 r:0.075  \move(8.5 2) \htext{$2$}
  \move(6 2) \fcir f:0 r:0.075  \move(5.8 2) \htext{$1$}
  \move(8 0) \lvec (10 1) \fcir f:0 r:0.075 \move (10.5 1)
  \htext{$5$} \move (5.8 -.3) \htext{$4$}
\linewd 0.01 \lpatt(0.2 0.2)
  \lpatt(0.01 100)
\end{texdraw}
\end{figure}

If we now consider the face $F=\{3,5\}$ of the Alexander Dual
(this is so as $\{1,2,3,4,5\} \setminus F $ contains an edge of
$G$). Therefore $F$ corresponds to the induced subgraph of $G$ on
vertices
$$\{1,2,3,4,5\} \setminus \{3,5\} =\{1,2,4\}.$$
\begin{figure}[h]
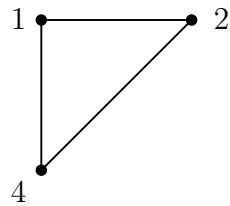
\label{fig4.2}
\begin{texdraw}
  \drawdim cm
  \textref h:R v:C
  \move(0 0)
  \move(6 0) \fcir f:0 r:0.075
  \lvec(8 2) \lvec(6 2) \lvec(6 0)
  \move(8 2) \fcir f:0 r:0.075  \move(8.5 2) \htext{$2$}
  \move(6 2) \fcir f:0 r:0.075  \move(5.8 2) \htext{$1$}
  \move (5.8 -.3) \htext{$4$} \linewd 0.01 \lpatt(0.2 0.2)
  \lpatt(0.01 100)
\end{texdraw}
\caption{Induced subgraph of $G$ on vertices 1,2 and 4.}
\end{figure}

We can write the link of $F$ as
$$\Link_{\AD(G)} F=\e( \{1,2\},\{1,4\},\{2,4\} ;\{1,2,4\}) $$
         (this is the simplicial complex which consists of three isolated
         vertices).

\begin{thm}[\rm]{\Def}
Let $H$ be an induced subgraph of the graph $G$. If $H$ is
associated to the face $F$ of the Alexander Dual as described
above we write $\e(H)$ for $\e(e_1,\dots,e_s;V)$, where
$e_1,\dots,e_s$ are the edges of $H$ and $V$ is the vertex set of
$V(G) \setminus F$ (or equivalently the vertex set of $H$).
\end{thm}

\begin{thm}{Theorem}\label{ER2}
For a graph G the $\mathbb{N}$-graded Betti numbers are
$$\b_{i,d}(G)=\sum _{H \subset G, |V(H)|=d} {\dim}_k
\Hred_{i-2}(\e(H);k)$$and the total Betti numbers are
$$\b_{i}(G)=\sum _{H \subset G} {\dim}_k
\Hred_{i-2}(\e(H);k)$$where the sum is over all induced subgraphs
of $G$ which contain at least one edge.
\end{thm}

\begin{proof}
This is just Theorem (\ref{ER}) with the $\linkad F$ replaced by
$\e(H)$ where $H$ is the subgraph of $G$ associated to the face
$F$ of $\AD$ as described above. Note that $|F|=n-d \iff
|V(H)|=n-|F|=d$.
\end{proof}

\section{Reduced Homology of Links}
We now prove a few lemmas about the reduced homology of the type
of simplicial complexes we will use later.

\begin{thm}{\Lem}\label{lemA}
Let $e_1,\dots,e_t$ be subsets of the finite set $V$. Let
$E=\e(e_1,\dots,e_t;V)$. If $V \setminus \bigcup_{i=1}^t e_i \neq
\emptyset$, \ie if there is an element of $V$ which is in none of
$e_1,\dots,e_t$, then $\Hred_i(E)=0$ for all $i$.
\end{thm}

\begin{proof}
If $V \setminus \bigcup_{i=1}^t e_i \neq \emptyset$ then there
exists some vertex common to all maximal faces of $E$, \ie $E$ is
a cone and so has zero reduced homology everywhere.
\end{proof}

\begin{thm}{\Cor}\label{isolink}
Let $H$ be a graph with at least one edge and at least one
isolated vertex. Then $\dim_k \Hred_i (\e(H);k)=0$ for all i.
\end{thm}

\begin{proof}
Let $a$ be an isolated vertex of $H$ and let $e_1,\dots,e_r$ be
the edges of $H$. Then $\e(H)=\e(e_1,\dots,e_r;V(H))$. As $a\in
V(H) \setminus \bigcup_{j=1}^r e_j$ we have $\dim_k \Hred_i
(\e(H);k)=0$ for all $i$ by Lemma (\ref{lemA}).
\end{proof}

\begin{thm}{\Lem}\label{lemB}
Suppose the simplicial complexes $E_1$ and $E_2$ are defined as
follows. Let $E_1=\e(e_1\dots,e_t;V)$ and let $E_2=\e(f;V)$. The
intersection may be expressed as $E_1 \cap E_2 =\e(f \cup
e_1,\dots,f\cup e_t;V )$. If $f \cap \bigcup_{i=1}^t e_i =
\emptyset$ then $E_1 \cap E_2 =\e(e_1,\dots,e_t;V \setminus f )$.
\end{thm}

\begin{proof}
The maximal faces of $E_1 \cap E_2$ are the intersections of the
maximal faces of $E_1$, that is $V \setminus e_1, \dots V
\setminus e_t$, with the maximal face of $E_2$, $V \setminus f$.
These are the sets $V \setminus (e_i \cup f)$ for $i=1,\dots,t$.
Hence we may write $$E_1 \cap E_2 =\e(f \cup e_1,\dots,f\cup e_t;V
).$$

If $f \cap \bigcup_{i=1}^t e_i = \emptyset$ then we have $$E_1
\cap E_2= \e((e_1 \cup f) \setminus f,\dots,(e_t \cup f) \setminus
f;V \setminus f)=\e(e_1,\dots,e_t;V \setminus f ).$$
\end{proof}

\begin{thm}{\Lem}\label{lemC}
If $f \subseteq g$ then
$\e(f,g,e_1,\dots,e_t;V)=\e(f,e_1,\dots,e_t;V).$
\end{thm}

\begin{proof}
The simplicial complex $\e(f,g,e_1,\dots,e_t;V)$ has faces $$V
\setminus f, V \setminus g, V \setminus e_1,\dots, V \setminus
e_t$$ and all of their subsets. If $f \subseteq g$ then $V
\setminus g$ is a subset of $V \setminus f$, as are all the
subsets of $V \setminus g$, so we may write the simplicial complex
as $\e(f,e_1,\dots,e_t;V).$
\end{proof}

\begin{thm}[\rm]{Remark}
A tool which will be of much use is the Mayer-Vietoris sequence,
see \cite{Ha} chapter $2$ or \cite{S} chapter $4$. This is the
long exact sequence of reduced homology modules
\begin{eqnarray*}
\dots \to \Hred_i(\D_1 \cap \D_2;k) \to \Hred_{i}(\D_1;k) \oplus
\Hred_i(\D_2;k)\to \Hred_i(\D;k)&\to&\\
\to \Hred_{i-1}(\D_1 \cap \D_2;k) \to \Hred_{i-1}(\D_1;k) \oplus
\Hred_{i-1}(\D_2;k)\to \Hred_{i-1}(\D;k)\to\dots &\ & \
\end{eqnarray*}
where the simplicial complex $\D$ is the union of the subcomplexes
$\D_1$ and $\D_2$.

\end{thm}

\begin{thm}{\Lem}\label{lemD}
We define the simplicial complex, $E$, as follows. Let $E=\e(
\{a\},e_1,\dots,e_t;V)$ where $\{a\},e_1,\dots,e_t$ are subsets of
the set $V$. If $ a \notin \bigcup_{i=1}^t e_i$ then
$$\Hred_i(E)=\Hred_{i-1}(\e(e_1,\dots,e_t;V \setminus \{a\}))$$
for all $i$.
\end{thm}

\begin{proof}
Let $\e_1=\e(\{a\};V)$ and let $\e_2=\e(e_1,\dots,e_t;V)$. It is
easily seen that $E=\e_1 \cup \e_2$. The simplicial complex $\e_1$
is in fact just a simplex (the $|V|-2$ dimensional simplex on the
vertices $V \setminus \{a\}$) and so has zero reduced homology
everywhere. Also we note that $a \in V \setminus \bigcup_{i=1}^t
e_i$ so by Lemma (\ref{lemA}) $\Hred_i(\e_2)=0$ for all $i$. We
now make use of the Mayer-Vietoris sequence
$$\dots \to \Hred_{i}(\e_1) \oplus \Hred_i(\e_2)\to \Hred_i(E)\to \Hred_{i-1}(\e_1 \cap \e_2)\to \dots$$
which becomes
$$\dots \to 0 \to \Hred_i(E) \to \Hred_{i-1}(\e_1 \cap \e_2)\to
0\to \dots$$ and therefore we have the isomorphism of reduced
homology modules $\Hred_i(E) \cong \Hred_{i-1}(\e_1 \cap \e_2)$
for all $i$. By Lemma (\ref{lemB}) $$\e_1 \cap \e_2
=\e(e_1\cup\{a\},\dots,e_t\cup\{a\};V)=\e(e_1,\dots,e_t;V
\setminus \{a\})$$ as $a \notin \bigcup_{i=1}^t e_i.$ Therefore we
conclude that $$\Hred_i(E) \cong \Hred_{i-1}(\e(e_1,\dots,e_t;V
\setminus \{a\} ))$$ for all $i$.
\end{proof}

\begin{thm}{\Cor}\label{CorE}
Let $a_1,\dots,a_s$ be $s$ distinct elements of $V$ and let
$E=\e(\{a_1\},\dots,\{a_s\},e_1,\dots,e_t;V)$. If
$$\{a_1,\dots,a_s \} \cap \bigcup_{i=1}^t e_i = \emptyset$$ then
$\Hred_i(E)=\Hred_{i-s}(\e(e_1,\dots,e_t;V'))$, where $V'=V
\setminus \{a_1,\dots,a_s \}$, for all $i$.
\end{thm}

\begin{proof}
Because $a_j \notin (\bigcup_{i=j+1}^s \{a_i\}) \cup
(\bigcup_{i=1}^t e_i)$ for $j=1,\dots, s$ we may repeatedly apply
Lemma (\ref{lemD}) to obtain
\begin{eqnarray*}
\Hred_i(E)&=&\Hred_{i-1}(\e(\{a_2\},\dots,\{a_s\},e_1,\dots,e_t;V
\setminus \{a_1 \}))\\
&=&\Hred_{i-2}(\e(\{a_3\},\dots,\{a_s\},e_1,\dots,e_t;V
\setminus \{a_1,a_2 \}))\\
&=& \dots \\
&=&\Hred_{i-s}(\e(e_1,\dots,e_t;V')).
\end{eqnarray*}

\end{proof}

In the chapters which follow we will make use of these facts about
the reduced homology of links to investigate Betti numbers and
projective dimensions of cycles and trees.
\chapter{Betti Numbers of Cycles}
\setcounter{satz}{0}
\begin{figure}[h]
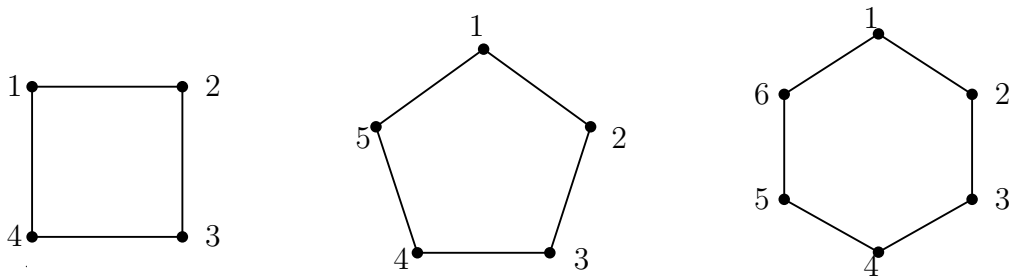
\label{cycle}
\begin{texdraw}
  \drawdim cm
  \textref h:R v:C
  \move(0 0) \fcir f:0 r:0.075  \htext{$4 \ $} \lvec(2 0)
  \lvec(2 2) \lvec(0 2)
  \move(0 2) \fcir f:0 r:0.075  \htext{$1 \ $} \lvec(0 0)
  \move(7.9 1.3) \htext{$2$} \move(7.4 -0.3) \htext{$3$} \move(5 -0.3)
  \htext{$4$} \move(4.5 1.3)\htext{$5$}
  \move(2 0) \fcir f:0 r:0.075  \move(2.5 0) \htext{$3$}
  \move(2 2) \fcir f:0 r:0.075  \move(2.5 2) \htext{$2$}
  \move(6 2.8) \htext{$1$}
  \move (6 2.5) \fcir f:0 r:0.075 \lvec(7.4265844774 1.463525492)
  \move (7.4265844774 1.463525492) \fcir f:0 r:0.075 \lvec(6.881677878 -0.213525491)
  \move (6.881677878 -0.213525491) \fcir f:0 r:0.075 \lvec(5.118322422 -0.213525491)
  \move (5.118322422 -0.213525491) \fcir f:0 r:0.075 \lvec(4.573415226 1.463525492)
  \move (4.573415226 1.463525492)  \fcir f:0 r:0.075 \lvec(6 2.5)

  \move (10 0.5) \fcir f:0 r:0.075 \lvec(10 1.9)
  \move (10 1.9) \fcir f:0 r:0.075 \lvec (11.25 2.7)
  \move (12.5 0.5) \fcir f:0 r:0.075 \lvec (11.25 -0.2)
  \move (12.5 1.9) \fcir f:0 r:0.075 \lvec (12.5 0.5)
  \move (11.25 2.7) \fcir f:0 r:0.075 \lvec (12.5 1.9)
  \move (11.25 -0.2) \fcir f:0 r:0.075 \lvec(10 0.5)
  \move (11.25 2.9) \htext{$1$}
  \move (13 1.9) \htext{$2$}
  \move (13 0.5) \htext{$3$}
  \move (11.25 -0.4) \htext{$4$}
  \move (9.8 0.5) \htext{$5$}
  \move (9.8 1.9) \htext{$6$}
  \linewd 0.01 \lpatt(0.2 0.2) \lpatt(0.01 100)

\end{texdraw}
\caption{$C_4$, $C_5$ and $C_6$}
\end{figure}
Throughout this chapter $C_n$ will denote the cycle graph on $n$
vertices. We fix a field $k$ and let $R=R_k(C_n)=k[x_1,\dots,x_n]$
and let $I=I(C_n)=\langle x_1x_2,x_2x_3,\dots,x_{n-1}x_n,x_nx_1
\rangle.$ To calculate the Betti numbers of $I(C_n)$ we will use
the formula given in (\ref{ER}) in the form from (\ref{ER2})
$$\b_i^k(C_n)= \sum_{H \subset C_n} {\dim}_k \Hred_{i-2}(\e(H);k).$$
To find a more explicit description of the Betti numbers of $C_n$
we must consider the induced subgraphs of $C_n$.

\begin{thm}[\rm]{Remark}
The use of Theorem (\ref{HochstersFormula}) to calculate Betti
numbers of complete, complete multipartite and star graphs relied
on their associated simplicial complexes being simple enough to
easily understand their reduced homology. The simplicial complex
$\Delta(C_n)$ consists of subsets of the vertices of $C_n$ such
that no two vertices are adjacent ( $1$ and $n$ are also
considered to be adjacent). The reduced homology of this complex
is not so easily comprehended so we approach the problem using the
alternative formula.
\end{thm}

\section{Induced Subgraphs of Cycles}

\begin{thm}[\rm]{\Def}
Let $W \subsetneq V= \{ 1, \dots, n \}$. The induced subgraph of
$C_n$ on the elements of $W$ will consist of one or more connected
components which are straight line graphs. We will call such a
connected subgraph which has $m$ vertices a {\it run of length
$m$}.
\end{thm}

For example, if $n=8$ and $W=\{ 1,2,4,5,6 \} \subseteq
\{1,2,3,4,5,6,7,8 \}=V(C_8)$, then the induced subgraph of $C_8$
on $W$ consists of one run of length 2 and one run of length 3.
\begin{figure}[h]
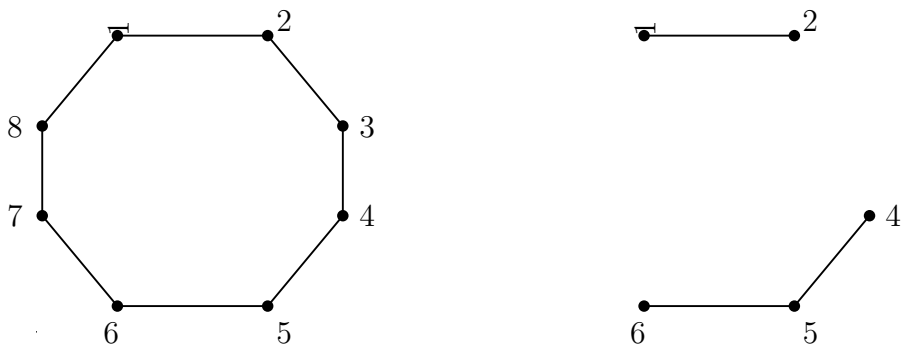
\label{cycle}
\begin{texdraw}
  \drawdim cm
  \textref h:R v:C
  \move(0 0) \fcir f:0 r:0.075 \move(0 -0.35)\htext{$6$} \move (2 0)  \fcir f:0 r:0.075 \move(2.3 -0.35)\htext{5}
  \move(0 3.6)\fcir f:0 r:0.075 \move(0 3.8) \vtext{$1$}\move (2 3.6)  \fcir f:0 r:0.075 \move(2.3 3.8)\htext{2}
  \move(3 1.2) \fcir f:0 r:0.075 \move(3.4 1.2)\htext{4}\move (3 2.4)  \fcir f:0 r:0.075 \move(3.4 2.4)\htext{3}
  \move(-1 2.4) \fcir f:0 r:0.075 \htext{8 \ }\move (-1 1.2)  \fcir f:0 r:0.075 \htext{7 \ }
  \lvec (-1 2.4) \lvec(0 3.6) \lvec(2 3.6) \lvec(3 2.4) \lvec(3 1.2) \lvec(2 0) \lvec(0 0)\lvec(-1 1.2)

  \move(7 0) \fcir f:0 r:0.075 \move(7 -0.35)\htext{$6$}\move (9 0)  \fcir f:0 r:0.075 \move(9.3 -0.35)\htext{5}
  \move(7 3.6) \fcir f:0 r:0.075 \move(7 3.8) \vtext{$1$}\move(7 3.6)\lvec (9 3.6)\move (9 3.6)  \fcir f:0 r:0.075
  \move(9.3 3.8)\htext{2}
  \move(10 1.2) \fcir f:0 r:0.075 \move(10.4 1.2)\htext{4}
\move(10 1.2)  \lvec(9 0) \lvec(7 0)

  \linewd 0.01 \lpatt(0.2 0.2) \lpatt(0.01 100)

\end{texdraw}
\caption{$C_8$ and an Induced Subgraph.}
\end{figure}

\begin{thm}[\rm]{Remark}
All the proper (\ie all the induced subgraphs except $C_n$ itself)
induced subgraphs of $C_n$ are composed of a collection of runs.
If such a subgraph, $H$, has any isolated vertices then by
(\ref{isolink}) $\e(H)$ has zero homology everywhere. Hence we
need only consider induced subgraphs with at least one edge and no
isolated vertices.

Suppose that two such induced subgraphs, $H$ and $G$ say, are both
composed of runs of length $s_1, \dots, s_r$, \ie they are
isomorphic graphs. The simplicial complexes $\e(H)$ and $\e(G)$
are homeomorphic and so have the same reduced homology. As we can
write the link of any face of the Alexander Dual as $\e(H)$ for
some induced subgraph of $C_n$ by Proposition (\ref{linkepsilon}),
all the information we require about the links of faces of the
Alexander dual is dependent on the number of runs and their
lengths of the associated subgraph of $C_n$. We may make the
following definition
\end{thm}

\begin{thm}[\rm]{\Def}
Suppose that the induced subgraph $H$ of $C_n$, which has edges
$e_1,\dots,e_t$ and vertex set $V$, is comprised of $r$ runs which
are of lengths $s_1,\dots,s_r$. We will write
$E(s_1,\dots,s_r)=\e(e_1,\dots,e_t;V)$. So for any $\emptyset \neq
F \in \AD(C_n)$ we  may write $\Link_{\AD(C_n)} F=
E(s_1,\dots,s_r)$ for some $s_1,\dots, s_r$.
\end{thm}

\section{Homology of Links}

In this section we show how to find the reduced homology of links
associated to subgraphs of cycles. We will see that these reduced
homology modules are a shift of the reduced homology modules of
links of `smaller' subgraphs of $C_n$.

\begin{thm}[\rm]{Remark}
In what follows, for a simplicial complex $\G$, will we often
simply write $\Hred_i(\G)$ rather than $\Hred_i(\G;k)$. None of
the results in this chapter which follow will depend on the
characteristic of $k$.
\end{thm}

\begin{thm}{Theorem}\label{homshift}
For $F \in \AD(C_n)$ , with $F\neq \emptyset$, such that $\linkad
F=E(s_1,\dots,s_r)$ with $r\geqslant1$ and $s_r\geqslant 5$, we
have the isomorphism of reduced homology modules $\Hred_{i}
(E(s_1,\dots,s_r)) \cong \Hred_{i-2}(E(s_1,\dots,s_r-3))$ for all
$i$.
\end{thm}

\begin{proof}
Suppose the edges of $G$ which are disjoint from $F$ are
$e_1,\dots,e_t$, ordered such that $e_1,\dots, e_{s_r-1}$ are the
edges in the run of length $s_r$. Let
$E=E(s_1,\dots,s_r)=\e(e_1,\dots,e_t;V)$. For convenience we label
the first five vertices in the run $1,2,3,4,5$ and so write
$e_1=\left\{1,2\right\}$, $e_2=\left\{2,3\right\}$,
$e_3=\left\{3,4\right\}$ and $e_4=\left\{4,5\right\}$. We can
write $E=E_1 \cup E_2$ where $E_1=\e(e_1;V)$ and
$E_2=\e(e_2,\dots,e_t;V)$. We now use the Mayer-Vietoris sequence
$$\dots\to\Hred_i(E_1)\oplus\Hred_i(E_2)\to\Hred_i(E)\to\Hred_{i-1}(E_1\cap
E_2)\to\Hred_{i-1}(E_1)\oplus\Hred_{i-1}(E_2)\to\dots$$ We know
that $\Hred_{i}(E_1)=0$ for all $i$ because $E_1$ is a simplex.
Also because $ 1 \in V \setminus \bigcup_{i=2}^t e_i$ we can use
Lemma (\ref{lemA}) to see that $\Hred_{i}(E_2)=0$ for all $i$.
Hence the Mayer-Vietoris sequence becomes
$$\dots\to 0 \to \Hred_i(E) \to \Hred_{i-1}(E_1\cap
E_2)\to 0 \to\dots$$ which implies $\Hred_i(E)\cong\Hred_{i-1}(E_1
\cap E_2)$ for all $i$. We now examine the simplicial complex $E_1
\cap E_2$. By Lemma (\ref{lemB})
\begin{eqnarray*}
E_1 \cap E_2 &=&\e( e_1 \cup e_2,e_1 \cup e_3,\dots,e_1 \cup
e_t;V)\\
&=&
\e(\left\{1,2,3\right\},\left\{1,2,3,4\right\},\left\{1,2,4,5\right\},\dots,\left\{1,2\right\}\cup
e_t;V)\\
&=& \e(\left\{3\right\},\left\{3,4\right\},\left\{4,5\right\},\dots,e_t;V')\\
&\ & (\rm{using \ Lemma} \ (\ref{lemC}) \ \rm{as} \  \{3\} \subseteq\{3,4\} ) \\
&=& \e(\left\{3\right\},\left\{4,5\right\},\dots,e_t;V')\\
&=& \e(\left\{3\right\},e_4,\dots,e_t;V')
\end{eqnarray*}
where $V'=V \setminus e_1$. Because $3 \notin \bigcup_{i=4}^t e_i$
we use Lemma (\ref{lemD}) to see that
$\Hred_i(\e(\left\{3\right\},e_4,\dots,e_t;V'))=\Hred_{i-1}(\e(e_4,\dots,e_t;V''
))$ where $V''=V \setminus  \{ 1,2,3 \}$. Putting these results
together
\begin{eqnarray*}
\Hred_i(E)&\cong&\Hred_{i-1}(\e(\left\{3\right\},e_4,\dots,e_t;V'))\\
          &\cong&\Hred_{i-2}(\e(e_4,\dots,e_t;V''))\\
          &\cong&\Hred_{i-2}(E(s_1,\dots,s_r-3)).
\end{eqnarray*}
\end{proof}

\begin{thm}[\rm]{Remark} There is no ordering of $s_1,\dots,s_r$
in $E(s_1,\dots,s_r)$, \ie
$$E(s_1,\dots,s_r)=E(\sigma(s_1),\dots,\sigma(s_r))$$ for any
permutation, $\sigma$, of $\{s_1,\dots,s_r\}$. So long as $r
\geqslant 1$ and $s_j \geqslant 5$, with $ 1\leqslant j \leqslant
r$, then we have an isomorphism of reduced homology groups
$\Hred_{i} (E(s_1,\dots,s_j,\dots,s_r)) \cong
\Hred_{i-2}(E(s_1,\dots,s_j-3,\dots,s_r))$ for all $i$.

\end{thm}

We now consider what happens when we have runs of lengths less
than 5.
\begin{thm}{\Prop}\label{homshift2}
Let $E=E(s_1,\dots,s_r)$.\\
(i) If $s_r=2$ and $r\geqslant 2$ then
$\Hred_i(E(s_1,\dots,s_r))=\Hred_{i-1}(E(s_1,\dots,s_{r-1}))$.\\
(ii) If $s_r=3$ and $r\geqslant 2$ then
$\Hred_i(E(s_1,\dots,s_r))=\Hred_{i-2}(E(s_1,\dots,s_{r-1}))$.\\
(iii) If $s_r=4$ and $r\geqslant 1$ then
$\Hred_i(E(s_1,\dots,s_r))=0$.
\end{thm}

\begin{proof}

$(i)$ Let
$E=E(s_1,\dots,s_{r-1},2)=\e(\left\{1,2\right\},e_2,\dots,e_t;V)$
 for some edges $e_2,\dots,e_t$ and where $\left\{1,2\right\}$ is the run of length 2. Let
$E_1=\e(\left\{1,2\right\};V)$ and $E_2=\e(e_2,\dots,e_t;V)$ so
that $E=E_1 \cup E_2$. We have $\Hred_{i}(E_1)=0$ for all $i$
because $E_1$ is a simplex, and $\Hred_i(E_2)=0$ for all $i$ by
Lemma (\ref{lemA}) as $1,2 \in V \setminus \bigcup_{i=2}^t e_i$.
Use of the Mayer-Vietoris sequence implies that $\Hred
_i(E)=\Hred_{i-1}(E_1\cap E_2)$ for all $i$. We now see that
\begin{eqnarray*}
E_1 \cap E_2 &=&\e(\left\{1,2\right\} \cup e_2,\dots,
\left\{1,2\right\} \cup e_t;V) \\
&\ & (\rm{using} \  \rm{Lemma}\ (\ref{lemB}))\\
 &=& \e(e_2,\dots,e_t;V')
\end{eqnarray*}
where $V'=V \setminus \left\{1,2\right\}$. Hence
$\Hred_i(E(s_1,\dots,s_r))=\Hred_{i-1}(\e(e_2,\dots,e_t;V'))$ for
all $i$. The simplicial complex
$\e(e_2,\dots,e_t;V')=E(s_1,\dots,s_{r-1})$, as this is the
complex associated with the induced subgraph of $C_n$ which is
composed of runs of lengths $s_1,\dots,s_{r-1}$. Hence
$\Hred_i(E(s_1,\dots,s_r))=\Hred_{i-1}(E(s_1,\dots,s_{r-1}))$ for
all $i$.

$(ii)$ Now let
$E=E(s_1,\dots,s_{r-1},3)=\e(\left\{1,2\right\},\left\{2,3\right\},e_3,\dots,e_t;V)$
where $\left\{1,2\right\},\left\{2,3\right\}$ is a run of length 3
and so none of 1,2 or 3 is in any of $e_3,\dots,e_t$. Let
$E_1=\e(\left\{1,2\right\};V)$ and
$E_2=\e(\left\{2,3\right\},e_3,\dots,e_t;V)$. As before
$\Hred_{i}(E_1)=\Hred_{i}(E_2)=0$ for all $i$ ($E_1$ is a simplex
and $1 \in  V \setminus (\{2,3\} \cup \bigcup_{i=3}^t e_i)$ so
Lemma (\ref{lemA}) applies). The Mayer-Vietoris sequence implies
$\Hred _i(E) \cong \Hred_{i-1}(E_1\cap E_2)$, for all $i$. We have
\begin{eqnarray*}
E_1 \cap E_2 &=& \e(\left\{1,2,3\right\},\left\{1,2\right\} \cup
e_3,\dots,\left\{1,2\right\} \cup e_t;V)\\
&=& \e(\left\{3\right\},e_3,\dots,e_t;V')
\end{eqnarray*}
where $V'= V \setminus \left\{1,2\right\}.$ As $3 \notin
\bigcup_{i=3}^t e_i$ we can make use of Lemma (\ref{lemD}) to
obtain, for all $i$,
\begin{eqnarray*}
\Hred_i(E_1 \cap E_2 )&=& \Hred_{i}(\e(\left\{3\right\},e_3,\dots,
e_t;V')) \\
&=& \Hred_{i-1}( \e(e_3,\dots,e_t;V''))
\end{eqnarray*}
where $V''=V' \setminus \left\{3\right\}$. We note that
$\e(e_3,\dots,e_t;V'')=E(s_1,\dots,s_{r-1})$, the simplicial
complex associated with the induced subgraph of $C_n$ which is
$r-1$ runs of lengths $s_1,\dots,s_{r-1}$. Hence
$$\Hred_i(E(s_1,\dots,s_{r-1},3))=\Hred_{i-2}(\e(e_3,\dots,e_t;V''))=\Hred_{i-2}(E(s_1,\dots,s_{r-1})).$$

$(iii)$ Finally suppose that
$$E=E(s_1,\dots,s_{r-1},4)=\e(\left\{1,2\right\},\left\{2,3\right\},\left\{3,4\right\},e_4,\dots,e_t;V).$$
with $\left\{1,2\right\},\left\{2,3\right\},\left\{3,4\right\}$
being the run of length 4. Let $E_1=\e(\left\{1,2\right\};V)$ and
$E_2=\e(\left\{2,3\right\},\left\{3,4\right\},e_4,\dots,e_t;V)$.
For similar reasons to the previous parts we have
$\Hred_{i}(E_1)=\Hred_{i}(E_2)=0$ for all $i$ and so $\Hred
_i(E)=\Hred_{i-1}(E_1\cap E_2)$ for all $i$. We see that
\begin{eqnarray*}
E_1 \cap E_2 &=& \e(\left\{1,2,3\right\},\left\{1,2,3,4\right\},
\left\{1,2\right\}\cup e_4,\dots,\left\{1,2\right\} \cup e_t;V)\\
&=& \e(\left\{3\right\},\left\{3,4\right\}, e_4,\dots,\ e_t;V')\\
&=& \e(\left\{3\right\},e_4,\dots,\ e_t;V')\\
&\ &(\rm{using \ Lemma\ } (\ref{lemC}))
\end{eqnarray*}
where $V'=V \setminus \left\{1,2\right\}$. Since $4 \in V'
\setminus ( \{3 \} \cup \bigcup_{i=4}^t e_i)$ we can use Lemma
(\ref{lemA}) to see that $\Hred_i(\e(\left\{3\right\},e_4,\dots,\
e_t;V'))=0$ for all $i$. Therefore we also have $\Hred_i(E)=0$ for
all $i$.
\end{proof}

We now find the reduced homology of the simplicial complexes
$E(2)$ and $E(3)$. All the simplicial complexes we are interested
in here (those of the form $E(s_1,\dots,s_r)$) will have reduced
homology which is some shift of the reduced homology of $E(2)$ or
$E(3)$.

\begin{thm}{\Prop}
For $E=E(2)$ we have
$$\Hred_i(E(2))=
\left\{
\begin{array}{ll}
k \ & \ if \  i=-1\\
0 \ &  \ if \ i\neq -1.
\end{array}
\right.$$

\end{thm}

\begin{proof} The result follows from observing that $E$ is in
fact the empty simplicial complex. We may write
\begin{eqnarray*}
E &=& \e(\{1,2\};\{1,2\})\ \  (\rm{this\  is\  the\  simplicial\
complex} \\
& \ & \rm{ \ with \ one \ maximal \ face \ } \{1,2\} \setminus
\{1,2\}) \\
  &=& \{\emptyset \}.
\end{eqnarray*}
\end{proof}

\begin{thm}{\Prop}\label{E(3)}
For $E=E(3)$ we have
$$\Hred_i(E)=
\left\{
\begin{array}{ll}
k \ & \ if \  i=0\\
0 \ & \ if \ i\neq 0.
\end{array}
\right.$$

\end{thm}

\begin{proof}
We may write $E$ in the following way
\begin{eqnarray*}
E &=&\e( \{1,2 \},\{2,3\} ;\{ 1,2,3 \} ) \\
  &=&\lc \lc 3 \rc, \lc 1 \rc \rc.
\end{eqnarray*}
This is the simplicial complex which is made up of two disjoint 0
dimensional faces hence the result.
\end{proof}

\begin{thm}{\Prop}\label{1mod3runs}
If $E=E(s_1,\dots,s_r,3m+1)$, for some $m\geqslant1$, then
$\Hred_{i}(E)=0$ for all $i$.
\end{thm}

\begin{proof}We may describe $E$ as follows. Let
$$E=\e(\left\{v_1,v_2\right\},\left\{v_2,v_3\right\},\dots,\left\{v_{3m},v_{3m+1}\right\},e_1,\dots,e_t;V).$$
where $v_1, \dots , v_{3m+1}$ is the run of length $3m+1$. We will
split $E$ into two complexes and use the Mayer-Vietoris sequence.
Let $E_1=\e(\{v_{3m}, v_{3m+1} \};V)$ and let
$$E_2=\e(\left\{v_1,v_2\right\},\left\{v_2,v_3\right\},\dots,\left\{v_{3m-1},v_{3m}\right\},e_1,\dots,e_t;V).$$
We have $\Hred_i(E_1)=0$ for all $i$ because $E_1$ is a simplex
and $\Hred_i(E_2)=0$ for all $i$ by Lemma (\ref{lemA}) because
$v_{3m+1}$ is in every maximal face of $E_2$. We also have $E=E_1
\cup E_2$ and so we may use the Mayer-Vietoris sequence
$$\dots \to \Hred_i(E_1) \oplus \Hred_i(E_2) \to \Hred_i(E)\to
\Hred_{i-1}(E_1 \cap E_2)\to \dots$$ which becomes
$$\dots \to 0 \to \Hred_i(E)\to \Hred_{i-1}(E_1 \cap E_2)\to 0 \dots$$
from which we see that $\Hred_i(E) \cong \Hred_{i-1}( E_1 \cap
E_2)$ for all $i$. So we now examine the complex $E_1 \cap E_2$.

\begin{eqnarray*}
E_1 \cap E_2 &=&
\e(\{v_1,v_2,v_{3m},v_{3m+1}\},\{v_2,v_3,v_{3m},v_{3m+1}\},\dots\\
&\ &\dots,\{v_{3m-2},v_{3m-1},v_{3m},v_{3m+1}\},
\{v_{3m-1},v_{3m},v_{3m+1}\},e_1,\dots,e_t;V)\\
             &=&
\e(\{v_1,v_2\},\{v_2,v_3\},\dots,\{v_{3m-3},v_{3m-2}\},\\
& &\{v_{3m-2},v_{3m-1}\}, \{v_{3m-1} \},e_1,\dots,e_t;V')\\
             &=&
\e(\{v_1,v_2\},\{v_2,v_3\},\dots,\{v_{3m-3},v_{3m-2}\}, \{v_{3m-1} \},e_1,\dots,e_t;V')\\
\end{eqnarray*}
where $V'=V \setminus \{v_{3m},v_{3m+1} \}$. By Lemma
(\ref{lemD}), since $$v_{3m-1} \notin (\bigcup_{j=1}^t e_j )\cup
(\{v_1,v_2\} \cup \dots \cup \{v_{3m-3}, v_{3m-2} \}).$$ we deduce
that, for all $i$,
\begin{eqnarray*}
\Hred_i(E) &\cong& \Hred_{i-1}(E_1 \cap E_2)\\
           &\cong& \Hred_{i-2}(\e(\{v_1,v_2\},\{v_2,v_3\},\dots,\{v_{3m-3},v_{3m-2}\},
           e_1,\dots,e_t;V' \setminus \{ v_{3m-1} \}) \\
           &=&\Hred_{i-2}(\e(\{v_1,v_2\},\{v_2,v_3\},\dots,\{v_{3(m-1)},v_{3(m-1)+1}\},
           e_1,\dots,e_t;V' \setminus \{ v_{3m-1} \}).
\end{eqnarray*}
The simplicial complex
$$\e(\{v_1,v_2\},\{v_2,v_3\},\dots,\{v_{3(m-1)},v_{3(m-1)+1}\},
e_1,\dots,e_t;V' \setminus \{ v_{3m-1} \})$$ may be described as
$E(s_1,\dots,s_r,3(m-1)+1)$. We have shown that the reduced
homology of $E$ is a shift of the homology of $$E(s_1,\dots,s_r,
3(m-1)+1),$$\ie
$\Hred_i(E(s_1,\dots,s_r,3m+1))\cong\Hred_{i-2}(E(s_1,\dots,s_r,3(m-1)+1))$
for all $i$. We can make repeated use of the above method to
obtain
$$\Hred_i(E) \cong \Hred_{i-2m+2}(E(s_1,\dots,s_r,4))$$ for all
$i$. And so, by Proposition (\ref{homshift2}), we have
$\Hred_i(E)=0$ for all $i$.
\end{proof}

We have now deduced that if $F$ is a non-empty face of $\AD(C_n)$
such that $\Link_{\AD(C_n)} F$ has some non-zero reduced homology
modules then$$\Link_{\AD(C_n)}
F=E(3p_1\dots,3p_\alpha,3q_1+2,\dots,3q_\beta+2).$$ That is to say
$\Link_{\AD(C_n)} F$ corresponds to a subgraph of $C_n$ with runs
only of lengths 0 and 2 mod  3. If we use the the formula
$$\b_i(C_n)=\sum_{H \subset C_n} {\dim}_k \Hred_{i-2}(\e(H))$$
(where the sum is over induced subgraphs of $C_n$ which contain at least one edge)
we
may exclude all the induced subgraphs of $C_n$ with any runs of
length 1 mod 3 from consideration as these will contribute $0$.

We will now describe the reduced homology modules of
$\Link_{\AD(C_n)} F$ as follows.

\begin{thm}{\Lem}\label{E(2,3)}
Write $E(3^{\a},2^{\b})$ to denote the simplicial complex
associated with the graph which is $\a$ runs of length $2$ and
$\b$ runs of length $3$.
$$\Hred_i(E(3^{\a},2^{\b}))=
\left\{
\begin{array}{ll}
k \ & \ if \  i=2\alpha +\beta-2\\
0 \ &  \ if \ i\neq 2\alpha +\beta-2.
\end{array}
\right.$$
\end{thm}

\begin{proof}
We first make repeated use of (\ref{homshift2}) part $(i)$ to see
that $\Hred_i(E(3^{\a},2^{\b})) \cong \Hred_{i-\b}(E(3^{\a}))$ for
all $i$. Now we use (\ref{homshift2}) part $(ii)$ to obtain
\begin{eqnarray*}
\Hred_i(E(3,\dots,3,2\dots,2)) &\cong& \Hred_{i-\b-2(\alpha
-1)}(E(3))\\
                               &=&\left\{
                               \begin{array}{ll}
                               k \ & \ if \  i=2\alpha +\beta-2\\
                               0 \ &  \ if \ i\neq 2\alpha +\beta-2
                               \end{array}
                               \right.
\end{eqnarray*}
using Lemma (\ref{E(3)}), which gives the reduced homology of
$E(3)$.
\end{proof}

\begin{thm}{\Prop}\label{homofE}
Let $E=\nzlink$, for integers $p_1,\dots,p_{\a},q_1,\dots,q_{\b}$,
and let $P=\sum_{j=1}^{\alpha}p_j$ and $Q=\sum_{j=1}^{\beta}q_j.$
Then
$$\dim_k(\Hred_i(E))=
\left\{
\begin{array}{ll}
1 \ if \  i=2(P+Q)+\beta -2\\
0 \ if \  i \neq2(P+Q)+\beta -2.
\end{array}
\right.$$
\end{thm}

\begin{proof}
Repeatedly applying (\ref{homshift}) we see that
$$
\begin{array}{lll}
\Hred_i(E)
&=&\Hred_{i-2}(E(3p_1-3,\dots,3p_{\alpha},3q_1+2,\dots,3q_{\beta}+2))\\
&=&\Hred_{i-2(p_1-1)}(E(3,3p_2,\dots,3p_{\alpha},3q_1+2,\dots,3q_{\beta}+2))\\
& &({\rm{Applying}} \ (\ref{homshift})\ p_1-1 \  {\rm{times \ to \ the \ first \ run}})\\
&=&\Hred_{i-2(p_1-1)-2(p_2-1)}(E(3,3,3p_3,\dots,3p_{\alpha},3q_1+2,\dots,3q_{\beta}+2))\\
& &({\rm{Applying}} \ (\ref{homshift})\ p_2-1 \  {\rm{times \ to \ the \ second \ run}})\\
&=&\Hred_{i-2P+2\alpha}(E(3^{\a},3q_1+2,\dots,3q_{\beta}+2))\\
&=&\Hred_{i-2P+2\alpha-2}(E(3^{\a},3(q_1-1)+2,\dots,3q_{\beta}+2))\\
& &({\rm{Applying}} \ (\ref{homshift})\ {\rm{\ to \ the \  next \ run}})\\
&=&\Hred_{i-2P+2\alpha-2q_1}(E(3^{\a},2,\dots,3q_{\beta}+2))\\
&=&\Hred_{i-2(P+Q)+2\alpha}(E(3^{\a},2^{\b})).\\
\end{array}$$
We now use this together with Lemma (\ref{E(2,3)}), which shows us
that $\Hred_i(E(3^{\a},2^{\b}))=k$ if and only if $i=2\alpha +\b
-2$ (and is zero otherwise), to conclude that
\begin{eqnarray*}
\Hred_i(E))&=& \left\{
\begin{array}{ll}
k \ \ if \ \  i-2(P+Q)+2\alpha=2\alpha + \b -2\\
0 \ \ if \ \ i-2(P+Q)+2\alpha \neq 2\alpha + \b -2\\
\end{array}
\right.\\
&=& \left\{
\begin{array}{ll}
k \ \ if \ \  i=2(P+Q)+ \b -2\\
0 \ \ if \ \  i \neq2(P+Q)+ \b -2.\\
\end{array}
\right.\\
\end{eqnarray*}
and the result is proved.
\end{proof}

\section{Betti Numbers from Counting Induced Subgraphs}

To find the Betti numbers of $C_n$ we use the formula
$$\beta_i(C_n)=\sum_{F \in \AD(C_n)} \dim_k \Hred_{i-2}(\Link_{\AD(C_n)}
 F;k)$$ from (\ref{ER}). Along with the above result this shows
that $\beta_i(C_n)$ is equal to the number of faces of $\AD(C_n)$
such that $\dim_k \Hred_{i-2}(\Link_{\AD(C_n)} F)=1$. This is the
number of faces of $\AD(C_n)$ with link of the form $$\nzlink$$
such that $i=2(P+Q) + \beta .,$ where $P=\sum_{j=1}^{\alpha}p_j$
and $Q=\sum_{j=1}^{\beta}q_j.$

It is now convenient to consider the $\mathbb{N}$-graded Betti
numbers of $C_n$ separately. This is because when using the
formula of (\ref{ER2}) the sum which gives $\b_{i,d}$ is taken
over all induced subgraphs which have precisely $d$ vertices.

\begin{thm}[\rm]{Remark} We will use Theorem (\ref{ER2}).
The $i$th  $\mathbb{N}$-graded Betti number of degree $d$ of $I$
is
\begin{eqnarray*}
\beta_{i,d}(C_n)&=&\sum_{F \in \AD(C_n): \  |V(C_n) \setminus
F|=d} {\dim}_k
\Hred_{i-2}(\Link_{\AD(C_n)} F;k)\\
           &=&\sum_{H \subset C_n, |V(H)|=d} {\dim}_k
           \Hred_{i-2}(\e(H))
\end{eqnarray*}
Note that $$\beta_i(C_n)=\sum_{d \in \mathbb{N}}\beta_{i,d}(C_n)$$
\end{thm}

\begin{thm}{\Lem}\label{counting}
The $i${\rm th} Betti number of degree $d$, $\beta_{i,d}$, may be
obtained by counting the number of faces of $\AD(C_n)$ which have
link of the form $$\nzlink,$$for integers
$p_1,\dots,p_{\a},q_1,\dots,q_{\b}$, such that
$$
\begin{array}{lll}
&(i)& i=2(P+Q)+\beta\\
&(ii)& d=3(P+Q)+2\beta\\
\end{array}$$
(where $P=\sum_{j=1}^{\a}p_j$ and $Q=\sum_{j=1}^{\b}q_j$).
Equivalently, $\b_{i,d}$ is the number of induced subgraphs of
$C_n$ which are composed of runs of lengths
$$3p_1,\dots,3p_{\alpha},3q_1+2,\dots,3q_{\b}+2$$ satisfying the above two conditions.
\end{thm}

\begin{proof}
As previously noted the only links of faces of $\AD(C_n)$ which
may have non-zero reduced homology are of the form $\nzlink$ which
corresponds to a subgraph of $C_n$ which has runs of lengths
$$3p_1,\dots,3p_{\alpha},3q_1+2,\dots,3q_{\b}+2.$$ Such a subgraph
has $d$ vertices, \ie
$d=3p_1+\dots+3p_{\alpha}+3q_1+2+\dots+3q_{\b}+2=3(P+Q)+2\b$,
which provides the condition $(ii)$. By Proposition (\ref{homofE})
the reduced homology module, $\Hred_{i-2}(\nzlink)$, is non-zero
(and 1-dimensional) if and only if $i=2(P+Q)+\b$, which is
condition $(i).$ Therefore each such induced subgraph  of $C_n$
contributes `1' to $\b_{i,d}$ in the formula $\b_{i,d}=\sum_{H
\subset C_n,|V(H)|=d}{\dim}_k\Hred_{i-2}(\e(H))$.
\end{proof}

\begin{thm}{\Prop}
If $d>2i$ then $\beta_{i,d}(C_n)=0.$
\end{thm}

\begin{proof}
This is a consequence of Theorem (\ref{greater2i}) as $I(C_n)$ is
generated by monomials of degree 2.
\end{proof}

\begin{thm}{\Prop}\label{betai2i}
The $i${\rm th} Betti number of degree $d=2i$, $\b_{i,2i}(C_n)$,
is the number of induced subgraphs of $C_n$ which are comprised of
$i$ runs of length $2$.

\end{thm}

\begin{proof}
We consider when a face of $\AD(C_n)$ will contribute `1' to
$\b_{i,2i}$ in the sum from (\ref{ER2}). From Lemma
(\ref{counting}) the face $F$ of $\AD(C_n)$ will contribute to
$\b_{i,2i}(C_n)$ if and only if $\Link_{\AD(C_n)} F=\nzlink$, for
some integers $p_1,\dots,p_{\alpha},q_1,\dots,q_{\b}$ such that
$i=2(P+Q)+\b$ and $2i=3(P+Q) + 2\b$ where
$P=\sum_{j=1}^{\alpha}p_j$, $Q=\sum_{j=1}^{\b}q_j$. These
conditions imply that $2(P+Q)+\b =i=2i-i=P+Q+\b$. As all these
numbers are non-negative integers, this is only possible if
$P=Q=0$ which in turn implies that
$p_1=\dots=p_{\alpha}=q_1=\dots=q_\b=0$ and $i=\b$. Hence
$\Link_{\AD(C_n)} F=E(2^i)$. The $i$th Betti number of degree $2i$
is thus obtained by counting all the induced subgraphs of $C_n$
which correspond to the links of such faces. These are the
subgraphs which are $i$ runs of length 2.
\end{proof}

\section{Counting Runs}

As shown above calculating Betti numbers of cycles is equivalent
to counting particular types of induced subgraphs. Here we
demonstrate how this may be done.

\begin{thm}[\rm]{\Def}\label{defC}
(i) Define $B(n,l,m)$ to be the number of ways of choosing an
induced subgraph of
$C_n$ which consists of $m$ runs of length $l$. \\
(ii) Define $C(n,l,m)$ to be the number of such subgraphs in which
neither 1 nor $n$ feature in a run.
\end{thm}

\begin{thm}{\Lem}\label{Cnlm}
The number of choices of $m$ runs of length $l$ in $C_n$ such that
neither 1 nor n are featured in any run is
$$C(n,l,m)={{n-lm-1}\choose{m}}.$$
\end{thm}

\begin{proof}
We can find $C(n,l,m)$ by considering the number of ways of
separating the vertices which are not in any run by the $m$ runs.
There are $n-lm$ vertices (including the vertices $1$ and $n$)
which cannot be in any runs. We can construct all the appropriate
subgraphs by placing the $m$ runs of length $l$ in the spaces
between the adjacent vertices which are not in any run (except
between 1 and $n$). There are $n-lm-1$ spaces into each we can
choose to place, or not place, a run. There are
${{n-lm-1}\choose{m}}$ ways of doing this.
\end{proof}

\begin{thm}{\Lem}\label{Bnlm}
The number of choices of $m$ runs of length $l$ in $C_n$ is
$$B(n,l,m)={{n-lm}\choose{m}}+l{{n-lm-1}\choose{m-1}}.$$
\end{thm}

\begin{proof}
We first count those subgraphs which have at least one of 1 and
$n$ in a run. There are $l+1$ possible runs of length $l$ which
contain $1$ or $n$ (or both). The other $m-1$ runs are selected
from the remaining $n-l$ vertices. Since the vertices adjacent to
the run containing 1 or $n$ cannot be included we have
$C(n-l,l,m-1)$ choices. Hence
\begin{eqnarray*}
B(n,l,m)&=&C(n,l,m)+(l+1)C(n-l,l,m-1) \\
\\
&=&{{n-lm-1}\choose{m}} +(l+1){{n-lm-1}\choose{m-1}} \\
\\
&=&{{n-lm}\choose{m}} +l{{n-lm-1}\choose{m-1}}.
\end{eqnarray*}
\end{proof}

\begin{thm}{\Cor}
The i{\rm th} Betti number of degree $2i$ of $C_n$ is

$$\b_{i,2i}(C_n)={{n-2i}\choose i}+2{{n-2i-1}\choose {i-1}}$$
\end{thm}

\begin{proof}
By (\ref{betai2i}) $\beta_{i,2i}(C_n)$ is the number of subgraphs
of $C_n$ which are $i$ runs of length 2. That is,
$\b_{i,2i}(C_n)=B(n,2,i)$. Using (\ref{Bnlm})
$$\b_{i,2i}(C_n)=B(n,2,i)={{n-2i}\choose i}+2{{n-2i-1}\choose {i-1}}.$$
\end{proof}

\begin{thm}{\Lem}\label{constructGraph}
Let $A$ be an induced subgraph of $C_n$ composed of $i$ runs of
length $2$ and let $j \leqslant i$. Select $j$ of the vertices of
$C_n$ which are adjacent to, and clockwise from, the runs of $A$.
Let $A'$ be the induced subgraph of $C_n$ which has the $j$
vertices selected together with the vertices of $A$. The graph
$A'$ is comprised of runs of lengths $0 \mod 3$ and $2 \mod 3$
only.
\end{thm}

\begin{proof}
A run of $A'$ will be of one of the following forms.

(i) It was constructed from a string of $r$ runs of length 2 of
$A$ which were separated by single vertices. The new run includes
these formerly separating vertices along with all the vertices of
the $r$ runs of length 2. It will have length $2r+(r-1)=3r-1
\equiv 2\mod 3$.

(ii) It was constructed as in part (i) but also includes the
vertex immediately clockwise of the clockwise most of the $r$
runs. In this case the new run will be of length $2r+r=3r \equiv 0
\mod 3.$

Note the case $r=1$ in the above gives a run of length 2 or 3.
\end{proof}

\begin{thm}{\Prop}\label{bet1}
For $2i+j \neq n$
\begin{eqnarray*}
\b_{i+j,2i+j}(C_n) &=& {i\choose j} \b_{i,2i} \\
              &=& {i\choose j} \lc {{n-2i}\choose i}+2{{n-2i-1}\choose
              {i-1}} \rc
\end{eqnarray*}
\end{thm}

\begin{proof}
First observe that any $\mathbb{N}$-graded Betti number can be
written in the above way since $\b_{l,d}=\b_{i+j,2i+j}$ where
$i=d-l$ and $j=2l-d$. By Lemma (\ref{counting}) $\b_{i+j,2i+j}$ is
the number of faces of $\AD(C_n)$ which have link of the form
$E(3p_1,\dots,3p_{\alpha},3q_1+2,\dots,3q_{\b}+2)$, for some
integers $p_1,\dots,p_{\a},q_1,\dots,q_{\b}$, such that
$$
\left\{
\begin{array}{lll}
i+j &=& 2(P+Q) +\b \\
2i+j &=& 3(P+Q) +2\b
\end{array}
\right.
$$
where $P=\sum_{j=1}^{\alpha}p_j$ and $Q=\sum_{j=1}^{\b}q_j$. This
is equivalent to counting the induced subgraphs of $C_n$ which
consist of runs of lengths
$3p_1,\dots,3p_{\alpha},3q_1+2,\dots,3q_{\b}+2$, for some
$p_1,\dots,p_{\a},q_1,\dots,q_{\b}$, satisfying the above
conditions.

We can construct any such subgraph starting from a subgraph of $i$
runs of length 2 in $C_n$. Let $A\subset C_n$ denote a subgraph
which consists of $i$ runs of length 2.

There are $i$ vertices of $C_n$ which are adjacent to, and
clockwise from, the runs of $A$. Select any $j$ of these vertices
and add them to the run (or runs) to which they are adjacent. By
Lemma (\ref{constructGraph}) the new subgraph of $C_n$ which this
defines, $A'$, say, will be comprised of runs of lengths $0 \mod
3$ and $2 \mod 3$ only.

Say the runs are of lengths
$3p_1,\dots,3p_{\alpha},3q_1+2,\dots,3q_{\b}+2.$ The number of
vertices in $A'$ is clearly $2i+j=3P+3Q +2\b$ where
$P=\sum_{i=1}^{\alpha}p_i$ and $Q=\sum_{i=1}^{\b}q_i$. To show
$A'$ is a subgraph of the type contributing to $\b_{i+j,2i+j}$ it
remains to demonstrate that $i+j=2(P+Q) + \b$.

When $j=0$ we see that $P+Q=0$ as we have only runs of length 2.
In this case $i+j=i=\b=2(P+Q) +\b$. Now suppose that $j>0$. Each
selection of one of the $j$ vertices either forms a run of length
3 or else increases a run by 3. Whichever of these happens $P+Q$
increases by one. Adding $j$ vertices will increase $P+Q$ by $j.$
Because $P+Q=0$ when $j=0$ this implies that $j=P+Q$. Hence
\begin{eqnarray*}
2i&=&3P+3Q+2\b-j\\
  &=&3P+3Q +2\b -P-Q \\
  &=&2P+2Q +2\b. \\
\end{eqnarray*}
And so $i=P+Q+\b$. Therefore we have $i+j=2(P+Q)+\b$ as required.

We now show that every subgraph of $C_n$ which contributes to
$\b_{i+j,2i+j}$ can be constructed in this way. Suppose that
$$E=E(3p_1,\dots,3p_{\alpha},3q_1+2,\dots,3q_{\b}+2)$$ is the link
of a face of $\AD(C_n)$ satisfying the above conditions. Consider
the associated subgraph of $C_n$, which is comprised of runs of
lengths $$3p_1,\dots,3p_{\alpha},3q_1+2,\dots,3q_{\b}+2.$$ For a
run of length $l \equiv 0 \mod 3$ with vertices
$a_1,a_2,\dots,a_l$ remove the vertices $a_3,a_6,\dots,a_{l}$.
Repeat this with every run of length $0 \mod 3$. For a run of
length $m \equiv 2\mod 3$ with vertices $b_1,b_2,\dots,b_m$ remove
the vertices $b_3,b_6,b_9,\dots,b_{m-2}$. This leaves a graph of
$i$ runs of length 2 from which our earlier construction gives the
graph associated with $E$.

We count the graphs we can construct in this way. For any graph
which is $i$ runs of length 2 we are choosing $j$ of the $i$
vertices which are adjacent and clockwise to the runs. So we have
${i \choose j}$ graphs constructed for each arrangement of $i$
runs of length 2, of which there are $\b_{i,2i}(C_n)$, by
Proposition (\ref{betai2i}). Hence
\begin{eqnarray*}
\b_{i+j,2i+j}(C_n) &=& {i\choose j} \b_{i,2i}(C_n) \\
              &=& {i\choose j}\lc {{n-2i}\choose i}+2{{n-2i-1}\choose {i-1}}\rc.
\end{eqnarray*}
\end{proof}

\section{Betti Numbers of Degree $n$}

It only remains to calculate the reduced homology of
$\Link_{\AD(C_n)} \emptyset$. It can be easily seen from
definition (\ref{link}) that $\Link_{\AD(C_n)} \emptyset
=\AD(C_n)$. Note that, using the correspondence between faces of
$\AD(C_n)$ and induced subgraphs of $C_n$, the face $\emptyset$ is
associated with $C_n$. This is the only subgraph of $C_n$ which
cannot be thought of as a collection of runs and so must be
considered separately.

\begin{thm}[\rm]{Notation} For any positive integer $m$ we will
write $$[m]=\{1,2,\dots,m \}.$$
\end{thm}

\begin{thm}{\Lem}\label{adhom}
Define the simplicial complex $\Gamma$ as follows
$$\Gamma = \e(\left\{1,2\right\},\left\{2,3\right\},\dots,\left\{m-1,m
\right\},\left\{m+1\right\},\left\{m+2\right\};[m+2])$$ Then
$\Gamma$ has reduced homology which is a shift of the reduced
homology of $E(m)$, the simplicial complex associated with the
graph which is a single run of length $m$, as follows
\begin{eqnarray*}
\Hred_i(\Gamma) &=&
\Hred_{i-2}(\e(\left\{1,2\right\},\left\{2,3\right\},\dots,\left\{m-1,m
\right\};[m]) \\
&=&\Hred_{i-2}(E(m)).
\end{eqnarray*}
\end{thm}

\begin{proof}
This is just Lemma (\ref{lemD}) used twice. Note that the vertex
$m+2$ is in none of the sets
$\{1,2\},\{2,3\},\dots,\{m-1,m\},\{m+1\}$ so Lemma (\ref{lemD})
implies
\begin{eqnarray*}
\Hred_i(\Gamma)&=&\Hred_i(\e(\left\{1,2\right\},\left\{2,3\right\},\dots,\left\{m-1,m
\right\},\left\{m+1\right\},\left\{m+2\right\};[m+2]) \\
               &=&\Hred_{i-1}(\e(\left\{1,2\right\},\left\{2,3\right\},\dots,\left\{m-1,m
\right\},\left\{m+1\right\};[m+1])
\end{eqnarray*}
for all i. Now we note that the vertex $m+1$ is not in any of the
sets $\{1,2\},\{2,3\},\dots,\{m-1,m\}$ and we apply Lemma
(\ref{lemD}) once more to obtain
\begin{eqnarray*}
\Hred_i(\Gamma)&=&\Hred_{i-1}(\e(\left\{1,2\right\},\left\{2,3\right\},\dots,\left\{m-1,m
\right\},\left\{m+1\right\};[m+1])\\
               &=&\Hred_{i-2}(\e(\left\{1,2\right\},\left\{2,3\right\},\dots,\left\{m-1,m
\right\};[m])\\
               &=&\Hred_{i-2}(E(m)).
\end{eqnarray*}
for all i.

\end{proof}

\begin{thm}{\Lem}\label{lemin}
The i{\rm th} Betti number of degree n is
\begin{eqnarray*}
\b_{i,n}(C_n)&=&{\dim}_k \Hred_{i-2}(\AD(C_n);k)\\
             &=&{\dim}_k \Hred_{i-2}(\e(C_n);k).
\end{eqnarray*}
\end{thm}

\begin{proof}
Using the formula from (\ref{ER}) we have
$$\b_{i,n}(C_n) =\sum_{F
\in \AD(C_n) : |V(C_n) \setminus F|=n} {\dim}_k
\Hred_{i-2}(\Link_{\AD(C_n)}F;k)=\sum_{H \subset C_n} {\dim}_k
\Hred_{i-2}(\e(H);k).$$ The only face of $\AD(C_n)$ which is such
that its complement is of size $n$ is the empty set $\emptyset$.
Also $\Link_{\AD(C_n)} \emptyset = \AD(C_n)$ so we have
$\b_{i,n}(C_n)={\dim}_k \Hred_{i-2}(\AD(C_n))$. We can also write
$\AD(C_n)$ as $\e(C_n)$.
\end{proof}

\begin{thm}{\Prop}\label{prop2}
The non-zero Betti numbers of $C_n$ of degree $n$ are as follows:
$$
\begin{array}{lll}
(i) &If& n\equiv 1 \mod 3, \  \b_{\frac{2n+1}{3},n}(C_n)=1 \\
(ii) &If& n\equiv 2 \mod 3, \  \b_{\frac{2n-1}{3},n}(C_n)=1 \\
(iii) &If& n \equiv 0 \mod 3, \ \ \ \ \
\b_{\frac{2n}{3},n}(C_n)=2.
\end{array}
$$
\end{thm}

\begin{proof}
As commented above $\Link_{\AD(C_n)} \emptyset =\AD(C_n)$. We can
write the Alexander dual as
\begin{eqnarray*}
\AD(C_n) &=& \e(\lc 1,2 \rc,\lc 2,3 \rc ,\dots, \lc n-1,n \rc, \lc n,1 \rc ;[n]) \\
&=& E_1 \cup E_2
\end{eqnarray*}
where $E_1=\e(\lc 1,2 \rc,\lc 2,3 \rc ,\dots, \lc n-1,n \rc ;[n])$
and $E_2=\e( \lc n,1 \rc ;[n]).$ The intersection of these
simplicial complexes is
\begin{eqnarray*}
E_1 \cap E_2&=&\e(\lc 1,2,n \rc,\lc 1,2,3,n \rc
,\{3,4,1,n\},\dots \\
& \ & \dots, \lc n-3,n-2,1,n \rc,\{n-2,n-1,n,1\}
\lc n-1,n,1 \rc ;[n] )\\
&=&\e(\lc 2 \rc,\lc 3,4 \rc ,\dots, \lc n-3,n-2 \rc,
\lc n-1 \rc ;[n] \setminus \lc 1,n \rc)\\
& &{\rm (Using \ Lemmas \ (\ref{lemB})\  and\  (\ref{lemC})).}
\end{eqnarray*}
 By Lemma (\ref{adhom}) $\Hred_i(E_1 \cap
E_2)=\Hred_{i-2}(E(n-4)).$

$(i)$ We now consider the first of these cases. Here we have
$n=3m+1$ for some $m$. By Proposition (\ref{1mod3runs})
$\Hred_{i}(E_1)=0$ for all $i$. We also have $\Hred_{i}(E_2)=0$
for all $i$ since $E_2$ is a simplex. Hence we can use the
Mayer-Vietoris sequence
$$\dots\to\Hred_i(E_1)\oplus\Hred_i(E_2)\to\Hred_i(\AD)\to\Hred_{i-1}(E_1\cap
E_2)\to\Hred_{i-1}(E_1)\oplus\Hred_{i-1}(E_2)\to\dots$$ to see
that
\begin{eqnarray*}
\Hred_i(\AD(C_n))&=&\Hred_{i-1}(E_1 \cap E_2)) \\
&=&\Hred_{i-3}(E(n-4))\\
&=&\Hred_{i-3}(E(3(m-1))).
\end{eqnarray*}
By Proposition (\ref{homofE})
\begin{eqnarray*}
\Hred_i(\AD(C_n)) &=& \Hred_{i-3}(E(3(m-1))) \\
             &=& \lc \begin{array}{lll}
                     k &if& i-3=2(m-1)-2 \\
                     0 & & otherwise
                     \end{array}
                     \right. \\
             &=& \lc \begin{array}{lll}
                     k &if& i=2m-1 \\
                     0 & & otherwise.
                     \end{array}
                     \right.
\end{eqnarray*}

We used this along with Lemma (\ref{lemin}) to find
$\b_{i,n}(C_n)$:
\begin{eqnarray*}
\b_{i,n}(C_n) &=& \dim_k \Hred_{i-2}(\Link_{\AD(C_n)} \emptyset) \\
         &=& \dim_k \Hred_{i-2}(\AD(C_n)) \\
         &=& \lc \begin{array}{lll}
                 1 &if& i=2m+1 \\
                 0 &if& i \neq 2m+1.
                 \end{array}
                 \right.
\end{eqnarray*}

$(ii)$ In this case $n=3m+2$ for some $m$. Lemma (\ref{adhom}) and
Proposition (\ref{1mod3runs}) imply
\begin{eqnarray*}
\Hred_i(E_1 \cap E_2) &=& \Hred_{i-2}(E(n-4)) \\
                      &=& \Hred_{i-2}(E(3(m-1) +1)) \\
                      &=& 0
\end{eqnarray*}
for all $i$. The Mayer-Vietoris sequence
$$\dots\to\Hred_i(E_1)\oplus\Hred_i(E_2)\to\Hred_i(\AD)\to\Hred_{i-1}(E_1\cap
E_2)\to\Hred_{i-1}(E_1)\oplus\Hred_{i-1}(E_2)\to\dots$$
becomes
$$\dots\to 0
\to\Hred_i(E_1)\oplus\Hred_i(E_2)\to\Hred_i(\AD(C_n))\to 0 \to
\dots$$ and so we obtain
\begin{eqnarray*}
\Hred_i(E_1) \oplus \Hred_i(E_2) &=&\Hred_i(E_1) \\
& &(\rm{as \ } E_2 \rm{\ is \ a \ simplex\ })\\
 &=& \Hred_i(\AD(C_n)).
\end{eqnarray*}
By Proposition (\ref{homofE})
\begin{eqnarray*}
\Hred_i(\AD(C_n)) &=& \Hred_i(E_1) \\
             &=& \lc \begin{array}{lll}
                     k &if& i=2m-1 \\
                     0 &if& i \neq 2m-1.
                     \end{array}
                     \right.
\end{eqnarray*}
As in the previous case
\begin{eqnarray*}
\b_{i,n}&=& \dim_k\Hred_{i-2}(\AD(C_n)) \\
        &=& \lc \begin{array}{lll}
                1 &if& i=2m+1 \\
                0 &if& i \neq 2m+1.
                \end{array}
                \right.
\end{eqnarray*}
$(iii)$ Finally we have the case where $n=3m$ for some $m$. We now
use the Mayer-Vietoris sequence
$$\dots \to \Hred_i( E_1 \cap E_2) \to \Hred_i(E_1) \oplus
\Hred_i(E_2) \to \Hred_i(\AD)\to \dots$$ As before
$\Hred_i(E_2)=0$ for all $i$ and $\Hred_i(E_1 \cap E_2) =
\Hred_{i-2}(E(n-4))=\Hred_{i-2}(E(3m-4)$ for all $i$. So the
sequence becomes
$$\dots \to \Hred_{i-2}(E(3m-4)) \to \Hred_i(E(3m)) \to
\Hred_i(\AD)\to \dots$$ From Theorem (\ref{homshift}) we obtain
\begin{eqnarray*}
\Hred_i(E(3m))&=&\Hred_{i-2(m-1)}(E(3)) \\
             &=& \lc \begin{array}{lll}
                     k &if& i-2(m-1)=0 \\
                     0 &if& i-2(m-1) \neq 0
                     \end{array}
                     \right.\\
             &=& \lc \begin{array}{lll}
                     k &if& i=2m-2 \\
                     0 &if& i\neq 2m-2
                     \end{array}
                     \right.
\end{eqnarray*}
and also using (\ref{homshift})
\begin{eqnarray*}
\Hred_i(E(3m-4))    &=& \Hred_i(E(3(m-2)+2))\\
             &=&\Hred_{i-2(m-2)}(E(2)) \\
             &=& \lc \begin{array}{lll}
                     k &if& i-2(m-2)=-1 \\
                     0 &if& i-2(m-2) \neq -1
                     \end{array}
                     \right.\\
             &=& \lc \begin{array}{lll}
                     k &if& i=2m-5 \\
                     0 &if& i\neq 2m-5.
                     \end{array}
                     \right.
\end{eqnarray*}
If we now look at the appropriate part of the Mayer-Vietoris
sequence we find
 $$\dots\to \Hred_{2m-4}(E(3m-4)) \to
\Hred_{2m-2}(E(3m)) \to \Hred_{2m-2}(\AD(C_n)) \to $$
$$\ \ \ \ \to\Hred_{2m-5}(E(3m-4)) \to
\Hred_{2m-3}(E(3m))\to\dots$$ From the above comments about the
reduced homology of $E(3m)$ and $E(3m-4)$ we convert this into a
short exact sequence of $k$ vector spaces
$$0\to k \to \Hred_{2m-2}(\AD(C_n)) \to k \to 0$$
So we must have $\Hred_{2m-2}(\AD(C_n))=k^2$. We also have that
$\Hred_i(\AD)=0$ for all $i \neq 2m-2$ because the Mayer-Vietoris
sequence is
$$\dots \to 0 \to \Hred_i(\AD(C_n)) \to 0 \to \dots$$
wherever $i \neq 2m-2$. Hence the $i$th Betti number of degree $n$
is
\begin{eqnarray*}
\b_{i,n}(C_n)&=& {\dim}_k \Hred_{i-2}(\Link_{\AD(C_n)} \emptyset;k) \\
        &=& {\dim}_k \Hred_{i-2}(\AD(C_n))\\
        &=& \lc \begin{array}{lll}
                     2 &if& i=2m \\
                     0 &if& i\neq 2m.
                     \end{array}
                     \right.
\end{eqnarray*}
\end{proof}

\section{Betti numbers of Cycles}

We now combine these results to obtain

\begin{thm}{Theorem}\label{betticycle}
The non-zero $\mathbb{N}$-graded Betti numbers of $C_n$ are all in
degree less than or equal to $n$ and for $2i+j < n$
\begin{eqnarray*}
\b_{i+j,2i+j}(C_n) &=& {i\choose j} \b_{i,2i} \\
              &=& {i\choose j} \lc {{n-2i}\choose i}+2{{n-2i-1}\choose
              {i-1}} \rc\\
              &=& \frac{n}{n-2i} {i \choose j}{{n-2i} \choose i}
\end{eqnarray*}
Equivalently, for $d<n$ and $2l \geqslant d$,

$$\b_{l,d}(C_n) = \frac{n}{n-2(d-l)} {{d-l} \choose {2l-d}}{{n-2(d-l)} \choose
{d-l}}$$

and

\begin{eqnarray*}
\b_{2m+1,n}(C_n)=1 &if& n=3m+1 \\
\b_{2m+1,n}(C_n)=1 &if& n=3m+2\\
\b_{2m,n}(C_n)=2 &if& n=3m.
\end{eqnarray*}

\end{thm}
\begin{proof}
This is a combination of Propositions (\ref{bet1}) and
(\ref{prop2}) together with $i=d-l$ and $j=2l-d$.
\end{proof}

\begin{thm}[\rm]{Remark}
Notice that the Betti numbers of cycles do not depend on our
choice of field.
\end{thm}

\begin{thm}{\Cor}
The projective dimension of the cycle graph is independent of the
characteristic of the chosen field and is

$$\pd (C_n)=\left\{
{\begin{array}{lll}
\frac{2n}{3} & if & n \equiv0 \mod 3 \\
\frac{2n+1}{3} & if & n \equiv 1 \mod 3\\
\frac{2n-1}{3} & if & n \equiv 2 \mod 3.
\end{array}} \right.
 $$
\end{thm}

\begin{proof}
The expression given in (\ref{betticycle}), $$\frac{n}{n-2(d-l)}
{{d-l} \choose {2l-d}}{{n-2(d-l)} \choose {d-l}},$$ is non-zero
only if $d-l \geqslant 2l-d$ ,\ie if $2d \geqslant 3l$.

We first consider the case $n=3m$ for some integer $m$. By Theorem
(\ref{betticycle}) $\b_{2m}(C_n) \neq 0$ and so $\pd(C_n)
\geqslant 2m$. Suppose that $l \geqslant 2m+1$. Assume that
$\b_{l,d}(C_n)\neq 0$. We have $2d \geqslant 3l \geqslant
3(2m+1)=6m+3$. Therefore $d \geqslant 3m+\frac{3}{2}>n$ which is
not possible. Hence $\b_{l,d}(C_n)= 0$ and so
$\pd(C_n)=2m=\frac{2n}{3}.$

Now suppose that $n=3m+1$ for some integer $m$. Theorem
(\ref{betticycle}) shows that $\b_{2m+1}(C_n) \neq 0$ and so
$\pd(C_n) \geqslant 2m+1$. Suppose that $l \geqslant 2m+2$. Assume
that $\b_{l,d}(C_n) \neq 0$. This implies that $2d \geqslant 3l
\geqslant 6m+6$ and so $d \geqslant 3m+3 > n$, a contradiction.
Hence $\b_{l,d}(C_n)=0$ and $\pd(C_n)=2m+1=\frac{2n+1}{3}.$

Finally we consider the case $n=3m+2$ for some $m$. Theorem
(\ref{betticycle}) shows that $\b_{2m+1}(C_n) \neq 0$ and so
$\pd(C_n) \geqslant 2m+1$. Now suppose that $l \geqslant 2m+2$ and
that $\b_{l,d}(C_n) \neq 0$. This implies that $2d \geqslant 3l
\geqslant 6m+6$ and so $d \geqslant 3m+3 > n$, a contradiction.
Hence $\b_{l,d}(C_n)=0$ and $\pd(C_n)=2m+1=\frac{2n-1}{3}.$
\end{proof}

\section{Betti Numbers of Lines}

We conclude this chapter by calculating the Betti numbers of a
family of graphs related to the cycles using similar methods.

\begin{thm}[\rm]{\Def}For $n \geqslant 2$, we let $L_n$ denote the line graph on $n$
vertices. This is the graph with vertices $1,2,\dots,n$ and edges
$\{j,j+1\}$ for all $j=1,2,\dots,n-1$. Hence, for some field $k$,
we have $R(L_n)=k[x_1,\dots,x_n]$ and $I(L_n)=\langle
x_1x_2,x_2x_3,\dots,x_{n-1}x_n \rangle .$
\end{thm}

\begin{thm}[\rm]{Remark}
As with cycles, the field $k$ will play no role in calculating the
Betti numbers of lines.
\end{thm}

\begin{figure}[h]
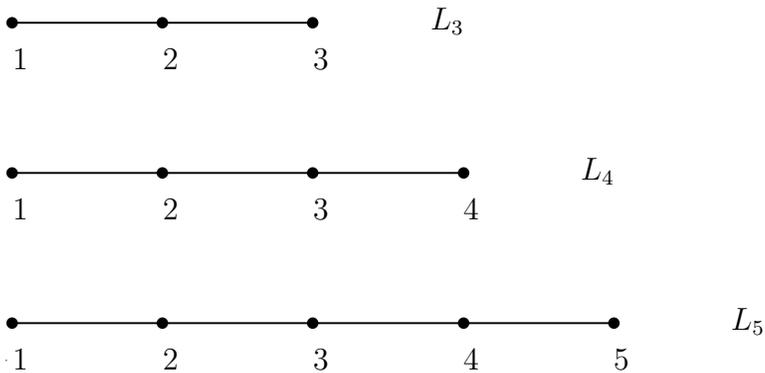
\label{fig}
\begin{texdraw}
  \drawdim cm
  \textref h:R v:C
  \move(0 0) \fcir f:0 r:0.075 \lvec(2 0)\move(0.2 -.5)\htext{$1$}
  \move(2 0) \fcir f:0 r:0.075 \lvec(4 0)\move(2.2 -.5)\htext{$2$}
  \move(4 0) \fcir f:0 r:0.075 \move(4.2 -.5)\htext{$3$} \move(6
  0) \htext{$L_3$}
  \move(0 -2) \fcir f:0 r:0.075 \lvec(2 -2)\move(0.2 -2.5)\htext{$1$}
  \move(2 -2) \fcir f:0 r:0.075 \lvec(4 -2)\move(2.2 -2.5)\htext{$2$}
  \move(4 -2) \fcir f:0 r:0.075 \lvec(6 -2)\move(4.2 -2.5)\htext{$3$}
  \move(6 -2) \fcir f:0 r:0.075 \move(6.2 -2.5)\htext{$4$}
  \move(8 -2) \htext{$L_4$}
  \move(0 -4) \fcir f:0 r:0.075 \lvec(2 -4)\move(0.2 -4.5)\htext{$1$}
  \move(2 -4) \fcir f:0 r:0.075 \lvec(4 -4)\move(2.2 -4.5)\htext{$2$}
  \move(4 -4) \fcir f:0 r:0.075 \lvec(6 -4)\move(4.2 -4.5)\htext{$3$}
  \move(6 -4) \fcir f:0 r:0.075 \lvec(8 -4)\move(6.2 -4.5)\htext{$4$}
  \move(8 -4) \fcir f:0 r:0.075 \move(8.2 -4.5)\htext{$5$}
  \move(10 -4) \htext{$L_5$}

\linewd 0.01 \lpatt(0.2 0.2)
  \lpatt(0.01 100)
\end{texdraw}
\caption{Line Graphs}
\end{figure}

These graphs are very similar to the cycles.\ie $L_n$ is the
subgraph of $C_n$ which has all the vertices of $C_n$ and all
edges except $\{1,n\}$. This allows us to apply much of what we
have deduced about the Betti numbers of cycles to find those of
Lines.

We use the formula from Theorem (\ref{ER2})
$$\b_i(L_n)=\sum _{H \subset G} {\dim}_k \Hred_{i-2}(\e(H);k)$$
where the sum is over all induced subgraphs of $L_n$ which contain
at least one edge. All of the induced subgraphs of $L_n$
(including $L_n$ itself) are collections of runs, much as in the
case for cycles.

\begin{thm}{\Prop}
The value of $\b_{i,2i}(L_n)$ is the number of induced subgraphs
of $L_n$ which are $i$ runs of length $2$. It follows that
$$\b_{i,2i}(L_n)={ {n-2i+1}\choose i}.$$
\end{thm}

\begin{proof}
Exactly the same reasoning for cycles from Proposition
(\ref{betai2i}) shows the first part. Using the notation of
definition (\ref{defC}) the number of induced subgraphs which are
$i$ runs of length $2$ is $C(n+2,2,i)$ and by Lemma (\ref{Cnlm})
this is ${ {n-2i+1}\choose i}$.
\end{proof}

\begin{thm}{Theorem}\label{betlines}
The non-zero Betti numbers of the line on $n$ vertices are, for
$i>j$,
$$\b_{i+j,2i+j}(L_n)={i \choose j}{ {n-2i} \choose i}+{{i-1}
\choose j}{{n-2i} \choose {i-1}}. $$

Or equivalently for $2d \geqslant l$
$$\b_{l,d}(L_n)= {{d-l}\choose {2l-d}}{{n-2d+2l} \choose {d-l}}+{{d-l-1}\choose {2l-d}}
{{n-2d+2l}\choose {d-l-1}}. $$

\end{thm}

\begin{proof}
This is similar to Proposition (\ref{bet1}). As in the case of
cycles we need to count the number of induced subgraphs of $L_n$
which are runs of lengths $0$ and $2 \mod 3$, \ie made of runs of
lengths $$3p_1,\dots,3p_{\a},3q_1+2,\dots,3q_{\b}+2$$ such that
$$
\left\{
\begin{array}{lll}
i+j &=& 2(P+Q) +\b \\
2i+j &=& 3(P+Q) +2\b
\end{array}
\right.
$$
with $P=\sum_{k=1}^{\a}p_k$ and $Q=\sum_{k=1}^{\b}q_k$.

 We may construct any induced subgraph of $L_n$ starting from
an induced subgraph which is $i$ runs of length $2$. Suppose that
$A$ is such a subgraph of $L_n$ comprised of $i$ runs of length
$2$ which does not feature the vertex $1$. We can select any $j$
of the $i$ vertices which are adjacent to, and left of, the runs
to construct a new induced subgraph of $L_n$. This implies there
are $${i \choose j} C(n+1,2,i)={i \choose j}{ {n-2i} \choose i}$$
subgraphs of this type contributing to $\b_{i+j,2i+j}(L_n)$. The
remaining subgraphs to count are those which contain vertices $1$
and $2$ but do not contain the vertex $3$. The number of such
subgraphs which are $i$ runs of length $2$ is $C(n-1,2,i-1)={
{n-2i} \choose {i-1}}$. We can select any $j$ of the $i-1$
vertices adjacent and left of the $i-1$ runs of length $2$ (which
don't include the run on vertices $1$ and $2$) to get ${{i-1}
\choose j}{{n-2i} \choose {i-1}}$ possible induced subgraphs.

Altogether we obtain
\begin{eqnarray*}
\b_{i+j,2i+j}(L_n)&=&{i \choose j}{ {n-2i} \choose i}+{{i-1}
\choose j}{{n-2i} \choose {i-1}}.\\
& &\\
\end{eqnarray*}
We can find $\b_{l,d}(L_n)$ as given in the statement of the
theorem by setting $i=d-l$ and $j=2l-d$.
\end{proof}

\begin{thm}{\Cor}\label{linedim}
The projective dimension of the line graph is independent of the
characteristic of the chosen field and is

$$\pd (L_n)=\left\{
{\begin{array}{lll}
\frac{2n}{3} & if & n \equiv 0 \mod 3 \\
\frac{2n-2}{3} & if & n \equiv 1 \mod 3\\
\frac{2n-1}{3} & if & n \equiv 2 \mod 3.
\end{array}} \right.
 $$

\end{thm}

\begin{proof}
We first consider the case $n=3m$ for some $m$. From Theorem
(\ref{betlines}) the $2m$th Betti number of degree $3m$ is
\begin{eqnarray*}
\b_{2m,3m}(L_{3m})&=&{ {m} \choose {m}}{{m} \choose {m}}+{ {m-1}
\choose {m}}{{m} \choose {m-1}}\\
&=&1.
\end{eqnarray*}
 Hence $\pd(L_n) \geqslant 2m$. Suppose now that $l > 2m$ If
 $\b_{l,d}(L_n) \neq 0$ then Theorem (\ref{betlines}) implies that $d-l
 \geqslant 2l-d$ , \ie that $2d \geqslant 3l$. As we are supposing
 that $l >2m$ this implies $d>3m=n$ which is not possible.
 Therefore $\b_{l,d}(L_n)=0$ for all $ l >2m$ and so
 $\pd(L_{3m})=2m=\frac{2n}{3}.$

 Now suppose that $n=3m+1$ for some $m$. From Theorem
(\ref{betlines}) the $2m$th Betti number of degree $3m$ is
\begin{eqnarray*}
\b_{2m,3m}(L_{3m+1})&=&{ {m} \choose {m}}{{m+1} \choose {m}}+{
{m-1}
\choose {m}}{{m+1} \choose {m-1}}\\
&=&m+1 \\
&\neq& 0.
\end{eqnarray*}
 Hence $\pd(L_n) \geqslant 2m$. We now suppose that $l > 2m$. If
 $\b_{l,d}(L_n) \neq 0$ then Theorem (\ref{betlines}) implies that
$d-l
 \geqslant 2l-d$ , \ie that $2d \geqslant 3l$. We are assuming that $l \geqslant 2m+1$
 which implies that $2d \geqslant 3l \geqslant 3(2m+1)=6m+3$ and
 so $d \geqslant 3m+ \frac{3}{2}$. Because $d$ is an integer this
 would mean that $d \geqslant 3m +2>n$ which is not possible.
 Therefore $\b_{l,d}(L_n)=0$ for all $l > 2m$ so $\pd(L_n)=2m= \frac{2n-2}{3}.$

 Finally we suppose that $n=3m +2$ for some integer $m$. From
 Theorem (\ref{betlines}) the $(2m+1)$th Betti number of degree
 $3m+2$ is
\begin{eqnarray*}
 \b_{2m+1,3m+2}(L_{3m+2})&=&{ {m+1} \choose {m}}{{m} \choose {m+1}}+{
{m}
\choose {m}}{{m} \choose {m}}\\
&=&1.
\end{eqnarray*}
Hence $\pd(L_n) \geqslant 2m+1$. Suppose now that $l \geqslant
2m+2.$ If $\b_{l,d}(L_n) \neq 0$ then we must have $d-l \geqslant
2l-d$, \ie $2d \geqslant 3l \geqslant 3(2m+2) =6m+6$. This implies
that $d \geqslant 3m+3 >n.$ which is not possible. Hence
$\b_{l,d}(L_n)=0$ for all $l>2m+1$ and so
$\pd(L_n)=2m+1=\frac{2n-1}{3}$.

\end{proof}

\chapter{Betti Numbers of Graphs of Degree 2}
\setcounter{satz}{0}

We have found the Betti numbers of cycles and lines and that they
are independent of our choice of field. We now show that the Betti
numbers of any graph which has vertices only of degree at most 2
are independent of the field.

To show this we use Theorem (\ref{ER2}) which gives the Betti
numbers in terms of induced subgraphs $$\b_{i,d}(G)=\sum _{H
\subset G, |V(H)|=d} {\dim}_k \Hred_{i-2}(\e(H);k),$$ where the
sum is over all induced subgraphs of $G$ which contain at least
one edge. If a graph is such that $\dim_k\Hred_i(\e(H);k)$ does
not depend on $k$ for all $i$ and all induced subgraphs which
contain an edge $H$ then the Betti numbers must also be
independent of $k$. If a graph is such that all its vertices have
degree at most 2 then this condition is also true of all its
induced subgraphs. Therefore we need to show that the dimension of
$\Hred_i(\e(G);k)$ is independent of $k$ for all $i$ and all
graphs which only have vertices of degree 2 or less.

Such a graph must be the disjoint union of cycles, lines and
isolated vertices. If the graph features any isolated vertices
then $\Hred_i(\e(G);k)=0$ for all $i$ and so has dimension which
is independent of $k$. Therefore we may assume that the graph is a
disjoint union of cycles and lines. To prove such a graph has
Betti numbers independent of choice of field we require the
following sequence of lemmas.

\section{Reduced Homology of Unions of Lines}

Throughout this chapter $L_n$ will denote the line graph on $n$
vertices where $n \geqslant 2$. If $G$ and $H$ are graphs we will
denote their disjoint union by $G \cup H$.

\begin{thm}{\Lem}\label{lemma0}
For a graph $G$ and $n \geqslant 5$ there is an isomorphism of
reduced homology modules
$$\Hred_i(\e(L_n \cup G);k) \cong \Hred_{i-2}(\e(L_{n-3} \cup
G);k)$$
\end{thm}

\begin{proof}
Let $E=\e(L_n \cup
G)=\e(\{1,2\},\{2,3\},\dots.\{n-1,n\},e_1,\dots,e_r;V)$ where
$V=V(L_n \cup G)$ and $e_1,\dots,e_r$ are the edges of $G$. We can
write $E=E_1 \cup E_2$ where $E_1=\e( \{1,2\};V)$ and
$E_2=\e(\{2,3\},\dots.\{n-1,n\},e_1,\dots,e_r;V)$. We now use the
Mayer-Vietoris sequence

\begin{eqnarray*}\cdots \to
\Hred_i(E_1;k) \oplus
\Hred_i(E_2;k)\to \Hred_i(E;k)&\to&\Hred_{i-1}(E_1 \cap E_2;k)\to \\
\to \Hred_{i-1}(E_1;k) \oplus \Hred_{i-1}(E_2;k)\to\cdots
\end{eqnarray*}

 We
have $\Hred_i(E_1;k)=0$ for all $i$ because $E_1$ is a simplex and
$\Hred_i(E_2;k)=0$ for all $i$ by Lemma (\ref{lemA}) because $1
\in V \setminus (\{2,\dots,n\} \cup \bigcup_{i=2}^re_r)$. Hence we
obtain from the above sequence $\Hred_i(E;k) \cong \Hred_{i-1}(E_1
\cap E_2;k)$ for all $i$. We now examine $E_1 \cap E_2$. By Lemma
(\ref{lemB})
\begin{eqnarray*}
E_1 \cap E_2 &=&\e( \{1,2\} \cup \{2,3\},\{1,2\}\cup
\{3,4\},\dots\\
& &\dots \{1,2\} \cup \{n-1,n\},\{1,2\} \cup e_1,\dots,
\{1,2\}\cup
e_r;V)\\
&=&
\e(\left\{1,2,3\right\},\left\{1,2,3,4\right\},\left\{1,2,4,5\right\},\dots\\
& &\dots \{1,2,n-1,n\},\left\{1,2\right\}\cup
e_1,\dots,\left\{1,2\right\}\cup
e_r;V)\\
&=& \e(\left\{3\right\},\left\{3,4\right\},\left\{4,5\right\},\dots,\{n-1,n\},e_1,\dots,e_r;V')\\
&\ & (\rm{using \ Lemma} \ (\ref{lemC}) \ \rm{as} \  \{3\} \subseteq\{3,4\} ) \\
&=&
\e(\left\{3\right\},\left\{4,5\right\},\dots,\{n-1,n\},e_1,\dots,e_r;V')
\end{eqnarray*}
where $V'=V \setminus \{1,2\}$. Because $3 \notin
(\{4,5,\dots,n-1,n\} \cup \bigcup_{i=1}^r e_r)$ we can use Lemma
(\ref{lemD}) to see that
$\Hred_i(\e(\left\{3\right\},\left\{4,5\right\},\dots,\{n-1,n\},e_1,\dots,e_r;V');k)
=\Hred_{i-1}(\e(\left\{4,5\right\},\dots,\{n-1,n\},e_1,\dots,e_r;V'');k)$
where $V''=V \setminus \{1,2,3\}$. Putting these results together
\begin{eqnarray*}
\Hred_i(E;k)&\cong&\Hred_{i-1}(\e(\left\{3\right\},\left\{4,5\right\},\dots,\{n-1,n\},e_1,\dots,e_r;V');k)\\
          &\cong&\Hred_{i-2}(\e(\left\{4,5\right\},\dots,\{n-1,n\},e_1,\dots,e_r;V'');k)\\
          &\cong&\Hred_{i-2}(\e(L_{n-3} \cup G);k).
\end{eqnarray*}
\end{proof}

\begin{thm}{\Lem}\label{lemm0.5}
There are isomorphisms of reduced homology modules\\
(i) $\Hred_i(\e(L_2 \cup G);k) \cong \Hred_{i-1}(\e(G);k),$\\
(ii) $\Hred_i(\e(L_3 \cup G);k) \cong \Hred_{i-2}(\e(G);k),$\\
(iii) $\Hred_i(\e(L_4 \cup G);k) \cong  0.$\\
\end{thm}

\begin{proof}
$(i)$. Let $E=\e(L_2 \cup G)=E_1 \cup E_2$ where
$E_1=\e(\{1,2\};V)$,  $E_2=\e(e_1,\dots,e_r;V)$,  $e_1,\dots,e_r$
are the edges of $G$ and $V=V(L_2 \cup G)$. Note that
$\Hred_i(E_1;k)=0$ for all $i$ as $E_1$ is a simplex and
$\Hred_i(E_2;k)=0$ for all $i$ by (\ref{lemA}) as $1 \in V
\setminus \bigcup_{i=1}^r e_i.$ The Mayer-Vietoris sequence,
\begin{eqnarray*}\cdots \to
\Hred_i(E_1;k) \oplus
\Hred_i(E_2;k)\to \Hred_i(E;k)&\to&\Hred_{i-1}(E_1 \cap E_2;k)\to \\
\to \Hred_{i-1}(E_1;k) \oplus \Hred_{i-1}(E_2;k)\to\cdots
\end{eqnarray*}
 thus provides the isomorphism
$\Hred_i(E;k) \cong \Hred_{i-1}(E_1 \cap E_2;k)$ for all $i$. By
Lemma (\ref{lemB}) we can write
\begin{eqnarray*}
E_1 \cap E_2 &=& \e(\{1,2\} \cup e_1, \dots, \{1,2\} \cup e_r;V) \\
             &=& \e(e_1, \dots,e_r;V \setminus \{1,2\}) \\
             &=& \e(G)
\end{eqnarray*}
from which we conclude $\Hred_i(\e(L_2 \cup G);k) =
\Hred_{i-1}(\e(G);k)$ for all $i$.

$(ii)$ Now we let $E=\e(L_3 \cup G)=E_1 \cup E_2$ where
$E_1=\e(\{1,2\};V)$, \\ $E_2=\e(\{2,3\},e_1,\dots,e_r;V)$,
$e_1,\dots,e_r$ are the edges of $G$ and $V=V(L_3 \cup G)$. Again
we have $\Hred_i(E_1;k)=\Hred_i(E_2;k)=0$ for all $i$, so the
Mayer-Vietoris sequence implies there is an isomorphism
$\Hred_i(E;k) \cong \Hred_{i-1}(E_1 \cap E_2;k)$. By Lemma
(\ref{lemB}) we write the simplicial complex $E_1 \cap E_2$ as
\begin{eqnarray*}
E_1 \cap E_2 &=&\e(\{1,2,3\},\{1,2\} \cup e_1,\dots,\{1,2\} \cup
e_r;V) \\
&=&\e(\{3\},e_1,\dots,e_r;V \setminus \{1,2\}) \\
\end{eqnarray*}
By Lemma (\ref{lemD})
\begin{eqnarray*}
\Hred_i(E_1 \cap E_2;k) &=& \Hred_{i-1}(\e(e_1,\dots,e_r;V
\setminus
\{1,2,3\};k) \\
&=& \Hred_{i-1}(\e(G);k)
\end{eqnarray*}
Hence $\Hred_i(E;k)=\Hred_{i-1}(E_1 \cap
E_2;k)=\Hred_{i-2}(\e(G);k).$

$(iii)$ Suppose that
$$E=\e(L_4 \cup G)=\e(\left\{1,2\right\},\left\{2,3\right\},\left\{3,4\right\},e_1,\dots,e_r;V).$$
where $e_1,\dots,e_r$ are the edges of $G$ and $V=V(L_4 \cup G)$.

 Let $E_1=\e(\left\{1,2\right\};V)$ and
$E_2=\e(\left\{2,3\right\},\left\{3,4\right\},e_1,\dots,e_r;V)$.
For similar reasons to the previous parts we have
$\Hred_{i}(E_1;k)=\Hred_{i}(E_2;k)=0$ for all $i$ and so $\Hred
_i(E;k)=\Hred_{i-1}(E_1\cap E_2;k)$ for all $i$. Using Lemma
(\ref{lemB}) we see that
\begin{eqnarray*}
E_1 \cap E_2 &=& \e(\left\{1,2,3\right\},\left\{1,2,3,4\right\},
\left\{1,2\right\}\cup e_1,\dots,\left\{1,2\right\} \cup e_r;V)\\
&=& \e(\left\{3\right\},\left\{3,4\right\}, e_1,\dots,\ e_r;V')\\
&=& \e(\left\{3\right\},e_1,\dots,\ e_r;V')\\
&\ &(\rm{using \ Lemma\ } (\ref{lemC}))
\end{eqnarray*}
where $V'=V \setminus \left\{1,2\right\}$. Since $4 \in V'
\setminus ( \{3 \} \cup \bigcup_{i=1}^r e_i)$ we can use Lemma
(\ref{lemA}) to see that $\Hred_i(\e(\left\{3\right\},e_4,\dots,\
e_t;V');k)=0$ for all $i$. Therefore we also have $\Hred_i(E;k)=0$
for all $i$.

\end{proof}

\begin{thm}{\Lem}\label{lemma1}
The reduced homology of the disjoint union of the line $L_n$ and
any graph $G$ may be expressed
$$\Hred_i(\e(L_n \cup G);k)=   \left\{ \begin{array}{lll}
                                \Hred_{i-\frac{2n}{3}}(\e(G);k)   &if& n \equiv0 \mod 3 \\
                                0                                &if& n \equiv 1 \mod 3 \\
                                \Hred_{i-\frac{2n-1}{3}}(\e(G);k) &if& n \equiv 2 \mod 3
                                \end{array} \right  .$$
\end{thm}

\begin{proof}
First we suppose that $n=3m$ for some integer $m$. By Lemma
(\ref{lemma0})
\begin{eqnarray*}
\Hred_i(\e(L_n \cup G);k) &=& \Hred_{i-2}(\e(L_{n-3} \cup G);k)\\
                        &=& \Hred_{i-2(m-1)}(\e(L_{3} \cup G);k)\\
                        &=&\Hred_{i-2m}(\e(G);k)\\
                        & &({\mathrm{By\  Lemma \ (\ref{lemm0.5}))}}\\
                        &=&\Hred_{i-\frac{2n}{3}}(\e(G);k).\\
\end{eqnarray*}

Now suppose that $n=3m+1$. By Lemma (\ref{lemma0})
\begin{eqnarray*}
\Hred_i(\e(L_n \cup G);k) &=& \Hred_{i-2}(\e(L_{n-3} \cup G);k)\\
                        &=& \Hred_{i-2(m-1)}(\e(L_{4} \cup G);k)\\
                        &=&0 \\
                        & &({\mathrm{By\  lemma \
                        (\ref{lemm0.5})}})
\end{eqnarray*}

Finally we have the case $n=3m+2$. By Lemma (\ref{lemma0})
\begin{eqnarray*}
\Hred_i(\e(L_n \cup G);k) &=& \Hred_{i-2}(\e(L_{n-3} \cup G);k)\\
                        &=& \Hred_{i-2(m)}(\e(L_{2} \cup G);k)\\
                        &=& \Hred_{i-(2m+1)}(\e(G);k)\\
                        & &({\mathrm{By\  Lemma \
                        (\ref{lemm0.5})}})\\
                        &=& \Hred_{i-\frac{2n-1}{3}}(\e(G);k)\\
\end{eqnarray*}
\end{proof}

\begin{thm}{\Cor}\label{cor0}
If, for a graph $H$, $\dim_k\Hred_i(\e(H);k)$ is independent of
$k$ for all i then so is $\dim_k\Hred_i(\e(H\cup L_n);k)$. Hence
$\dim_k(\Hred_i(\e(H \cup L_{n_1} \cup \cdots \cup L_{n_s});k)$ is
also independent of $k$.
\end{thm}

\begin{proof}
This follows from Lemma (\ref{lemma1}) as $\Hred_i(\e(H\cup
L_n);k)$ is either $0$ or a shift of $\Hred_i(\e(H);k)$.
\end{proof}

\begin{thm}{\Cor} \label{cor01}
If $G$ is the disjoint union of lines then the dimension of the
vector space $\Hred_i(\e(G);k)$ is independent of $k$.
\end{thm}

\begin{proof}
This is Corollary (\ref{cor0}) when $H$ is a line. We know that
$\dim_k(\Hred_i(L_n);k)$ is independent of $k$ from Theorem
(\ref{betlines}).
\end{proof}

\section{Reduced Homology of Unions of Cycles}

\begin{thm}{\Lem}\label{lemma2}
Let $H$ be any graph and let $G=C_n \cup H$. There is a long exact
sequence
$$\cdots\to \Hred_i(\e(L_n \cup H);k)\to
\Hred_i(\e(G);k)\to\Hred_{i-3}(\e(L_{n-4} \cup H);k)\to$$ $$\to
\Hred_{i-1}(\e(L_n \cup H);k)\to
\Hred_{i-1}(\e(G);k)\to\Hred_{i-4}(\e(L_{n-4} \cup
H);k)\to\cdots$$
\end{thm}

\begin{proof}
Let the vertex set of $C_n$ be $\{1,\dots,n\}$ and let the edges
be $$\{1,2\},\{2,3\},\dots,\{n-1,n\},\{1,n\},$$ Suppose the edges
of $H$ are $e_1,\dots,e_r$. We may write $\e(G)=\e_1 \cup \e_2$
where $\e_1=\e(\{1,n\};V(G))$ and
$$\e_2=\e(\{1,2\},\dots,\{n-1,n\},e_1,\dots,e_r;V(G))\simeq \e(L_n
\cup H).$$ Note that $\Hred_i(\e_1;k)=0$ for all $i$ as $\e_1$ is
a simplex. Using Lemmas (\ref{lemB}) and (\ref{lemC}) we see that
the intersection of $\e_1$ and $\e_2$ is
$$\e_1 \cap \e_2=\e(\{2\},\{3,4\},\dots,\{n-3,n-2\},\{n-1\},e_1,\dots,e_r;V(G)\setminus\{1,n\}).$$
By Corollary (\ref{CorE})
\begin{eqnarray*}
\Hred_i(\e_1 \cap
    \e_2)&=&\Hred_{i-2}(\e(\{3,4\},\dots,\{n-3,n-2\},e_1,\dots,e_r;V');k)\\
      &=&\Hred_{i-2}(\e(L_{n-4} \cup H);k),
\end{eqnarray*}
where $V'=V(G)\setminus\{1,2,n-1,n\}$. Now we use the
Mayer-Vietoris sequence
$$\cdots\to \Hred_i(\e_1;k) \oplus \Hred_i(\e_2;k)\to
\Hred_i(\e(G);k)\to\Hred_{i-1}(\e_1 \cap \e_2);k)\to$$
$$\to \Hred_{i-1}(\e_1;k) \oplus \Hred_{i-1}(\e_2;k)\to
\Hred_{i-1}(\e(G);k)\to\Hred_{i-2}(\e_1 \cap \e_2);k)\to\cdots$$

which becomes
$$\cdots\to \Hred_i(\e(L_n \cup H);k)\to
\Hred_i(\e(G);k)\to\Hred_{i-3}(\e(L_{n-4} \cup H);k)\to$$
$$\to \Hred_{i-1}(\e(L_n \cup H);k)\to
\Hred_{i-1}(\e(G);k)\to\Hred_{i-4}(\e(L_{n-4} \cup
H);k)\to\cdots$$
\end{proof}

\begin{thm}{\Lem}\label{lemma3}
Let $H$ be a graph and let $G=C_n \cup H$.
$$\Hred_i(\e(G);k)=   \left\{ \begin{array}{lll}
                               \Hred_{i-\frac{2n+1}{3}}(\e(H);k)   &if& n \equiv 1 \mod 3 \\
                               \Hred_{i-\frac{2n-1}{3}}(\e(H);k)   &if& n \equiv 2 \mod 3
                                \end{array} \right  .$$

\end{thm}

\begin{proof}
(i) Suppose that $n \equiv 1 \mod 3$. From Lemma (\ref{lemma2}) we
have the long exact sequence
$$\cdots\to \Hred_i(\e(L_n \cup H);k)\to
\Hred_i(\e(G);k)\to\Hred_{i-3}(\e(L_{n-4} \cup H);k)\to\cdots$$ We
use Lemma (\ref{lemma1}) to obtain $\Hred_i(\e(L_n \cup H);k)=0$
for all $i$ and $\Hred_{i-3}(\e(L_{n-4} \cup
H);k)=\Hred_{i-\frac{2n+1}{3}}(\e(H);k)$ for all $i$. Hence the
sequence tells us that
$\Hred_i(\e(G);k)=\Hred_{i-\frac{2n+1}{3}}(\e(H);k)$ for all $i$.

(ii) Now suppose that $n \equiv 2 \mod 3$. In this case Lemma
(\ref{lemma1}) is used to obtain $\Hred_i(\e(L_n \cup
H);k)=\Hred_{i-\frac{2n-1}{3}}(\e(H);k)$ for all $i$ and
$\Hred_{i-3}(\e(L_{n-4} \cup H);k)=0$ for all $i$. From the exact
sequence we deduce that
$\Hred_i(\e(G);k)=\Hred_{i-\frac{2n-1}{3}}(\e(H);k).$
\end{proof}

\begin{thm}{\Cor}\label{Cor1}
Let $H$ be a graph and let $G=C_n \cup H$. If $n \equiv 1 \mod 3$
or $n\equiv 2 \mod 3$ and $\dim_k \Hred_i(\e(H);k)$ is independent
of $k$ for all $i$ then $\dim_k \Hred_i(\e(G);k)$ is also
independent of $k$ for all $i$.

Furthermore, if $G=C_{n_1}\cup \cdots\cup C_{n_s} \cup H$, where
none of $n_1,\dots,n_s$ are divisible by $3$, and $\dim_k
\Hred_i(\e(H);k)$ is independent of $k$ for all $i$ then $\dim_k
\Hred_i(\e(G);k)$ is also independent of $k$ for all $i$.

\end{thm}

\begin{proof}
This follows from Lemma (\ref{lemma3}) as $\dim_k
\Hred_i(\e(G);k)$ is a shift of $\dim_k \Hred_i(\e(H);k)$ in both
cases.
\end{proof}

\begin{thm}{\Cor}\label{cor12}
If $G$ is the disjoint union of cycles, each with number of
vertices not divisible by $3$, then $\dim_k(\Hred_i(\e(G);k))$ is
independent of $k$.
\end{thm}

\begin{proof}
This is Corollary (\ref{Cor1}) when $H$ is a cycle. By Theorem
(\ref{betticycle}) $\dim_k(\Hred_i(\e(C_n);k)$ is independent of
$k$.
\end{proof}

\begin{thm}{\Cor}\label{cor13}
If $G$ is the disjoint union of lines and cycles (which all have
number of vertices not divisible by $3$) then $\dim_k
\Hred_i(\e(G);k))$ is independent of $k$.
\end{thm}

\begin{proof}
This is Corollary (\ref{Cor1}) when $H$ is a union of lines. By
Corollary (\ref{cor01}) $\dim_k \Hred_i(\e(H);k)$ is independent
of $k$.
\end{proof}

Let $G$ be any graph which is the disjoint union of lines and
cycles. Let $G=H \cup L_{n_1} \cup \dots \cup L_{n_s}$, where $H$
is a disjoint union of cycles, then by Corollary (\ref{cor0}) if
$\dim_k \Hred_i(\e(H);k)$ is independent of $k$ then so is $\dim_k
\Hred_i(\e(G);k)$. Suppose that $H=J \cup C_{m_1} \cup \dots \cup
C_{m_t}$, where none of $m_1,\dots,m_t$ are divisible by 3 and $J$
is the disjoint union of cycles each with number of vertices
divisible by 3. By Corollary (\ref{Cor1}) if $\dim_k
\Hred_i(\e(J);k)$ is independent of $k$ then so is $\dim_k
\Hred_i(\e(H);k)$. If we show that $\dim_k \Hred_i(\e(J);k)$ is
independent of $k$ it will follow that $\dim_k \Hred_i(\e(G);k)$
is independent of $k$. Therefore the remaining case to consider is
when $G$ is the union of cycles all of which have number of
vertices divisible by 3.

We first consider the union of two such cycles.

\begin{thm}{\Lem} \label{lemma4}
Let $G=C_n \cup C_m$ where $n=3p$ and $m=3q$ for $p$ and $q$
integers. Then
\begin{eqnarray*}
\Hred_i(\e(C_n \cup C_m);k)=\left\{ \begin{array}{lll}
                                      k^4 &if& i=2(p+q-1)\\
                                      0   &if& i \neq 2(p+q-1).
                                      \end{array}
                                      \right.
\end{eqnarray*}
\end{thm}

\begin{proof}
The exact sequence of Lemma (\ref{lemma2})
$$\cdots\to \Hred_i(\e(L_n \cup C_m);k)\to
\Hred_i(\e(G);k)\to\Hred_{i-3}(\e(L_{n-4} \cup C_m);k)\to\cdots$$
when used in conjunction with Lemma (\ref{lemma1}) becomes
$$\cdots\to \Hred_{i-2p}(\e(C_m);k)\to
\Hred_i(\e(G);k)\to\Hred_{i-2p}(\e(C_m);k)\to\cdots$$ From the
proof of Proposition (\ref{prop2}) we know that
$$ \Hred_i(\e(C_m;k))= \left\{ \begin{array}{lll}
                            k^2 &if& i=2q-2\\
                            0   &if& i \neq 2q-2.
                            \end{array} \right.
                            $$
So the above exact sequence becomes zero everywhere except for the
section
$$\cdots \to 0 \to k^2 \to \Hred_{2p+2q-2}(\e(G);k) \to k^2 \to 0
\to \cdots$$ from which the result follows.
\end{proof}

\begin{thm}{\Lem}\label{lemma5}
Let $n=3p$ for some integer $p$. Suppose that $r>0$ and $H$ is a
graph such that
$$\dim_k\Hred_{j}(\e(H);k)=\left\{ \begin{array}{lll}
                                    r &if& j=i-2p\\
                                    0 &if& j\neq i-2p.
                                    \end{array}
                                    \right.
                                    $$
                                    Then
$$\dim_k\Hred_{j}(\e(C_n \cup H);k)=\left\{ \begin{array}{lll}
                                    2r &if& j=i\\
                                    0 &if& j\neq i.
                                    \end{array}
                                    \right.
                                    $$

\end{thm}

\begin{proof}
Again we use the exact sequence of Lemma (\ref{lemma2})
$$\cdots\to \Hred_i(\e(L_n \cup H);k)\to
\Hred_i(\e(C_n \cup H);k)\to\Hred_{i-3}(\e(L_{n-4} \cup
H);k)\to\cdots$$ Using Lemma (\ref{lemma1}) this becomes

$$\cdots\to \Hred_{i-2p}(\e(H);k)\to
\Hred_i(\e(C_n \cup H);k)\to\Hred_{i-2p}(\e(H);k)\to\cdots$$

The facts about the reduced homology of $\e(H)$ in the statement
of the Lemma mean that the sequence is zero everywhere except for
the section
$$\cdots\to 0\to k^r \to \Hred_i(\e(C_n \cup H);k)\to k^r \to 0
\to \cdots $$ from which the result follows.
\end{proof}

\begin{thm}{\Cor}\label{Cor2}
Let $G=C_{3p_1}\cup \cdots \cup C_{3p_r}$ for integers
$p_1,\dots,p_r$. Then
$$ \dim_k\Hred_i(\e(G);k)=\left\{ \begin{array}{lll}
                            2^r &if& i=2(p_1+\cdots+p_r)-2\\
                            0 &if& i\neq 2(p_1+\cdots+p_r)-2.\\
                            \end{array}
                            \right.
                            $$
Consequently $\dim_k\Hred_i(\e(G);k)$ is independent of $k$.
\end{thm}

\begin{proof}
This follows from repeated application of Lemma (\ref{lemma5}) and
the fact that, for an integer $q$,
$$ \dim_k \Hred_i(\e(C_{3q};k))= \left\{ \begin{array}{lll}
                            2 &if& i=2q-2\\
                            0   &if& i \neq 2q-2.
                            \end{array} \right.
                            $$
from the proof of Proposition (\ref{prop2}).
\end{proof}

\section{Betti Numbers of Graphs of Degree 2}

The above lemmas allow us to conclude the following.

\begin{thm}{Theorem}
The Betti numbers of graphs which only have vertices of degree $2$
or less are independent of choice of field.
\end{thm}

\begin{proof}
To show that $$\b_{i,d}(G)=\sum _{H \subset G, |V(H)|=d} {\dim}_k
\Hred_{i-2}(\e(H);k)$$ is independent of $k$ for all $i$ for this
family of graphs we need to show that ${\dim}_k
\Hred_{i}(\e(G);k)$ is independent of $k$ for all graphs which are
disjoint unions of cycles and lines. By Corollary (\ref{cor13})
this is proven whenever none of the cycles have number of vertices
divisible by 3 and by Corollary (\ref{Cor2}) whenever the graphs
are unions only of cycles all of which have number of vertices
divisible by 3. Otherwise we use Corollary (\ref{cor0}) to show
that we only need prove the result for unions of cycles which we
do by using Corollary (\ref{Cor1}) when $H$ is the subgraph which
features all the cycles with number of vertices divisible by 3.
\end{proof}

\chapter{Betti Numbers of Forests}
\setcounter{satz}{0}

In this chapter we are concerned with the Betti numbers and
projective dimensions of forests. In some sense forests have less
rigid structure in general than the other families of graphs, such
as cycles and complete graphs, we have studied. Because of this we
cannot find an explicit formula for Betti numbers of forests in
the way we have for other families. Instead we find a way of
finding the Betti numbers of a forest in terms of the Betti
numbers of some of its subgraphs (which must always be forests
too). This provides us with an algorithm for finding the Betti
numbers by repeated application of the formula to successively
smaller subgraphs. This will also allow us to make statements
about the projective dimension of forests.

\section{Trees and Forests}

Trees and forests were defined in $(\ref{definitions})$. Note that
as forests have no cycles they must have terminal vertices.

\begin{thm}{\Prop} Every forest has a vertex, v, such that either
\begin{enumerate}

\item all but perhaps one of its neighbours have degree $1$ \item
if v is of degree $1$ then its unique neighbour is also of degree
$1$.
\end{enumerate}
\end{thm}

\begin{proof}
Let $T$ be a forest. We use induction on the number of vertices of
$T$. If $T$ has one vertex then the statement is true as this
vertex has no neighbours. Assume the statement of the proposition
is true for any forest with $m$ vertices. Now suppose that
$|V(T)|=m+1$.

If $T$ has an isolated vertex then we are done. Otherwise take a
connected component of $T$. This must be a tree and so will have a
vertex of degree 1. Call this vertex $w$ and let $H=T \setminus
\{w\}$. As $H$ is a forest with $m$ vertices we can select a
vertex, $u$ say, which has neighbours $u_1,\dots,u_n$ such that
$\deg u_1 = \cdots =\deg u_{n-1}=1$ (and if $n=1$ then $\deg
u_1=1$).

In $T$ we must have one of the following
\begin{enumerate}
\item $w$ has one of $u_1,\dots,u_{n-1}$ as its unique neighbour,
\item one of the edges of $T$ is $\{u,w\}$, \item $w$ is not a
neighbour of any of $u,u_1,\dots,u_{n-1}$.
\end{enumerate}
In case (i) $w$ has unique neighbour  $u_i$, for some $u_i \in
\{u_1,\dots,u_{n-1}\}.$ Therefore all but at most one of the
neighbours of $u_i$ has degree 1. This is because $u_i$ has two
neighbours, $w$ which has degree 1 and $u$. So the statement is
proved in this case.

For case (ii) $\{u,w\} \in E(T)$ and all but at most one of the
neighbours of $u$ have degree 1.

In case (iii) $w$ is not a neighbour of any of
$u,u_1,\dots,u_{n-1}$ and so all but perhaps one of the neighbours
of $u$ have degree 1.
\end{proof}

\begin{figure}[h]\label{tree}
\begin{texdraw}
  \drawdim cm
  \textref h:R v:C
  \move(0 3) \fcir f:0 r:0.08  \htext{$v_1\ $} \lvec( 1.5 1.5)
  \move(0 2) \fcir f:0 r:0.08  \htext{$v_2\ $}  \lvec( 1.5 1.5)
  \move(0 0) \fcir f:0 r:0.08  \htext{$v_{n-1}\ $}  \lvec( 1.5 1.5)
  \textref h:C v:D
  \move(1.5 1.5) \fcir f:0 r:0.08  \move(1.5 1.1) \htext{$v$}  \move(1.5 1.5) \lvec( 3 1.5)
  \move(3 1.5) \fcir f:0 r:0.08  \move(3 1.1) \htext{$v_n$}
  \move(5.250000000 2.799038106) \fcir f:0 r:0.08   \lvec( 3 1.5) \move(5.250000000 3) \htext{$w_m$}
  \move(4.700961894 2.249999999) \fcir f:0 r:0.08   \lvec( 3 1.5) \move(5.3 2) \htext{$w_{m-1}$}
  \move(5.249999999 .200961895) \fcir f:0 r:0.08    \lvec( 3 1.5) \move(5.249999999 -.24)  \htext{$w_1$}
  \linewd 0.01 \lpatt(0.2 0.2) \move(6 1.5) \lcir r:1.5 \fcir f:0.95 r:1.5
  \lpatt(0.01 100)
  \move(5.250000000 2.799038106) \fcir f:0 r:0.08   \lvec( 3 1.5) \move(5.250000000 3) \htext{$w_m$}
  \move(4.700961894 2.249999999) \fcir f:0 r:0.08   \lvec( 3 1.5) \move(5.3 2) \htext{$w_{m-1}$}
  \move(5.249999999 .200961895) \fcir f:0 r:0.08    \lvec( 3 1.5) \move(5.249999999 -.24)  \htext{$w_1$}
  \move(6 1.2) \htext{H}
  \move(6.5 -1.5) \htext {Figure 9.1}
\end{texdraw}
\end{figure}

\begin{thm}[\rm]{Notation}
Throughout this chapter $v$ will denote a vertex of $T$ which has
all but at most one of its neighbours of degree 1 (and if it has
exactly one neighbour then that neighbour also has degree 1). The
neighbours of $v$ will be denoted $v_1,\dots,v_n$ such that
$v_1,\dots,v_{n-1}$ all have degree 1. Also the neighbours of
$v_n$ other than $v$ will be denoted by $w_1,\dots,w_m$. This is
illustrated in Figure 9.1.

In this chapter let $T$ denote a forest and let $T'$ denote the
subgraph of $T$ which is obtained by deleting the vertex $v_1$ and
let $T''$ denote the subgraph of $T$ which is obtained by deleting
the vertices $v,v_1,\dots,v_n$. That is, $T'=T \setminus \lc v_1
\rc$ and $T''=T \setminus \lc v,v_1,\dots,v_n \rc$. Note that $T'$
and $T''$ must both be forests.

\end{thm}

To find the $\mathbb{N}$-graded Betti numbers, in terms of the
$\mathbb{N}$-graded Betti numbers of subgraphs of $T$, we will
again use the formula of (\ref{ER}).

\section{{$\mathbb{N}$}-graded Betti numbers of Forests}

We will use the formula
$$ \beta_{i,d}(T)=\sum_{F \in \AD(T): \ |V(T) \setminus F|=d} \dim_k
\Hred_{i-2}(\linkad F;k). $$

It will be convenient to split this sum into several cases for
faces of the Alexander dual.

\begin{thm}[\rm]{\Def}\label{bet1234}
Let
\begin{eqnarray*}
\b_{i,d}^{(1)} &=& \sum_{F \in \AD(T): \ |V(T) \setminus F|=d, \
v_1 \in F}
\dim_k \Hred_{i-2}(\Link_{\AD(T)} F;k).\\
\b_{i,d}^{(2)} &=& \sum_{F \in \AD(T): \ |V(T) \setminus F|=d, \
v_1\notin F, \ v \in F} \dim_k \Hred_{i-2}(\Link_{\AD(T)}
F;k).\\
\b_{i,d}^{(3)} &=& \sum_{F \in \AD(T): \ |V(T) \setminus F|=d, \
v,v_1\notin
F, \ v_n \in F }\dim_k \Hred_{i-2}(\Link_{\AD(T)} F;k).\\
\b_{i,d}^{(4)} &=& \sum_{F \in \AD(T): \ |V(T) \setminus F|=d,\
v,v_1,v_n \notin F} \dim_k \Hred_{i-2}(\Link_{\AD(T)} F;k).
\end{eqnarray*}
Hence $\b_{i,d}(T)=\sum_{j=1}^4 \b_{i,d}^{(j)}.$
\end{thm}

\begin{thm}{\Lem}\label{lem1}
The contribution to $\b_{i,d}$ from those faces of $\AD(T)$ which
contain $v_1$ is the $i$th Betti number of degree $d$ of $T'$, \ie
$$\b_{i,d}^{(1)}= \sum_{F \in \AD(T): \ |V(T) \setminus F|=d, \ v_1 \in F} \dim_k
\Hred_{i-2}(\Link_{\AD(T)} F;k)=\beta_{i,d}(T'). $$
\end{thm}
\begin{proof}
The faces of $\AD (T)$ which contribute to $\b_{i,d}^{(1)}$ are
those which include the vertex $v_1$. These are exactly the faces
of $\AD(T)$ which have link excluding $v_1$. If $F \in \AD(T)$
does not include $v_1$ then it has the same link as the face $F
\setminus \left\{v_1 \right\}$ of $\AD(T \setminus \left\{v_1
\right\}) = \AD(T')$. Hence
\begin{eqnarray*}
\b_{i,d}^{(1)}&=& \sum_{F \in \AD(T): \ |V(T) \setminus F|=d, \
v_1 \in F}
\dim_k \Hred_{i-2}({\Link}_{\AD(T)} F;k)\\
&=&\sum_{F \in \AD(T'): \ |V(T') \setminus F|=d }\dim_k
\Hred_{i-2}({\Link}_{\AD(T')}
F;k)\\
&=&\b_{i,d}(T')
\end{eqnarray*}
for all $i$.
\end{proof}
\begin{thm}{\Lem}\label{lem2}
The faces of $\AD(T)$ which contain $v$ but not $v_1$ contribute
nothing to $\b_{i,d}(T)$, \ie
$$\b_{i,d}^{(2)}=\sum_{F \in \AD(T): \ |V(T) \setminus F|=d, v_1 \notin F , v \in F} \dim_k
\Hred_{i-2}(\Link_{\AD(T)} F;k)=0. $$
\end{thm}
\begin{proof}
Using Proposition (\ref{linkepsilon}) we write $\Link_{\AD(T)} F
=\e(e_1,\dots e_r;V)$ for $F \in \AD(T)$, where $e_1,\dots, e_r$
are the edges of $T$ which are disjoint from $F$ and $V= V(T)
\setminus F$. The faces of $\AD(T)$ in the sum $\b_{i,d}^{(2)}$
include $v$ but not $v_1$. Fix an $F \in \AD(T)$. The vertex $v$
does not belong to any of the edges $e_1,\dots, e_r$. Neither does
$v_1$ belong to any of $e_1,\dots, e_r$ since the only edge in $T$
which includes $v_1$ is $\{v,v_1\}$. Because $v_1$ is an element
of the set $V$ but not an element of $e_1 \cup \dots \cup e_r$ it
is in every maximal face of $\Link_{\AD(T)} F$. Hence
$\Link_{\AD(T)} F$ has zero reduced homology everywhere.
\end{proof}

We need now to calculate the contribution to $\b_{i,d}(T)$ from
the faces of $\AD(T)$ which contain neither $v$ nor $v_1$. First
we suppose that $F$ is such a face of $\AD(T)$ and that $v_n \in
F$ .

We can write $\Link_{\AD(T)}F= \e( \lc v,v_1 \rc, \lc v,v_{i_1}
\rc ,\dots, \lc v,v_{i_j} \rc, e_1, \dots, e_r ;V)$ where $1<i_1,
\dots < i_j <n$ and $j\geqslant0$ ( and $v_l \notin e_m$ for all
$l,m$).

\begin{thm}{\Lem}\label{TreeHomLem}
Suppose that $F \in \AD(T)$ with $v,v_1 \notin F$, $v_n \in F$.
Let $E=\Link_{\AD(T)} F= \e( \lc v,v_1 \rc, \lc v,v_{i_1} \rc
,\dots, \lc v,v_{i_j} \rc, e_1, \dots, e_r ;V)$. We have an
isomorphism of reduced homology modules
$$\Hred_i(E) \cong \Hred_{i-(j+1)}(\e(e_1,\dots,e_r;V'))$$ where $V'=V
\setminus \lc v,v_1,v_{i_1},\dots, v_{i_j} \rc.$
\end{thm}

\begin{proof}
We write the simplicial complex $E$ as the union of two simplicial
complexes which can be seen to have zero reduced homology
everywhere and then use the Mayer-Vietoris sequence. Let $\e_1=\e(
\lc v,v_1 \rc ;V)$ and $$\e_2 =\e(\lc v,v_{i_1} \rc ,\dots, \lc
v,v_{i_j} \rc, e_1, \dots, e_r ;V).$$ Clearly $E= \e_1 \cup \e_2.$
The reduced homology of $\e_1$ is zero as this is just a simplex.
The vertex $v_1$ belongs to every maximal face of $\e_2$ and
therefore $\e_2$ has zero reduced homology everywhere. We now use
the Mayer-Vietoris sequence

\begin{eqnarray*}
\dots \to \Hred_i(\e_1 \cap \e_2)\to \Hred_i(\e_1)\oplus
\Hred_i(\e_2) \to \Hred_i(E) \to \\ \to \Hred_{i-1}(\e_1 \cap
\e_2) \to \Hred_{i-1}(\e_1)  \oplus \Hred_{i-1}(\e_2) \to
\Hred_{i-1}(E) \to \dots
\end{eqnarray*}
which becomes
\begin{eqnarray*}
\dots \to \Hred_i(\e_1 \cap \e_2) \to 0 \to \Hred_i(E) \to
\Hred_{i-1}(\e_1 \cap \e_2) \to 0 \to \Hred_{i-1}(E) \to \dots
\end{eqnarray*}
to see that $\Hred_i(E) \cong \Hred_{i-1}(\e_1 \cap \e_2)$ for all
$i$. The intersection of the simplicial complexes $\e_1$ and
$\e_2$ can, by Lemma (\ref{lemB}), be written
$$\e_1 \cap \e_2 =\e( \lc v_{i_1} \rc, \dots ,\lc v_{i_j} \rc, e_1,
\dots, e_r; V \setminus \lc v,v_1 \rc). $$ None of the vertices
$v_{i_1}, \dots, v_{i_j}$ can belong to any of the edges
$e_1,\dots,e_r$ and so by Corollary  (\ref{CorE}) $$\Hred_i(\e_1
\cap \e_2) \cong \Hred_{i-j}(\e(e_1,\dots,e_r;V \setminus \lc
v,v_1,v_{i_1},\dots,v_{i_j} \rc )$$ for all $i$. Putting these
together we obtain $$\Hred_i(E)=\Hred_{i-1}(\e_1 \cap
\e_2)=\Hred_{i-(j+1)}(\e(e_1,\dots,e_r;V \setminus
\{v,v_1,v_{i_1},\dots,v_{i_j} \})).$$
\end{proof}

\begin{thm}{\Lem}\label{Lem3}
The contribution to $\b_{i,d}(T)$ from the faces of $\AD(T)$ which
include $v_n$ but exclude $v$ and $v_1$ may be expressed in terms
of $\mathbb {N}$-graded Betti numbers of $T''$ as follows.

\begin{eqnarray*}
\b_{i,d}^{(3)} &=&\sum_{F \in \AD(T): \ |V(T) \setminus F|=d,v,v_1
\notin
F,v_n \in F} \dim_k \Hred_{i-2}(\Link_{\AD(T)} F;k)\\
               &=&\sum_{j=0}^{n-2} {{n-2} \choose j}
               \b_{i-(j+1),d-(j+2)}(T'').
\end{eqnarray*}
\end{thm}

\begin{proof}
If $F$ is a face of the Alexander dual, $\AD(T)$, which includes
$v_n$ but neither $v$ nor $v_1$ then we may write
$$\Link_{\AD(T)} F= \e( \lc v,v_1 \rc, \lc v,v_{i_1} \rc ,\dots, \lc
v,v_{i_j} \rc, e_1, \dots, e_r ;V)$$ where
$\{v_{i_1},\dots,v_{i_j} \} \subseteq \{v_2,\dots,v_{n-1} \}$,
$e_1,\dots,e_r$ are the edges which do not feature any of
$v,v_1,\dots,v_n$ and $V=V(T) \setminus F$. By (\ref{TreeHomLem})
$$\Hred_i(\Link_{\AD(T)} F) \cong \Hred_{i-(j+1)}(\e(e_1,\dots,e_r;V'))$$
for all $i$ where $V'=V \setminus \lc v,v_1,v_{i_1},\dots,v_{i_j}
\rc $. We may consider $\e(e_1,\dots,e_r;V')$ to be the link of a
face of $\AD(T \setminus  \lc v,v_1,v_2,\dots,v_{n} \rc )=
\AD(T'')$, \ie $$\e(e_1,\dots,e_r;V')={\Link}_{\AD(T'')} F'$$
where $F'=F \cap \AD(T'').$ For a fixed choice of
$v_{i_1},\dots,v_{i_j}$ we consider all the possible choices for
$F$. The edges $e_1,\dots,e_r$ can be any selection of the edges
of $T \setminus \{v,v_1,\dots,v_{n-1},v_n \}=T''$ such that
$|\bigcup_{i=1}^{r}e_i|=d-(j+2)$. Therefore we have a contribution
to $\b_{i,d}^{(3)}$ of $\b_{i-(j+1),d-(j+2)}(T'')$. For a fixed
$j$ we have ${{n-2} \choose j}$ choices for the vertices
$v_{i_1},\dots,v_{i_j}$. Summing over the values of $j$ from 0 to
$n-2$ we have $\b_{i,d}^{(3)}=\sum_{j=0}^{n-2} {{n-2} \choose j}
\b_{i-(j+1),d-(j+2)}(T'')$.

\end{proof}

The final case to consider is when $F \in {\AD(T)}$ is such that
$v,v_1,v_n \notin F$. For such a face of the Alexander Dual we may
write
\begin{eqnarray*}
E:&=&\Link_{\AD(T)} F \\
  &=& \e( \lc v,v_1 \rc, \lc v,v_n \rc, \lc v,v_{i_1}
\rc,\dots,\lc v,v_{i_j} \rc,  \\
 & & \ \ \ \ \ \lc v_n,w_{k_1} \rc,\dots,\lc v_n,w_{k_l} \rc,e_1,\dots,e_t
;V)
\end{eqnarray*}
for some choices $2 \leqslant i_1 < i_2< \dots <i_j \leqslant n-1$
and $\{w_{k_1},\dots,w_{k_l} \} \subseteq \{w_1,\dots,w_m\}$ and
where $e_1,\dots,e_t$ are edges that do not feature any of
$\{v,v_1,\dots,v_n\}.$

\begin{thm}{\Lem}\label{anotherlem}
Let  $V=V(G) \setminus F$, and let
$$E=\e( \lc v,v_1 \rc, \lc v,v_n \rc, \lc v,v_{i_1} \rc,\dots,\lc
v,v_{i_j} \rc,\lc v_n,w_{k_1} \rc,\dots,\lc v_n,w_{k_l}
\rc,e_1,\dots,e_t ;V),$$ as above. We have an isomorphism of
reduced homology modules
$$\Hred_i(E)\cong \Hred_{i-(j+2)}(\e(e_1,\dots,e_t; V \setminus
\lc v,v_1,v_n,v_{i_1},\dots,v_{i_j} \rc )).$$
\end{thm}

\begin{proof}
Let $\e_1=\e(\{v,v_1\};V)$ and let $$\e_2=\e(\lc v,v_n \rc, \lc
v,v_{i_1} \rc,\dots,\lc v,v_{i_j} \rc,\lc v_n,w_{k_1}
\rc,\dots,\lc v_n,w_{k_l} \rc,e_1,\dots,e_t ;V).$$ The union of
$\e_1$ and $\e_2$ is $E$. The reduced homology of $\e_1$ is zero
everywhere as this is just a simplex. Also $\e_2$ has zero reduced
homology everywhere as $v_1$ is in every maximal face. From the
Mayer-Vietoris sequence we deduce that $\Hred_i(E)=
\Hred_{i-1}(\e_1 \cap \e_2)$ for all $i$. Now we examine $\e_1
\cap \e_2$:

\begin{eqnarray*}
\e_1 \cap \e_2 &=&
\e(\{v,v_1,v_n\},\{v,v_1,v_{i_1}\},\dots,\{v,v_1,v_{i_j}\},\{v,v_1,v_n,w_{k_1}\},\dots\\
               & &\dots,\{v,v_1,v_n,w_{k_l}\},\{v,v_1\}\cup e_1,\dots, \{v,v_1\} \cup e_t;V)\\
&=&\e(\{v_n\},\{v_{i_1}\},\dots,\{v_{i_j}\},\{v_n,w_{k_1}\},\dots,\{v_n,w_{k_l}\},e_1,\dots,e_t;V \setminus \{v,v_1\})\\
& &{\rm (Using\  Proposition\  (\ref{lemB}))}\\
               &=&\e( \{v_n \},\{v_{i_1}\},\dots,\{v_{i_j}
               \},e_1,\dots,e_t;V \setminus \{v,v_1\})\\
               & &{\rm (Using\  Proposition\  (\ref{lemC}))}.
\end{eqnarray*}
As $\{ v_n,v_{i_1},\dots,v_{i_j} \} \cap \bigcup_{i=1}^t e_i =
\emptyset$ we may use Corollary (\ref{CorE}) to conclude
$\Hred_i(E)=\Hred_{i-1}(\e_1 \cap
\e_2)=\Hred_{i-(j+2)}(\e(e_1,\dots,e_t;V'))$ where $$V'=V
\setminus \{v,v_1,v_n,v_{i_1},\dots,v_{i_j}\}.$$
\end{proof}

\begin{thm}{\Lem}\label{lemm99}
Let $E$ be as above. If $\lc w_{k_1},\dots,w_{k_l} \rc \nsubseteq
\bigcup_{i=1}^t e_i$ then $\Hred_{i}(E)=0$ for all $i$.
\end{thm}

\begin{proof}
Assume that $\lc w_{k_1},\dots,w_{k_l} \rc \nsubseteq
\bigcup_{i=1}^t e_i$. We may select $w \in \lc
w_{k_1},\dots,w_{k_l} \rc$ which is such that $w \notin
\bigcup_{i=1}^t e_i$. By Lemma (\ref{anotherlem})
$$\Hred_i(E) =\Hred_{i-(j+2)}(\e(e_1,\dots,e_t;V'))$$ for all $i$, where $V'=V \setminus
\{v,v_1,v_n,v_{i_1},\dots,v_{i_j}\}$. We now note that $w \in V'
\setminus \bigcup_{i=1}^t e_i$ so by lemma (\ref{lemA})
$\Hred_i(E)=0$ for all $i$.
\end{proof}

\begin{thm}[\rm]{\Def}
Let $W=\lc w_{k_1},\dots,w_{k_l} \rc$, a subset of $\Omega =\lc
w_1,\dots w_m \rc$ and  $\Gamma=\AD(T'' \setminus \overline{W})$,
where $\overline{W}= \Omega \setminus W$.
\end{thm}

By Lemma (\ref{lemm99}) we may assume that for $F\in \AD(T)$ such
that $v,v_1,v_n \notin F$ with
\begin{eqnarray*}
\Link_{\AD(T)} F&=&\e( \lc v,v_1 \rc, \lc v,v_n \rc, \lc v,v_{i_1}
\rc,\dots\\
& &\dots,\lc v,v_{i_j} \rc,\lc v_n,w_{k_1} \rc,\dots,\lc
v_n,w_{k_l} \rc,e_1,\dots,e_t ;V),
\end{eqnarray*}
with $V=V(T) \setminus F$, then $\lc w_{k_1},\dots,w_{k_l} \rc
\subseteq \bigcup_{i=1}^t e_i$, where $e_1,\dots,e_t$ are edges
that do not feature any of $v,v_1,v_2,\dots,v_n$, for otherwise
$F$ will give no contribution to $\b_{i,d}(T)$. By lemma
(\ref{anotherlem}) $$\Hred_i(\Link_{\AD(T)}
F)=\Hred_{i-(j+2)}(\e(e_1,\dots,e_t; V \setminus \lc
v,v_1,v_n,v_{i_1},\dots,v_{i_j} \rc )).$$

We may consider $\e(e_1,\dots,e_t;V')$, where $V'=V \setminus \lc
v,v_1,v_n,v_{i_1},\dots v_{i_j} \rc$, to be the link of a face,
$G$ say, of $\Gamma$. Hence we may write $\Hred_i(\Link_{\AD(T)}
 F)=\Hred_{i-(j+2)}(\Link_{\G}  G).$ Note that
$w_{k_1},\dots,w_{k_l} \notin G$ and $|\G \setminus G|=d-(j+3).$

\begin{thm}{\Lem}\label{lem4}
Let $Y(W,j)$ denote the contribution to $\b_{i,d}$ from the faces
of $\AD(T)$ which contain $\overline{W}$, do not contain any of
$w_{k_1},\dots,w_{k_l}$ and contain exactly $(n-2)-j$ of the
vertices $v_2,\dots,v_{n-1}$. Let $x=x(j)=i-(j+2)$ and
$y=y(j)=d-(j+3)$.

\begin{eqnarray*}
Y(W,j)&=&\b_{x,y}(T''\setminus \overline{W})-\sum_{U\subseteq W,
|U|=1} \b_{x,y}(T''
\setminus(U \cup \overline{W})) \\
&+& \sum_{U\subseteq W, |U|=2} \b_{x,y}(T'' \setminus(U \cup
\overline{W}))\\
& -& \dots + (-1)^l \b_{x,y}(T'' \setminus \lc w_1,\dots,w_{m}
\rc).
\end{eqnarray*}

\end{thm}

\begin{proof}
As stated above any face, $F$ say, of $\AD(T)$ of this type which
makes a non-zero contribution to $\b_{i,d}$ has link which we may
write as $\linkad F=\e( \lc v,v_1 \rc, \lc v,v_n \rc, \lc
v,v_{i_1} \rc,\dots,\lc v,v_{i_j} \rc,\lc v_n,w_{k_1}
\rc,\dots,\lc v_n,w_{k_l} \rc,e_1,\dots,e_t ;V)$ where $V=V(G)
\setminus F$, $e_1,\dots,e_t$ are edges which do not feature any
of $v,v_1,\dots,v_n$ and $\{w_{k_1},\dots,w_{k_l} \} \subseteq
\bigcup_{i=1}^te_i$. Let $E_1=\e( \{v,v_1\};V)$ and let
$$E_2=\e( \lc v,v_n \rc, \lc v,v_{i_1} \rc,\dots,\lc v,v_{i_j}
\rc,\lc v_n,w_{k_1} \rc,\dots,\lc v_n,w_{k_l} \rc,e_1,\dots,e_t
;V).$$ Clearly $E=E_1 \cup E_2$. Because $E_1$ is a simplex we
have $\Hred_i(E_1)=0$ for all $i$ and because $v_1 \in V \setminus
(\{v,v_n,v_{i_1},\dots,v_{i_j},w_{k_1},\dots,w_{k_l} \} \cup
\bigcup_{i=1}^t e_i)$ Lemma (\ref{lemA}) implies $\Hred_i(E_2)=0$
for all $i$. From the Mayer-Vietoris sequence we obtain the
isomorphism of reduced homology modules $\Hred_i(E) \cong
\Hred_{i-1}(E_1 \cap E_2)$. Using Lemmas (\ref{lemB}) and
(\ref{lemC}) we may write
\begin{eqnarray*}
E_1\cap
E_2&=&\e(\{v,v_1,v_n\},\{v,v_1,v_{i_1}\},\dots,\{v,v_1,v_{i_j}\},\{v,v_1,v_n,w_{k_1}\},\dots
\\
   &\dots,&
   \{v,v_1,v_n,w_{k_l}\},\{v,v_1\}\cup e_1,\dots,\{v,v_1\}\cup
   e_t;V)\\
   &=&\e( \{v_n\},\{v_{i_1}\},\dots,\{v_{i_j}\},e_1,\dots,e_t;V
   \setminus \{v,v_1\}).
\end{eqnarray*}
We now use Corollary (\ref{CorE}) to see that $$\Hred_i(\e(
\{v_n\},\{v_{i_1}\},\dots,\{v_{i_j}\},e_1,\dots,e_t;V
   \setminus \{v,v_1\})\cong
   \Hred_{i-(j+1)}(\e(e_1,\dots,e_t;V')),$$ where $V'=V \setminus
   \{v,v_1,v_n,v_{i_1},\dots,v_{i_j} \}$. Hence we obtain
   $$\Hred_i(E) \cong \Hred_{i-(j+2)}(\e(e_1,\dots,e_t;V'))$$ for
   all $i$.
We consider the contribution to $\b_{i,d}$ from all faces of $\AD$
of this form for a fixed $W$ and a fixed $j$. This will be
$\sum_{H} {\dim}_k \Hred_{i-2} (\e(H))$ where the sum is over all
induced subgraphs, $H$,  of $T'' \setminus \overline{W}$ which
feature all of the vertices $\{w_{k_1},\dots,w_{k_l} \}$. We use
an inclusion / exclusion principle to obtain the conclusion of the
statement of the lemma. First we take $\b_{x,y}(T'' \setminus
\overline{W})$ which also includes the contribution from the
subgraphs of $T'' \setminus \overline{W}$ which do not contain all
of the vertices of $W=\{w_{k_1},\dots,w_{k_l} \}$. We then
subtract all the Betti numbers $\b_{x,y}(T'' \setminus (U \cup
\overline{W}))$ for subsets $U$, of $W$ which contain a single
element, add those for which $|U|=2$ and so on to conclude
$$Y(W,j)=\b_{x,y}(T''\setminus \overline{W})-\sum_{U\subseteq W,
|U|=1} \b_{x,y}(T'' \setminus(U \cup \overline{W}))$$\\
$$+ \sum_{U\subseteq W, |U|=2} \b_{x,y}(T'' \setminus(U \cup
\overline{W}))  - \dots + (-1)^l \b_{x,y}(T'' \setminus \lc
w_1,\dots,w_{m} \rc).$$
\end{proof}

\begin{thm}{\Prop}\label{aprop} The contribution to the i{\rm th} Betti number of
degree d from faces of $\AD(T)$ which do not feature $v,v_1$ or
$v_n$ is
$$\b_{i,d}^{(4)}=\sum_{ W \subseteq \Omega}
\sum_{j=0}^{n-2}{ {n-2} \choose j}Y(W,j).$$
\end{thm}

\begin{proof}
From definition (\ref{bet1234}) $$\b_{i,d}^{(4)}= \sum_{F \in
\AD(T): \ |V(T) \setminus F|=d,\ v,v_1,v_n \notin F} \dim_k
\Hred_{i-2}(\Link_{\AD(T)} F;k).$$ For a fixed $W$ and a fixed $j$
we get a contribution of ${{n-2} \choose j} Y(W,j)$. We sum this
over all $W$ such that $ W \subseteq \Omega$ and all $j$ from $0$
to $n-2$ to obtain the result.
\end{proof}

\section{Betti Numbers in terms of Betti \\ Numbers of Subgraphs}

In this section we use the above results to write the graded Betti
numbers of a forest in terms of Betti numbers of some of its
subgraphs. We show how we can greatly simplify our initial
formula, which is given in the following proposition. Let the
notation be as in sections 9.1 and 9.2.

\begin{thm}{\Prop}\label{treeform1}
Summing all the cases, we obtain
\begin{eqnarray*}
\b_{i,d}(T) &=& \b_{i,d}(T') + \sum_{j=0}^{n-2} {{n-2} \choose j}
\b_{i-(j+1),d-(j+2)}(T'') \\
& &+\sum_{W \subseteq \Omega} \sum_{j=0}^{n-2} {{n-2} \choose j}
Y(W,j).
\end{eqnarray*}
\end{thm}

\begin{proof}
We use Lemmas (\ref{lem1}), (\ref{lem2}), (\ref{Lem3}) and
(\ref{aprop}) to obtain this result.
\end{proof}

The above formula may be simplified by considering the
contribution to $\b_{i,d}$ from $\b_{i-(j+2),d-(j+3)}(T''\setminus
\overline{X})$ for each $X \subseteq \Omega$.

\begin{thm}{\Prop}\label{Yzero}
The contribution to $\b_{i,d}$ from
$\b_{i-(j+2),d-(j+3)}(T''\setminus \overline{X})$ is zero for each
$X \subset \Omega$, with $X \neq \Omega$.
\end{thm}

\begin{proof}
We let $\Lambda=\b_{x,y}(T'' \setminus \overline{X})$ where
$x=x(j)=i-(j+2)$ and $y=y(j)=d-(j+3)$ and consider how many times
the formula given in (\ref{treeform1}) counts $\Lambda$. These all
come from the summand $\sum_{W \subseteq \Omega} \sum_{j=0}^{n-2}
{ {n-2} \choose j} Y(W,j)$. We now consider how many contributions
of $\Lambda$ we acquire from $Y(W,j)$ for various $W$. From Lemma
(\ref{lem4})
\begin{eqnarray*}
Y(W,j)&=&\b_{x,y}(T''\setminus
\overline{W})-\sum_{U\subseteq W, |U|=1} \b_{x,y}(T''
\setminus(U \cup \overline{W})) \\
&+& \sum_{U\subseteq W, |U|=2} \b_{x,y}(T'' \setminus(U \cup
\overline{W}))\\
& -& \dots + (-1)^l \b_{x,y}(T'' \setminus \lc w_1,\dots,w_{m}
\rc).
\end{eqnarray*}
So we need only consider $W$ such that $\overline{W} \subseteq
\overline{X}$. Suppose $\overline{X}=\{w_{a_1},\dots,w_{a_p}\}$
and for a fixed such $W$ let $|\overline{W}|=y$. We obtain
contributions of $\Lambda$ from the summand of $Y(W,j)$ which sums
to $(-1)^{p-y}\sum_{U \subseteq W, |U|=p-y}\b_{x,y}(T''
\setminus(U \cup \overline{W})$ whenever $U \cup \overline{W} =
\overline{X}$. Hence we count how many such subsets, $U$, there
are. We must have $|U|=p-y$, so we have ${ p \choose {p-y}}={ p
\choose {y}}$ choices. Therefore we obtain $(-1)^{p-y}{ p \choose
{y}}{ {n-2} \choose {j}}\Lambda$ for $y=0,1,\dots,p$. We sum these
over these values of $y$ to find
$$\Lambda { {n-2} \choose j} \left\{(-1)^p{ p \choose 0}+(-1)^{p-1}{ p \choose 1}+ \dots -{ p \choose{p-1}} + { p \choose p}
\right\}.$$ We complete the proof by noting that
$$(-1)^p{ p \choose 0}+(-1)^{p-1}{ p \choose 1}+ \dots -{ p \choose{p-1}} + { p \choose p}
=0.$$
\end{proof}

\begin{thm}{Theorem}\label{BettiTrees}
The i{\rm th} Betti number of degree d of the forest $T$ may be
expressed in terms of Betti numbers of the subgraphs $T'$ and
$T''$ as follows,
$$\b_{i,d}(T)=\b_{i,d}(T') + \sum_{j=0}^{n-1} {{n-1} \choose {j}}
\b_{i-(j+1),d-(j+2)}(T'').$$

\end{thm}

\begin{proof}
We use Lemma (\ref{Yzero}) to rewrite the formula of Lemma
(\ref{treeform1}) as
\begin{eqnarray*}
\b_{i,d}(T) &=& \b_{i,d}(T') + \sum_{j=0}^{n-2} {{n-2} \choose j}
\b_{i-(j+1),d-(j+2)}(T'') \\
            & & \ \ \ \ \ \ \ \ \ \ + \sum_{j=0}^{n-2} {{n-2} \choose j}
\b_{i-(j+1),d-(j+3)}(T''). \\
\end{eqnarray*}
We collect together the terms in Betti numbers of $T''$ to obtain
\begin{eqnarray*}
\sum_{j=0}^{n-2} {{n-2} \choose j} \b_{i-(j+1),d-(j+2)}(T'') +
\sum_{j=0}^{n-2} {{n-2} \choose j} \b_{i-(j+1),d-(j+3)}(T'') \\
={{n-2} \choose 0} \b_{i-1,d-2}(T'')+ \left\{{ {n-2} \choose 1}+{
{n-2}
\choose 0}\right\} \b_{i-2,d-3}(T'')+ \dots\\
 +\left\{{ {n-2} \choose {n-2}}+{
{n-2} \choose {n-3}} \right\}\b_{i-(n-1),d-n}(T'')+{ {n-2} \choose
{n-2}} \b_{i-n,d-(n+1)}(T'')\\
=\sum_{j=1}^n\left\{ { {n-2} \choose {j-1}}+ { {n-2} \choose
{j-2}}
\right\} \b_{i-j,d-(j+1)}(T'')\\
=\sum_{j=0}^{n-1} { {n-1} \choose {j}}
\b_{i-(j+1),d-(j+2)}(T'').\\
\end{eqnarray*}

\end{proof}

\section{Projective Dimension of Forests}

We can use our results about the Betti numbers of forests to say
something about the projective dimension. The projective dimension
of $T$ can be found from the projective dimensions of the
subforests $T'$ and $T''$.

\begin{thm}[\rm]{Remark}
Note that the projective dimension of any forest is independent of
the characteristic of the field $k$. This follows from the fact
that the Betti numbers are independent of $k$.
\end{thm}

\begin{thm}{Theorem}\label{pdtrees}
Let $p=\pd(T)$,$p'=\pd(T')$ and $p''=\pd(T'')$. The projective
dimension of the forest $T$ is $$p= \max \lc p', p''+n \rc.$$
\end{thm}

\begin{proof}
To determine $p$ we must find the largest $i\in \mathbb{N}$ such
that $\b_i(T) \neq 0$. Using Theorem (\ref{BettiTrees}) we see
that
\begin{eqnarray*}
\b_i(T)&=&\sum_{d \in \mathbb{N}} (\b_{i,d}(T))\\
&=& \sum_{d \in \mathbb{N}}(\b_{i,d}(T')) + \sum_{d \in
\mathbb{N}}( \sum_{j=0}^{n-1} {{n-1} \choose j}
\b_{i-(j+1),d-(j+2)}(T'')).
\end{eqnarray*}
This implies that $\b_i(T)=0$ if and only if $\b_i(T')=0$ and
$$\b_{i-1}(T''),\b_{i-2}(T''),\dots,\b_{i-n}(T'')$$ are all zero.
If $\b_{i-n}(T'')=0$ then $\b_{i-n+1}(T''),\dots,\b_{i-1}(T'')$
are all zero too. Hence $\b_i(T)=0$ if and only if $\b_i(T')=0$
and $\b_{i-n}(T'')=0$. Therefore $i >p$ if and only if $i>p'$ and
$i-n>p''$. Hence $p= \max\{p',p''+n\}$.


%



\end{proof}

\begin{thm}[\rm]{Remark}
This enables us to find the projective dimension of a forest by
recursively applying the above theorem. We can find the projective
dimension, an algebraic invariant, by using an algorithm on a
forest, a combinatorial object. We can keep finding the projective
dimension of the forest in terms of the projective dimensions of
smaller forests until we arrive at a star graph, whose projective
dimension we know from Proposition (\ref{Star}).
\end{thm}

\begin{appendix}
\end{appendix}
\newpage
\end{document}